\documentclass[review,hidelinks,onefignum,onetabnum]{siamart220329}

\usepackage{lipsum}
\usepackage{amsfonts}
\usepackage{graphicx}
\usepackage{epstopdf}
\usepackage{algorithmic}
\ifpdf
  \DeclareGraphicsExtensions{.eps,.pdf,.png,.jpg}
\else
  \DeclareGraphicsExtensions{.eps}
\fi

\newsiamremark{remark}{Remark}
\newsiamremark{hypothesis}{Hypothesis}
\crefname{hypothesis}{Hypothesis}{Hypotheses}
\newsiamthm{claim}{Claim}

\headers{Least-squares neural network method}{Z. Cai, J. Choi, and M. Liu}

\title{Least-Squares Neural Network (LSNN) Method \\ for
Linear Advection-Reaction Equation: \\[1mm] Discontinuity Interface
\thanks{Submitted to the editors DATE.
\funding{This work was supported in part by the National Science Foundation under grant DMS-2110571.}}}

\author{Zhiqiang Cai\thanks{Department of Mathematics, Purdue University, 150 N. University Street, West Lafayette, IN 47907-2067 
  (\email{caiz@purdue.edu}, \email{choi508@purdue.edu}).}
\and Junpyo Choi\footnotemark[2]
\and Min Liu\thanks{School of Mechanical Engineering, Purdue University, 585 Purdue Mall,
West Lafayette, IN 47907-2088(\email{liu66@purdue.edu}). }}

\usepackage{amsopn}

\ifpdf
\hypersetup{
  pdftitle={Least-Squares Neural Network (LSNN) Method For
Linear Advection-Reaction Equation: Discontinuity Interface},
  pdfauthor={Z. Cai, J. Choi, and M. Liu}
}
\fi
\nolinenumbers
\usepackage[utf8]{inputenc}
\usepackage{bm}
\usepackage{mathtools}
\usepackage{amssymb}
\usepackage{indentfirst}
\usepackage{subfigure}
\usepackage{tcolorbox}
\usepackage{color}
\usepackage[export]{adjustbox}
\usepackage{scalerel}
\usepackage{tikz}
\usetikzlibrary{decorations.pathreplacing}

\usepackage{algorithmic}
\Crefname{ALC@unique}{Line}{Lines}

\newcommand{\R}{\mathbb{R}}

\usepackage{graphicx}
\usepackage{diagbox}
\usepackage{pifont}
\newcommand{\vertiii}[1]{{\left\vert\kern-0.25ex\left\vert\kern-0.25ex\left\vert #1 
    \right\vert\kern-0.25ex\right\vert\kern-0.25ex\right\vert}}

\newcommand{\bxi}{\mbox{\boldmath$\xi$}}

\usepackage{enumitem}
\setlist[itemize]{left=16pt} 

\def\bx{{\bf x}}

\def\cL{{\cal L}}
\def\cM{{\cal M}}

\def\cT{{\cal T}}

\begin{document}

\maketitle

\begin{abstract}
We studied the least-squares ReLU neural network (LSNN) method for solving linear advection-reaction equation with discontinuous solution in [Cai, Zhiqiang, Jingshuang Chen, and Min Liu. ``Least-squares ReLU neural network (LSNN) method for linear advection-reaction equation.'' Journal of Computational Physics 443 (2021), 110514]. The method is based on a least-squares formulation and uses a new class of approximating functions: ReLU neural network (NN) functions. A critical and additional component of the LSNN method, differing from other NN-based methods, is the introduction of a properly designed and physics preserved discrete differential operator.

In this paper, we study the LSNN method for problems with discontinuity interfaces. First, we show that ReLU NN functions with depth $\lceil \log_2(d+1)\rceil+1$ can approximate any $d$-dimensional step function on a discontinuity interface generated by a vector field as streamlines with any prescribed accuracy. By decomposing the solution into continuous and discontinuous parts, we prove theoretically that discretization error of the LSNN method using ReLU NN functions with depth $\lceil \log_2(d+1)\rceil+1$ is mainly determined by the continuous part of the solution provided that the solution jump is constant. Numerical results for both two- and three-dimensional test problems with various discontinuity interfaces show that the LSNN method with enough layers is accurate and does not exhibit the common Gibbs phenomena along discontinuity interfaces. 
\end{abstract}



\section{Introduction}
Let $\Omega$ be a bounded domain in ${\R}^d$ ($d\ge2$) 
with Lipschitz boundary $\partial \Omega$. Consider the linear advection-reaction equation 
\begin{equation}\label{pde}
    \left\{\begin{array}{rccl}
    u_{\bm\beta} + \gamma\, u &=&f, &\text{ in }\, \Omega, \\[2mm]
    u&=&g, &\text{ on }\, \Gamma_{-},
    \end{array}\right.
\end{equation}
where $\bm{\beta}(\bx) = (\beta_1, \cdots, \beta_d)^T \in C^0(\bar{\Omega})^d$ is a given advective velocity field, $u_{\bm\beta}=\bm{\beta}\cdot\nabla u$ denotes the directional derivative of $u$ along $\bm{\beta}$, and $\Gamma_{-}$ is the inflow part of the boundary $\Gamma=\partial \Omega$ given by
\begin{equation}
    \Gamma_- = \{\bx\in\Gamma :\, \bm{\beta}(\bx) \cdot \bm{n}(\bx) <0\}
\end{equation}
with $\bm{n}(\bx)$ being the unit outward normal vector to $\Gamma$ at $\bx\in \Gamma$. We assume that the reaction coefficient $\gamma \in C^0(\bar{\Omega})$, the source term $f \in L^2(\Omega)$, and $g \in L^2(\Gamma_-)$. 

When the inflow boundary data $g$ is discontinuous, so is the solution of \cref{pde} as the solution is entirely determined by the characteristic curves (see \cref{l:solution}), along which the discontinuities are propagated across the domain. The discontinuity interface may be determined by the characteristic curves emanating from where $g$ is discontinuous. By using the location of the interface, one may design an accurate mesh-based numerical method. However, this type of method is usually limited to linear problems and is difficult to be extended to nonlinear hyperbolic conservation laws. 

In \cite{Cai2021linear}, we studied the least-squares ReLU neural network method (LSNN) for solving \cref{pde} with discontinuous solution. The method is based on the $L^2(\Omega)$ norm least-squares formulation analyzed in \cite{de2004least, bochev2016least} and employs a new class of approximating functions: multilayer perceptrons with the rectified linear unit (ReLU) activation function, i.e., ReLU neural network (NN) functions. A critical and additional component of the LSNN method, differing from other NN-based methods, is the introduction of a properly designed discrete differential operator. 

One of the appealing features of the LSNN method is its ability of automatically approximating the discontinuous solution without using {\it a priori} knowledge of the location of the interface. Hence, the method is applicable to nonlinear problems (see \cite{Cai2021nonlinear, Cai2022nonlinear}). Compared to mesh-based numerical methods including various adaptive mesh refinement (AMR) algorithms that locate the discontinuity interface through local mesh refinement (see, e.g., \cite{dahmen2012adaptive,houston2000posteriori,liu2020adaptive}), the LSNN method, a meshfree and pointfree method, is much more effective in terms of the number of degrees of freedom. Theoretically, it was shown in \cite{Cai2021linear} that a two- or three-layer ReLU NN function in two dimensions is sufficient to well approximate the discontinuous solution of \cref{pde} without oscillation, provided that the interface consists of a straight line or two-line segments and that the solution jump along the interface is constant.

The assumption on at most two-line segments in \cite{Cai2021linear} is very restrictive even in two dimensions. In general, the discontinuous solution of \cref{pde} has interfaces that are hypersurfaces in $d$ dimensions. The purpose of this paper has two-fold. First,
we show that any step function with interface that is generated by a vector field as streamlines may be approximated by ReLU NNs with at most $\lceil \log_2(d+1)\rceil+1$ layers for achieving a given approximation accuracy $\varepsilon$ (see Lemma~4.3), which extends the approximation result in \cite{Cai2021linear}. This is done by constructing a continuous piecewise linear (CPWL) function with a sharp transition layer of $\varepsilon$ width and combining the main results in \cite{arora2016understanding,tarela1999region,wang2005generalization} (see Proposition 2.1). Question on approximating piecewise smooth functions by ReLU NNs also arises in data science applications such as classification, etc. Some convergence rates were obtained in \cite{petersen2018optimal, caragea2020neural, KOK19, IF18, IF22}. Particularly, for a given $C^{\beta}$ ($\beta >0$) interface, \cite{petersen2018optimal} established approximation rates of ReLU NNs with no more than $L(\beta,d)=(3+\lceil \log_2\beta\rceil)(11+2\beta/d)$ layers. This upper bound is not applicable to $C^0$ interface and becomes large for very smooth interface.




Second, we establish a new kind of {\it a priori} error estimates (see \cref{curve-theorem}) for the LSNN method in $d$ dimensions for discontinuity interface. To do so, we decompose the solution as the sum of the discontinuous and continuous parts (see \cref{decop}). The continuous part of the solution may be approximated well by (even shallow) ReLU NN functions with standard approximation property (see, e.g., \cite{DeVore1997,Petrushev1998,pinkus1999approximation, shen2019deep, DeVore2021}). The discontinuous part of the solution can be approximated accurately by the class of all ReLU NN functions from $\mathbb{R}^d$ to $\mathbb{R}$ with at most $\lceil \log_2(d+1)\rceil+1$ depth, provided that the solution jump is constant. Hence, the accuracy of the LSNN method is mainly determined by the continuous part of the solution.



The explicit construction in this paper indicates that a ReLU NN function with at most $\lceil \log_2(d+1)\rceil+1$ depth is sufficient to accurately approximate discontinuous solutions without oscillation. The necessary depth $\lceil \log_2(d+1)\rceil+1$ of a ReLU NN function is shown numerically through several test problems in both two and three dimensions (two-hidden layers for $d=2,3$). At the current stage, it is still very expensive to numerically solve the discrete least-squares minimization problem, which is high dimensional and non-convex, when using stochastic gradient descent, Adam \cite{kingma2015}, etc., even though the degrees of freedom of the LSNN method is much less than those of mesh-based numerical methods.

Followed by recent success of deep neural networks (DNNs) in machine learning and artificial intelligence tasks such as computer vision and pattern recognition, there have been active interests in using DNNs for solving partial differential equations (PDEs) (see, e.g., \cite{PNAS2019, Berg18, cai2020deep, Weinan18, raissi2019physics, Sirignano18}). Due to the fact that the collection of DNN functions is not a linear space, NN-based methods for solving PDEs may be categorized as the Ritz and least-squares (LS) methods. The former (see, e.g., \cite{Weinan18}) requires the underlying problem having a natural minimization principle and hence is not applicable to \cref{pde}. 

For a given PDE, there are many least-squares methods and their efficacy depends on norms used for the PDE and for its boundary and/or initial conditions. When using NNs as approximating functions, least-squares methods may be traced back at least to 1990s (see, e.g., \cite{Dissanayake94,LLF98}), where the discrete $L^2$ norm on a uniform integration mesh was employed for both PDEs of the strong form and their boundary/initial conditions. Along this line, it is the popular physics-informed neural networks (PINNs) by Raissi-Perdikaris-Karniadakis \cite{raissi2019physics} in 2019 which uses auto-differentiation for computing the underlying differential operator at each integration point. Since the solution of \cref{pde} is discontinuous, those NN-based least-squares methods are also not applicable.


The rest of the paper is organized as follows. In \Cref{relu dnn}, we describe ReLU NN functions and CPWL functions, and introduce a known result about their relationship. Then we further investigate the structure of ReLU NN functions. \Cref{lsnn method} reviews the LSNN method in \cite{Cai2021linear} and formulates the method based on the framework in \Cref{relu dnn}. Then we prove that the method is capable of locating any discontinuity interfaces of the problem in \Cref{error estimates}. Finally, \Cref{numerics} presents numerical results for both two- and three-dimensional test problems with various discontinuity interfaces.

\section{ReLU NN functions}\label{relu dnn}
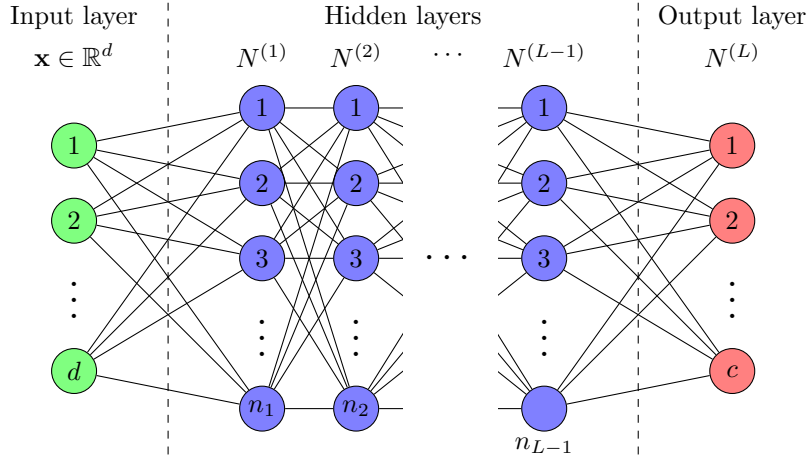
\begin{figure}[htbp]
\centering
\begin{tikzpicture}[]
\def\layersep{2.5cm}
\tikzset{%
  neuron missing/.style={
    draw=none, 
    scale=1.5,
    execute at begin node=$\vdots$
  },
}
\tikzstyle{neuron}=[circle,draw,minimum size=17pt]
\tikzstyle{input neuron}=[neuron, fill=green!50];
\tikzstyle{output neuron}=[neuron, fill=red!50];
\tikzstyle{hidden neuron}=[neuron, fill=blue!50];

\draw [dashed] (1/2*\layersep,0.9) -- (1/2*\layersep,-5.2);
\draw [dashed] (3*\layersep,0.9) -- (3*\layersep,-5.2);

\foreach \name / \y in {1,2,4}
\node[input neuron] (I-\name) at (0,-\y) {};
\node[neuron missing] (I-3) at (0,-2.9) {};

\foreach \name / \y in {1,2,3,5}
\path[yshift=0.5cm]
node[hidden neuron] (H1-\name) at (\layersep,-\y cm) {};
\path[yshift=0.5cm]
node[neuron missing] (H1-4) at (\layersep,-3.9cm) {};

\foreach \name / \y in {1,2,3,5}
\path[yshift=0.5cm]
node[hidden neuron] (H2-\name) at (3/2*\layersep,-\y cm) {};
\path[yshift=0.5cm]
node[neuron missing] (H2-4) at (3/2*\layersep,-3.9cm) {};
            
\foreach \name / \y in {1,2,3,5}
\path[yshift=0.5cm]
node[hidden neuron] (H3-\name) at (5/2*\layersep,-\y cm) {};
\path[yshift=0.5cm]
node[neuron missing] (H3-4) at (5/2*\layersep,-3.9cm) {};


\foreach \name / \y in {1,2,4}
\node[output neuron] (O-\name) at (7/2*\layersep,-\y) {};
\node[neuron missing] (O-3) at (7/2*\layersep,-2.9) {};

\foreach \source in {1,2,4}
\foreach \dest in {1,2,3,5}
\path (I-\source) edge (H1-\dest);

\foreach \source in {1,2,3,5}
\foreach \dest in {1,2,3,5}
\path (H1-\source) edge (H2-\dest);

\foreach \source in {1,2,3,5}
\foreach \dest in {1,2,3,5}
\path (H2-\source) edge (H3-\dest);
            
\foreach \source in {1,2,3,5}
\foreach \dest in {1,2,4}
\path (H3-\source) edge (O-\dest);
    
\fill [white] (3/2*\layersep+1/4*\layersep,0cm) rectangle (3/2*\layersep+3/4*\layersep,-5cm);
        
\path[yshift=0.5cm]
node[scale=1.5] at (2*\layersep,-3cm) {$\cdots$};

\node[] at (0cm,0.7cm) {Input layer};
\node[] at (0cm,0.2cm) {$\mathbf{x}\in \mathbb{R}^d$};

\node[] at (0,-1) {$1$};
\node[] at (0,-2) {$2$};
\node[] at (0,-4) {$d$};

\node[] at (\layersep,-0.5cm) {$1$};
\node[] at (\layersep,-1.5cm) {$2$};
\node[] at (\layersep,-2.5) {$3$};
\node[] at (\layersep,-4.5) {$n_1$};

\node[] at (3/2*\layersep,-0.5cm) {$1$};
\node[] at (3/2*\layersep,-1.5cm) {$2$};
\node[] at (3/2*\layersep,-2.5) {$3$};
\node[] at (3/2*\layersep,-4.5) {$n_2$};

\node[] at (5/2*\layersep,-0.5cm) {$1$};
\node[] at (5/2*\layersep,-1.5cm) {$2$};
\node[] at (5/2*\layersep,-2.5) {$3$};
\node[] at (5/2*\layersep,-5) {$n_{L-1}$};

\node[] at (7/2*\layersep,-1cm) {$1$};
\node[] at (7/2*\layersep,-2cm) {$2$};
\node[] at (7/2*\layersep,-4) {$c$};

\node[] at (\layersep+3/4*\layersep,0.7cm) {Hidden layers};
\node[] at (\layersep,0.2cm) {$N^{(1)}$};
\node[] at (3/2*\layersep,0.2cm) {$N^{(2)}$};
\node[] at (2*\layersep,0.2cm) {$\cdots$};
\node[] at (5/2*\layersep,0.2cm) {$N^{(L-1)}$};
\node[] at (7/2*\layersep,0.7cm) {Output layer};
\node[] at (7/2*\layersep,0.2cm) {$N^{(L)}$};

\end{tikzpicture}
\caption{The neural network function structure}
\label{relu dnn figure}
\end{figure}

First we begin with the definition of the rectified linear unit (ReLU) activation function. The ReLU activation function $\sigma$ is defined by
\begin{equation*}
\sigma(t) = \max\{0,t\}
=
\left\{\begin{array}{rl}
 0, & \text{if }  t\leq 0,\\[2mm]
 t, &  \text{otherwise.}
 \end{array}
 \right.
\end{equation*}
We say that a function $\mathcal{N}:\mathbb{R}^d\to\mathbb{R}^c$ with $c,d\in\mathbb{N}$ is a ReLU neural network (NN) function if the function $\mathcal{N}$ has a representation:
\begin{equation}\label{relu dnn def}
     \mathcal{N}=N^{(L)} \circ \cdots\circ N^{(2)}\circ N^{(1)}\text{ with }L>1,
    \end{equation}
where the symbol $\circ$ denotes the composition of functions, and for each $l=1,\ldots,L$, $N^{(l)}:\mathbb{R}^{n_{l-1}}\to\mathbb{R}^{n_{l}}$ with $n_l,n_{l-1}\in\mathbb{N}$ ($n_0=d$, $n_L=c$) given by:
\begin{enumerate}
    \item For $l=L$, $N^{(L)}(\mathbf{x})=\bm{\omega}^{(L)}\mathbf{x}-\mathbf{b}^{(L)}$ for all $\mathbf{x}\in\mathbb{R}^{n_{L-1}}$ for $\bm{\omega}^{(L)}\in\mathbb{R}^{n_L\times n_{L-1}}$, $\mathbf{b}^{(L)}\in\mathbb{R}^{n_L}$.
    \item For each $l=1,\ldots,L-1$, $N^{(l)}\left(\mathbf{x}\right)= \sigma  \left(\bm{\omega}^{(l)}\mathbf{x}-\mathbf{b}^{(l)}\right)$ for all $\mathbf{x}\in\mathbb{R}^{n_{l-1}}$ for $\bm{\omega}^{(l)}\in\mathbb{R}^{n_l\times n_{l-1}}$, $\mathbf{b}^{(l)}\in\mathbb{R}^{n_l}$, where $\sigma$ is applied to each component.
\end{enumerate}
We now establish some terminology as follows.
Let a ReLU NN function $\mathcal{N}$ have a representation $N^{(L)} \circ \cdots\circ N^{(2)}\circ N^{(1)}$ with $L>1$ (not unique) as in \cref{relu dnn def}. Then we say: 
\begin{enumerate}
\item $N^{(l)}$ is the $l^{\text{th}}$ layer (or also the $l^{\text{th}}$ hidden layer when $l<L$) of the representation, and the representation has $L$ layers or depth $L$, and $L-1$ hidden layers.
    \item The entries of $\bm{\omega}^{(l)}$ and $\mathbf{b}^{(l)}$ are the weights and biases, respectively, of the $l^{\text{th}}$ layer (or also the $l^{\text{th}}$ hidden layer when $l<L$).
    \item The natural number $n_{l}$ is the width or the number of neurons of the $l^{\text{th}}$ layer (or also the $l^{\text{th}}$ hidden layer when $l<L$).
\end{enumerate}
A motivation for this terminology is illustrated in \cref{relu dnn figure}.

For a given positive integer $n$, denote the set of all ReLU NN functions from $\mathbb{R}^d$ to $\mathbb{R}$ that have representations with depth $L$ and the total number of neurons of the hidden layers $n$ by
\[
\mathcal{M}(L,n)=\left\{\mathcal{N}:\mathbb{R}^d\to\mathbb{R} :\mathcal{N}=N^{(L)} \circ \cdots\circ N^{(2)}\circ N^{(1)} \text{ defined in \emph{\eqref{relu dnn def}}} :  n=\sum_{l=1}^{L-1} n_l
\right\}.
\]
Denote the set of all ReLU NN functions from $\mathbb{R}^d$ to $\mathbb{R}$ with $L$-layer representations by $\mathcal{M}(L)$. Then
\begin{equation}\label{relu dnn union}
  \mathcal{M}(L)=\bigcup_{n\in\mathbb{N}}\mathcal{M}(L,n).  
\end{equation}

We now introduce another function class, which is the set of continuous piecewise linear functions, and explore a theorem about the relationship between the two function classes. We say that a function $f\colon\mathbb{R}^d\to \mathbb{R}$ with $d\in\mathbb{N}$ is continuous piecewise linear (CPWL) if there exists a finite set of polyhedra with nonempty interior such that:
\begin{enumerate}
    \item The interiors of any two polyhedra in the set are disjoint.
    \item The union of the set is $\mathbb{R}^d$.
    \item $f$ is affine linear on each polyhedron in the set, i.e., on each polyhedron in the set, $f(\mathbf{x})=\mathbf{a}^T\mathbf{x}+b$ for all $\mathbf{x}\in\mathbb{R}^d$ for $\mathbf{a}\in\mathbb{R}^d$, $b\in\mathbb{R}$.
\end{enumerate}
Here by a polyhedron, we mean a subset of $\mathbb{R}^d$ surrounded by a finite number of hyperplanes, i.e.,  the solution set of a system of linear
inequalities
\begin{equation}\label{polyhedron}
    \{\mathbf{x}\in\mathbb{R}^d:\mathbf{A}\mathbf{x}\le \mathbf{b}\}\text{ for } \mathbf{A}\in\mathbb{R}^{m\times d}, \mathbf{b}\in\mathbb{R}^m\text{ with }m\in\mathbb{N},
\end{equation}
where the inequality is applied to each component. Thus the interior of the polyhedron in \cref{polyhedron} is
\begin{equation*}
    \{\mathbf{x}\in\mathbb{R}^d:\mathbf{A}\mathbf{x}< \mathbf{b}\}.
\end{equation*}

\begin{proposition}\label{cpwl=relu}
The set of all {\em CPWL} functions $f\colon\mathbb{R}^d\to\mathbb{R}$ is equal to $\cM(\lceil \log_2(d+1)\rceil+1)$, i.e., the set of all {\em ReLU NN functions} from $\mathbb{R}^d$ to $\mathbb{R}$ that have representations with depth $\lceil \log_2(d+1)\rceil+1$.
\end{proposition}

\begin{proof}
$\cM(\lceil \log_2(d+1)\rceil+1)$ is clearly a subset of the set of CPWL functions. Conversely, it is proved in \cite{arora2016understanding} that every CPWL function is a ReLU NN function from $\mathbb{R}^d$ to $\mathbb{R}$ that has a representation with depth at most $\lceil \log_2(d+1)\rceil+1$. Now, the result follows from the fact that $\cM(L)\subset \cM(\lceil \log_2(d+1)\rceil+1)$ for any $L\leq \lceil \log_2(d+1)\rceil+1$.
\end{proof}

\Cref{cpwl=relu} enables us to employ ReLU NN functions with a few layer representations to problems where CPWL functions are used, and to only control the number of neurons. Except the case $d=1$ (see, e.g., \cite{arora2016understanding}), there are currently no known results to give tight bounds on the number of neurons of the hidden layers. Therefore we suggest the following approach. The following proposition is a trivial fact.

\begin{proposition}\label{subset}
$\cM(L,n)\subseteq\cM(L,n+1)$.
\end{proposition}
Now \cref{cpwl=relu,relu dnn union,subset} suggest how we control the number of neurons of the hidden layers, i.e., when approximating a function $\mathbb{R}^d\to\mathbb{R}$ by a CPWL function, we start with the class $\mathcal{M}(\lceil \log_2(d+1)\rceil+1,n)$ with a small $n$ and the same width of each hidden layer, and then increase $n$ to have a better approximation.

Finally, in \cref{12 comparison figure,2d test1,2d test2,2d test3,2d test4,3d test1,3d test2,3d test3}, by the $l^{\text{th}}$- (hidden) layer breaking hyperplanes of a given representation as in \cref{relu dnn def} (with the output dimension being 1), we shall mean the set
\begin{itemize}
    \item $\{\mathbf{x}\in\Omega:\bm{\omega}^{(1)}\mathbf{x}-\mathbf{b}^{(1)}\text{ has a zero component}\}$ when $l=1$,
    \item $\{\mathbf{x}\in\Omega:\bm{\omega}^{(l)}(N^{(l-1)}\circ\cdots\circ N^{(2)}\circ N^{(1)}(\mathbf{x}))-\mathbf{b}^{(l)}\text{ has a zero component}\}$ when $2\le l<L$.
\end{itemize}
Breaking hyperplanes give a partition of $\mathbb{R}^d$, and on each element in the partition, the ReLU NN function is affine linear. (Breaking hyperplanes correspond to boundaries of linear regions of NNs as introduced in \cite{montufar2014number,pascanu2013number}.) They will be presented to help our understanding of the graphs of ReLU NN function approximations, especially along discontinuity interfaces.

\section{The LSNN method}\label{lsnn method}

We define the least-squares (LS) functional as
 \begin{equation}\label{ls}
    \mathcal{L}(v;{\bf f}) = \|v_{\bm\beta} +\gamma\, v-f\|_{0,\Omega}^2 +  \|v-g\|_{-\bm\beta}^2, 
\end{equation}
where ${\bf f} = (f,g)$, and $\|\cdot\|_{0,\Omega}$ and $\|\cdot\|_{-\bm\beta}$ denote, respectively, the $L^2(\Omega)$ norm and the weighted $L^2(\Gamma_{-})$ norm over the inflow boundary given by
 \[
 \|v\|_{-\bm{\beta}} 
 =\left<v,v\right>^{1/2}_{-\bm{\beta}} 
 =\left( \int_{\Gamma_-} |\bm{\beta}\! \cdot \!\bm{n}|\, v^2\,ds\right)^{1/2}.
 \]
Let $V_{\bm\beta} = \{v\in L^2(\Omega): v_{\bm{\beta}}\in L^2(\Omega)\}$ that is equipped with the norm as
 \[
 \vertiii{v}_{\bm\beta}= \left(\|v\|_{0,\Omega}^2 + \|v_{\bm\beta}
 \|_{0,\Omega}^2 \right)^{1/2}.
 \]
The least-squares formulation of problem \cref{pde} is to seek $u\in V_{\bm\beta}$ such that
\begin{equation}\label{minimization1}
    \mathcal{L}(u;{\bf f}) = \min_{v\in V_{\bm\beta}} \mathcal{L}(v;{\bf f}).
\end{equation}

\begin{proposition}[see \cite{bochev2016least, de2004least}]
Assume that either $\gamma=0$ or there exists a positive constant $\gamma_0$ such that
\begin{equation}\label{gamma}
    \gamma (\bx) -\frac{1}{2}\nabla \cdot \bm{\beta} (\bx) \geq \gamma_0 >0\quad  \text{ for almost all } \bx\in \Omega.
\end{equation}
Then the homogeneous {\em LS} functional $\mathcal{L}(v;{\bf 0})$ is equivalent to the norm $\vertiii{v}_{\bm\beta}^2$, i.e., there exist positive constants $\alpha$ and $M$ such that 
\begin{equation}\label{equiv}
\alpha\, \vertiii{v}_{\bm\beta}^2 
\leq \mathcal{L}(v;{\bf 0}) \leq M\, \vertiii{v}_{\bm\beta}^2 \quad\mbox{for all}\,\, v\in V_{\bm\beta}.
\end{equation}
\end{proposition}

The norm equivalence \cref{equiv} implies that problem \cref{minimization1} is well posed.
\begin{proposition}[see \cite{bochev2016least, de2004least}]
Problem \cref{minimization1} has a unique solution $u\in V_{\bm\beta}$ satisfying the following {\it a priori} estimate
\begin{equation}\label{stability}
\vertiii{u}_{\bm\beta} \leq C\, \left(\|f\|_{0,\Omega} + \|g\|_{-\bm\beta}\right).
\end{equation}
\end{proposition}

We note that 
$\cM(L,n)$ is a subset of $V_{\bm\beta}$. The least-squares approximation is then to find \\ $u_{_N} \in \cM(L,n)$ such that 
\begin{equation}\label{L-NN}
     \mathcal{L}\big(u_{_N};{\bf f}\big)
     = \min\limits_{v\in \cM(L,n)} \mathcal{L}\big(v;{\bf f}\big).
\end{equation}

\begin{lemma}[see \cite{Cai2021linear}]\label{Cea}
Let $u$ and $u_{_N}$ be the solutions of problems \cref{minimization1} and \cref{L-NN}, respectively. Then we have
\begin{equation}\label{Cea-L}
    \vertiii{u-u_{_N}}_{\bm\beta}
    \le \left(\dfrac{M}{\alpha}\right)^{1/2} \inf_{v\in \cM(L,n)} \vertiii{u-v}_{\bm\beta},
\end{equation}
where $\alpha$ and $M$ are constants in \cref{equiv}.
\end{lemma}

Optimization methods for solving the LSNN discretization problem in \cref{L-NN} repeatedly compute the following integration
\begin{equation}\label{int}
    \int_\Omega \left(v_{\bm\beta} +\gamma\, v-f\right)^2(\bx)\,d\bx
\end{equation}
for function $v$ in $\cM(L,n)$. In practice, the integration in \cref{int} is approximated by a numerical integration. Unlike conventional numerical methods using fixed meshes, designing an efficient and accurate numerical integration for \cref{int} is a non-trivial task. Apparently, the commonly used numerical integration of Monte Carlo type in scientific machine learning is inaccurate for problems with local features. It is also obvious that an accurate numerical integration should be based on the integral of the exact solution $u$, i.e., 
\begin{equation}\label{int2}
    \int_\Omega \left(u_{\bm\beta} +\gamma\, u-f\right)^2(\bx)\,d\bx.
\end{equation}
However, the exact solution $u$ and hence the integrand in \cref{int2} are unknown. 

To circumvent this difficulty, adaptive numerical integration was proposed and studied in the context of the deep Ritz method for linear elasticity equation in \cite{LiCaRa23}. Adaptive numerical integration is based on a composite numerical integration
\[
\sum_{K\in\mathcal{T}} \mathcal{Q}_K(w)\approx \int_\Omega w(\bx)\,d\bx =\sum_{K\in\mathcal{T}} \int_K w(\bx)\,d\bx,
\]
where $\mathcal{T}=\{K:\,\, K \text{ is an open subdomain of } \Omega\}$ is a partition of $\Omega$ and $\mathcal{Q}_K(w)\approx \int_K w(\bx)\,d\bx$ denotes a quadrature rule over $K$. First, $\mathcal{Q}_K$ may vary on $K\in \mathcal{T}$. Second, its choice is one of the standard quadrature rules like the Gaussian quadrature or Newton-Cotes formulas such as the midpoint, trapezoidal, or Simpson rule (see \cite{Cai2022nonlinear}). In the case of the midpoint rule for all $K\in \mathcal{T}$, $\mathcal{Q}_K(w)=w(\mathbf{x}_K)\lvert K\rvert$, where $\mathbf{x}_K$ is the centroid of $K$ and $\lvert K\rvert$ is the $d$-dimensional measure of $K$. 

\begin{remark}\label{int3}
The LSNN approximation $u_N\in \cM(L,n)$ defined in \cref{L-NN} is a continuous piece-wise linear with respect to a partition of the domain $\Omega$, referred to as the physical partition in \cite{LiuCai1, LiuCai2, Cai2021DeepAdaptive}. The partition $\mathcal{T}$ is a ``mesh'' for numerical integration and is completely different from the physical partition of $u_N$. Therefore, the partition $\mathcal{T}$ differs from meshes of traditional numerical methods. Nevertheless, the partition $\mathcal{T}$ and the corresponding quadrature $\mathcal{Q}_K$ are important for accuracy of the approximation $u_N$ by providing accurate information of the exact solution. 
\end{remark}

The integrand in \cref{int2} has a derivative term $u_{\bm\beta}(\bx)$. At where $u$ is differentiable, we have 
\begin{equation}\label{dd}
u_{\bm\beta}(\bx)=\sum_{i=1}^d\beta_i(\bx) \dfrac{\partial u(\bx)}{\partial x_i}.
\end{equation}
Obviously, \cref{dd} is invalid at where the solution $u$ is discontinuous, and hence any NN method such as the PINNs in \cite{raissi2019physics} using discrete or auto differentiation based on \cref{dd} would lead to an unreasonable approximation to a discontinuous solution. This phenomenon was already reported by several researchers, e.g., \cite{chen21} for \cref{pde} with ${\bm\beta}=(1,1)$ and \cite{FuksTchelepi2020} for scalar nonlinear hyperbolic conservation laws (HCLs). This essential difficulty was overcome by introduction of a physics preserved discrete differential operator: the discrete directional differentiation operator for \cref{pde} in \cite{Cai2021linear} and the discrete divergence operator for nonlinear HCLs in \cite{Cai2022nonlinear}. For any $\bx\in \Omega$, the discrete differential operator $D_{\bm\beta}$ is defined by
 \begin{equation}\label{finite_diff}
     D_{\bm\beta}v(\mathbf{x})\coloneqq \frac{v(\bx)-v\big(\bx - \rho\bar{\bm{\beta}}(\bx)\big)}{\rho/|{\bm\beta}(\mathbf{x})|}\approx v_{\bm\beta}(\bx),
\end{equation}
where $|{\bm\beta(\mathbf{x})}|$ is the magnitude of $\bm\beta(\mathbf{x})$, $\bar{\bm\beta}(\mathbf{x})=\frac{\bm\beta(\mathbf{x})}{|\bm\beta(\mathbf{x})|}$ is the unit vector along $\bm\beta(\mathbf{x})$, and $0<\rho \ll 1$. That is, the directional derivative $v_{\bm\beta}$ in the $\bm{\beta}$ direction is approximated by the backward finite difference quotient with the ``mesh'' size $\rho/|{\bm\beta}(\mathbf{x})|$. Fundamentally, the discrete differentiation operator $D_{\bm\beta}$ ensures that the derivative is computed without crossing the discontinuous interface.

For each $E\in\mathcal{E}_{-}=\{E=\partial K\cap \Gamma_{-}:K\in\mathcal{T}\}$, let $\mathcal{Q}_E(w)$ denote a quadrature rule for integrand $w$ defined on $E$. The discrete LS functional is defined by
\begin{equation}\label{dls}
\mathcal{L}_{_{{\cal T}}}\big(v;{\bf f}\big)=\sum_{K\in\mathcal{T}}\mathcal{Q}_K\left((D_{\bm{\beta}}v+\gamma v-f)^2\right)+\sum_{E\in\mathcal{E}_{-}}\mathcal{Q}_E\left(|\bm{\beta}\! \cdot \!\bm{n}|(v-g)^2\right).
\end{equation}
Then the discrete least-squares NN approximation of problem \cref{pde} is to find ${u}^{_N}_{_{{\cal T}}}\in \cM(L,n)$ such that
 \begin{equation}\label{discrete_minimization_functional}
  \mathcal{L}_{_{{\cal T}}} \big({u}^{_N}_{_{{\cal T}}};{\bf f}\big) 
  = \min\limits_{v\in \cM(L,n)} \mathcal{L}_{_{{\cal T}}}\big(v;{\bf f}\big).
\end{equation}
 
\begin{remark}\label{bc}
The discrete least-squares NN approximation defined in \cref{discrete_minimization_functional} enforces the inflow boundary condition through penalization. Instead, we may impose the inflow boundary condition through the discrete differentiation operator $D_{\bm\beta}$ defined in \cref{finite_diff} by choosing proper $\rho$ for integration point $\bx$, that is close to the inflow boundary, so that $\bx - \rho\bar{\bm{\beta}}(\bx)$ belongs to $\mathcal{E}_{-}$. In this way, the boundary term of $\mathcal{L}_{_{{\cal T}}}\big(v;{\bf f}\big)$ in \cref{dls} may be dropped.
\end{remark}

\section{Error estimates}\label{error estimates}
In this section, we provide error estimates for approximation by ReLU NN functions of the solution of the linear advection-reaction equation with a discontinuity interface. To this end, we note first that the solution of the problem is discontinuous if the inflow boundary data $g$ is discontinuous.

\begin{lemma}\label{l:solution}
For $d=2$, we assume that the inflow boundary data $g$ is discontinuous at $\bx_0\in \Gamma_{-}$ with values $g^+(\mathbf{x}_0)$ and $g^-(\mathbf{x}_0)$ from different sides. Let $I$ be the streamline of the vector field ${\bm\beta}$ emanating from $\mathbf{x}_0$ and let $\mathbf{x}(s)$ be a parameterization of $I$, i.e.,
\begin{equation}
    \frac{d\mathbf{x}(s)}{ds}={\bm\beta}(\mathbf{x}(s)),\quad \mathbf{x}(0)=\mathbf{x}_0.
\end{equation} 
Then the solution $u$ of \cref{pde} is discontinuous on $I$ with jump described as
\begin{multline}\label{jump}
    |u^+(\mathbf{x}(s))-u^-(\mathbf{x}(s))|=\\
    \exp\left(-\int_{0}^s\gamma(\mathbf{x}(t))\,dt\right)\left\lvert \int_0^s \exp\left(\int_{0}^t\gamma(\mathbf{x}(r))\,dr\right)(f^+(\mathbf{x}(t))-f^-(\mathbf{x}(t)))\,dt +g^+(\mathbf{x}_0)-g^-(\mathbf{x}_0)\right\rvert,
\end{multline}
where 
$u^+(\mathbf{x}(s))$ and $u^-(\mathbf{x}(s))$ are the solutions, and $f^+(\mathbf{x}(t))$ and $f^-(\mathbf{x}(t))$ are the values of $f$ of \cref{pde} along $I$ from different sides, respectively.
\end{lemma}

\begin{proof}
Along the interface $I$, by the definition of the directional derivative, we have
\[
u_{\bm\beta}(\mathbf{x}(s))=\frac{d}{ds}u(\mathbf{x}(s)).
\]
Thus the solutions $u^{\pm}(\mathbf{x}(s))$ along the interface $I$ satisfies the linear ordinary differential equations
\begin{equation}
\left\{\begin{array}{rcll}
\frac{d}{ds}u^{\pm}(\mathbf{x}(s))+\gamma(\mathbf{x}(s))u^{\pm}(\mathbf{x}(s))&=& f^\pm(\mathbf{x}(s)), &\text{ for }\, s>0, \\[2mm]
u^{\pm}(\mathbf{x}(0))&=&u^{\pm}(\mathbf{x}_0)=g^{\pm}(\mathbf{x}_0) &
    \end{array}\right.
\end{equation}
whose solutions are given by
\[ 
    u^{\pm}(\mathbf{x}(s))=\exp\left(-\int_{0}^s\gamma(\mathbf{x}(t))\,dt\right)\left[\int_0^s \exp\left(\int_{0}^t\gamma(\mathbf{x}(r))\,dr\right)f^\pm(\mathbf{x}(t))\,dt+g^{\pm}(\mathbf{x}_0)\right].
\] 
Hence, $u$ is discontinuous on $I$ with the jump \cref{jump}. This completes the proof of the lemma.
\end{proof}

\begin{remark}\label{ist}
For $d=3$, we assume that the inflow boundary data $g$ is discontinuous along a curve $C(t)\subset\Gamma_-$. In this case, the collection of the streamlines $\mathbf{x}(s)$ of the vector field ${\bm\beta}$ starting at all $\mathbf{x}_0=C(t)$ forms a surface $I(s,t)$. Then the solution $u$ of \eqref{pde} is discontinuous on the surface $I(s,t)$ with jump as in \eqref{jump} for every $\mathbf{x}_0=C(t)$.
\end{remark}

Let the discontinuity interface $I$ (in $\mathbb{R}^d$) as the union of streamlines of the vector field $\bm{\beta}$ as in \cref{ist} divide the domain $\Omega$ into two nonempty subdomains $\Omega_1$ and $\Omega_2$ (see \cref{Interface} for $d=2$):
\[\Omega=\Omega_1\cup\Omega_2\,\,\text{ and }\,\,I=\Omega_1\cap\Omega_2,\]
so that the solution $u$ is piecewise smooth with respect to the partition $\{\Omega_1,\Omega_2\}$. We assume that every streamline has a finite length. 

Furthermore, We assume that the jump of the solution is constant. Hence $u$ can be decomposed into
\begin{equation}\label{decop}
  u(\bx)=\hat{u}(\bx) + \chi(\mathbf{x}),  
\end{equation}
where $\hat{u}$ is continuous and piecewise smooth on $\Omega$, and $\chi(\bx)$ is the piecewise constant function defined by 
 \[
 \chi(\bx)=\left\{\begin{array}{rl}
 \alpha_1, & \bx \in \Omega_1,\\[2mm]
 \alpha_2, & \bx \in \Omega_2
 \end{array}
 \right.
 \]
with $\alpha_1=g^-(\bx_0)$ and $\alpha_2=g^+(\bx_0)$. Given $\varepsilon>0$, we approximate $I$ by a connected series of hyperplanes $\bm{\xi}_i\cdot\mathbf{x}-b_i=0$ in $\Omega_1$ for $i=1,2,\ldots, k$ such that $\bm{\xi}_i$ point toward $\Omega_2$ and the translation of $\bm{\xi}_i\cdot\mathbf{x}-b_i=0$ in the direction of $\bm{\xi}_i$ by $\varepsilon$ contains $I$ with $\bm{\xi}_i\cdot\mathbf{x}-b_i=0$ (see \cref{approximation}). By normalizing $\bm{\xi}_i$, we may assume $\lvert\bm{\xi}_i\rvert=1$. We now divide $\Omega$ by hyperplanes passing through the intersections of $\bm{\xi}_i\cdot\mathbf{x} =b_i$. Let $\Upsilon_i$ denote the subdomains determined by this process. (see \cref{subdomains}, two dotted lines divide the domain into three subdomains $\Upsilon_i$ for $i=1,2,3$.)

\begin{figure}[htbp]
\centering
\subfigure[An approximation of $I$ by a  connected series of hyperplanes ($k=3$)\label{approximation}]{
\begin{minipage}[t]{0.7\linewidth}
\centering
\begin{tikzpicture}[scale=0.8, transform shape]
\draw [] (0,5)-- (5,5);
\draw [] (0,5)-- (0,0);
\draw [] (0,0)-- (5,0);
\draw [] (5,0)-- (5,5);
\draw (3.5,0) arc[start angle=0, end angle=90, radius=3.5];
\draw [] (0,3.5)-- (2,2.872);
\draw [] (2,2.872)--(3.2,1.418);
\draw [] (3.5,0)--(3.2,1.418);

\draw [] (0,3.8)-- (2.116,3.045);
\draw [] (3.386,1.507)-- (2.116,3.045);
\draw [] (3.386,1.507)-- (3.704,0);

\node[] at (1,1) {$\Omega_{1}$};
\node[] at (4,4) {$\Omega_{2}$};
\node[] at (3,2.5) {$I$};
\draw[-latex] (3.2,2)--(7.7,2.8);

\node[] at (0.85,2.7) {$\bxi_1\cdot\bx =b_1$};
\node[] at (2,1.5) {$\bxi_2\cdot\bx =b_2$};
\node[] at (2.5,0.4) {$\bxi_3\cdot\bx =b_3$};

\draw (10,1) arc[start angle=0, end angle=90, radius=3.5];
\draw [] (6.5,4.5)-- (10,1);

\draw [] (7.5,5.5)--(11,2);
\draw [] (10,1)--(11,2);
\draw [] (6.5,4.5)--(7.5,5.5);

\draw[decoration={brace,raise=3pt,amplitude=3pt},decorate, semithick]
  (6.5,4.5) -- node[xshift=-0.3cm, yshift=0.35cm] {$\varepsilon$} (7.5,5.5);

\end{tikzpicture}
\end{minipage}
}%
\\
\subfigure[An interface $I$ in $\Omega\subset\mathbb{R}^{d}$\label{Interface}]{
\begin{minipage}[t]{0.4\linewidth}
\centering
\begin{tikzpicture}[scale=0.8, transform shape]
\draw [] (0,5)-- (5,5);
\draw [] (0,5)-- (0,0);
\draw [] (0,0)-- (5,0);
\draw [] (5,0)-- (5,5);
\draw (3.5,0) arc[start angle=0, end angle=90, radius=3.5];

\node[] at (1,1) {$\Omega_{1}$};
\node[] at (4,4) {$\Omega_{2}$};
\node[] at (3,2.5) {$I$};
\end{tikzpicture}
\end{minipage}%
}%
\hspace{0.2in}
\subfigure[A partition of $\Omega$ through the intersections of $\bxi_i\cdot\bx =b_i$ for $i=1,2,3$ ($k=3$)\label{subdomains}]{
\begin{minipage}[t]{0.4\linewidth}
\centering
\begin{tikzpicture}[scale=0.8, transform shape]
\draw [] (0,5)-- (5,5);
\draw [] (0,5)-- (0,0);
\draw [] (0,0)-- (5,0);
\draw [] (5,0)-- (5,5);
\draw (3.5,0) arc[start angle=0, end angle=90, radius=3.5];
\draw [] (0,3.5)-- (2,2.872);
\draw [] (2,2.872)--(3.2,1.418);
\draw [] (3.5,0)--(3.2,1.418);

\draw [] (0,3.8)-- (2.116,3.045);
\draw [] (3.386,1.507)-- (2.116,3.045);
\draw [] (3.386,1.507)-- (3.704,0);

\draw [dotted, very thick] (0,0.0723)--(3.446,5);
\draw [dotted, very thick] (0.256,0)--(5,2.267);

\node[] at (-0.4,0) {$\Omega_{1}$};
\node[] at (5.4,5) {$\Omega_{2}$};
\node[] at (3,2.5) {$I$};
\node[] at (1.4,4.2) {$\Upsilon_1$};
\node[] at (4,3.8) {$\Upsilon_2$};
\node[] at (4.3,0.9) {$\Upsilon_3$};
\node[] at (0.85,2.7) {$\bxi_1\cdot\bx =b_1$};
\node[] at (2,1.5) {$\bxi_2\cdot\bx =b_2$};
\node[] at (2.5,0.4) {$\bxi_3\cdot\bx =b_3$};
\end{tikzpicture}
\end{minipage}%
}%
\\
\subfigure[The hyperplanes $\bxi_i\cdot\bx=b_i$, $\bxi_i\cdot\bx =b_i+\varepsilon$ from $p_i(\mathbf{x})$\label{pi}]{
\begin{minipage}[t]{0.5\linewidth}
\centering
\begin{tikzpicture}[scale=0.8, transform shape]
\draw [dotted, very thick](3.5,0) arc[start angle=0, end angle=90, radius=3.5];
\draw [] (0,3.5)-- (3.5,0);
\draw [] (0.5,4)--(4,0.5);
\draw [] (1,4.5)--(4.5,1);
\draw [dotted, very thick] (3.5,0)--(4.5,1);
\draw [dotted, very thick] (0,3.5)--(1,4.5);

\draw[decoration={brace,raise=3pt,amplitude=3pt},decorate, semithick]
  (0,3.5) -- node[xshift=-0.3cm, yshift=0.35cm] {$\varepsilon$} (1,4.5);

\node[] at (4.5,-0.1) {$\bxi_i\cdot\bx =b_i$};
\node[] at (5.5,0.4) {$\bxi_i\cdot\bx =b_i+\varepsilon/2$};
\node[] at (5.8,1) {$\bxi_i\cdot\bx =b_i+\varepsilon$};

\end{tikzpicture}
\end{minipage}%
}%
\caption{An illustration of \cref{chi-curve lem}}
\end{figure}

\begin{lemma}\label{chi-curve lem}
Let $p(\bx)$ be the CPWL function \emph{(}see \cref{pi}\emph{)} defined by
\[
p(\mathbf{x})\coloneqq p_i(\bx)\coloneqq\alpha_1 +
\dfrac{\alpha_2-\alpha_1}{\varepsilon} \Big(\sigma(\bxi_i\cdot\bx -b_i) - \sigma(\bxi_i\cdot\bx -b_i-\varepsilon)\Big),\,\, \bx\in \Upsilon_i.
\]
(when there is a subdomain $\Upsilon_i$ either in $\Omega_1$ or $\Omega_2$ that does not contain any hyperplane $\bm{\chi}_i\cdot\mathbf{x}-b_i=0$, we define $p(\mathbf{x})=\alpha_1$ or $\alpha_2$ on $\Upsilon_i$. Then $\vertiii{\chi-p}_{\bm\beta}=0$ on $\Upsilon_i$. Hence, without loss of generality, we assume that each $\Upsilon_i$ contains a hyperplane $\bm{\chi}_i\cdot\mathbf{x}-b_i=0$ as in \cref{subdomains}.)
Then we have 
 \begin{equation}\label{chi-curve}
     \vertiii{\chi - p}_{\bm\beta}
     \leq \sqrt{2|I|}\, \big|\alpha_1-\alpha_2\big| \sqrt{\varepsilon}, 
 \end{equation}
where $|I|$ is the $d-1$ dimensional measure of the interface $I$.
\end{lemma}

\begin{proof}
For each $i$, let $|\mathcal{P}_i|$ denote the $d-1$ dimensional measure of the hyperplane $\bxi_i\cdot\bx=b_i$ in $\Upsilon_i$. It is easy to check that
\[
\|\chi- p_i\|^2_{0,\Upsilon_i}\le (\alpha_1-\alpha_2)^2|\mathcal{P}_i|\varepsilon,
\]
which implies
\begin{equation}\label{chi-curve a}
    \|\chi - p\|^2_{0,\Omega}=\sum_{i=1}^k\|\chi- p_i\|^2_{0,\Upsilon_i}\le\sum_{i=1}^k|\mathcal{P}_i|(\alpha_1-\alpha_2)^2\varepsilon\le |I| (\alpha_1-\alpha_2)^2 \varepsilon
\end{equation}
assuming $\sum_{i=1}^k|\mathcal{P}_i|\le |I|$. (we approximate the interface $I$ with a sufficiently large $k$ so that the assumption holds.)

We now prove
\begin{equation*}
    \|\chi_{\scriptsize\bm\beta} - p_{\scriptsize\bm\beta}\|^2_{0,\Omega}
     \leq |I| (\alpha_1-\alpha_2)^2 \varepsilon.
\end{equation*}
To this  end, for each $i$, let 
\[
\Upsilon_i^1= \{\mathbf{x}\in\Upsilon_i:0<\bxi_i\cdot\mathbf{x}-b_i<\varepsilon\}
\quad\mbox{and}\quad \Upsilon_i^2=\Upsilon_i\setminus \Upsilon_i^1.
\]
Clearly, $\chi_{\scriptsize\bm\beta}\equiv 0$ on $\Omega$, and $p_i$ is piecewise constant on $\Upsilon_i^2$. For each $i$, we construct a vector field $\bm{\beta}_i(\bx)$ on $\Upsilon_i$ such that for each $\bx\in \Upsilon_i^1$, $\bm{\beta}_i(\bx)$ is parallel to the hyperplane $\bxi_i\cdot\bx =b_i$ and that $\bm{\beta}(\bx)-\bm{\beta}_i(\bx)$ is parallel to $\bxi_i$. Then $(p_i)_{{\scriptsize\bm{\beta}}_i}\equiv0$ in $\Upsilon_i$ and
\begin{eqnarray*}
 \|(p_i)_{\scriptsize\bm\beta}\|_{0,\Upsilon_i}^2
 &= &\|(p_i)_{\scriptsize\bm\beta}-(p_i)_{{\scriptsize\bm\beta}_i}\|_{0,\Upsilon_i}^2
=\|(p_i)_{{\scriptsize\bm\beta}-{\scriptsize\bm\beta}_i}\|_{0,\Upsilon_i}^2=\|(p_i)_{{\scriptsize\bm\beta}-{\scriptsize\bm\beta}_i}\|_{0,{\Upsilon}^1_i}^2 \\[2mm]
&\le & \int_{{\Upsilon}^1_i}\left(\frac{\alpha_2-\alpha_1}{\varepsilon}\bxi_i\cdot\varepsilon{\bxi_i}\right)^2\,d\bx \le(\alpha_1-\alpha_2)^2|\mathcal{P}_i|\varepsilon,
\end{eqnarray*}
where for the first inequality, we used the fact that on ${\Upsilon}^1_i$, the gradient of $p_i$ is $\frac{\alpha_2-\alpha_1}{\varepsilon}\bxi_i$ and further assume that the magnitude of $\bm{\beta}(\mathbf{x})-\bm{\beta}_i(\mathbf{x})$ is less than or equal to $\varepsilon\bm{\xi}_i$. (Since $\bm{\beta}(\mathbf{x})$ is a tangent vector of a streamline, when the interface $I$ is approximated by a sufficiently large number of hyperplanes $\bm{\chi}_i\cdot\mathbf{x}-b_i=0$ for $i=1,2,\ldots,k$, so that the secant vector $\bm{\beta}_i(\mathbf{x})$ approaches $\bm{\beta}(\mathbf{x})$, the assumption holds.)
Thus
\begin{equation}\label{chi-curve b}
\|\chi_{\scriptsize\bm\beta} - p_{\scriptsize\bm\beta}\|^2_{0,\Omega}=\sum_{i=1}^k\|(p_i)_{\scriptsize\bm\beta}\|^2_{0,\Upsilon_i}\le\sum_{i=1}^k|\mathcal{P}_i|(\alpha_1-\alpha_2)^2\varepsilon\le|I|(\alpha_1-\alpha_2)^2 \varepsilon.
\end{equation}
Now \cref{chi-curve} follows from \cref{chi-curve a} and \cref{chi-curve b}.
\end{proof}

\begin{theorem}\label{curve-theorem}
Let $u$ and $u_{_{N}}$ be the solutions of problems \cref{minimization1} and \cref{L-NN}, respectively. Then we have 
 \begin{equation}\label{tau_1-error}
 \vertiii{u-u_{_{N}}}_{\bm\beta}
 \leq C\left(\big|\alpha_1-\alpha_2\big| \sqrt{\varepsilon} + \inf_{v\in \cM(L,n)} \vertiii{\hat{u}+p-v}_{\bm\beta}
 \right),
 \end{equation}
where $\hat{u}\in C(\Omega)$ and $p$ are given in \cref{decop} and \cref{chi-curve lem}, respectively. Moreover, if the depth of ReLU NN functions in \cref{L-NN} is at least $\lceil \log_2(d+1)\rceil+1$, then for a sufficiently large integer $n$, there exists an integer $\hat{n}\leq n$ such that 
\begin{equation}\label{tau_1-error-2}
 \vertiii{u-u_{_{N}}}_{\bm\beta}
 \leq C\,\left(\big|\alpha_1-\alpha_2\big|\, \sqrt{\varepsilon} + \inf_{v\in \cM(\log,n-\hat{n})} \vertiii{\hat{u}-v}_{\bm\beta}
 \right),
 \end{equation}
 where $\cM(\log,n-\hat{n})=\cM(\lceil \log_2(d+1)\rceil+1,n-\hat{n})$.
\end{theorem}

\begin{proof}
For any $v\in \cM(L,n)$, it follows from \cref{decop}, the triangle inequality, and \cref{chi-curve lem} that
\begin{align*}
    \vertiii{u-v}_{\bm\beta}&=\vertiii{\chi-p+\hat{u}+p-v}_{\bm\beta}
    \le \vertiii{\chi-p}_{\bm\beta}+\vertiii{\hat{u}+p-v}_{\bm\beta}\\[2mm]
    &\le \sqrt{2|I|}\,|\alpha_1-\alpha_2|\sqrt{\varepsilon}+\vertiii{\hat{u}+p-v}_{\bm\beta}. &&
\end{align*}
Taking the infimum over all $v\in\cM(L,n)$, \cref{tau_1-error} is then a direct consequence of \cref{Cea}. 

To show the last statement, first, we note that for a sufficiently large integer $n$, by \cref{cpwl=relu}, there exists an integer $\hat{n}\leq n$ such that 
\[
p\in \cM(\log,\hat{n})=\cM(\lceil \log_2(d+1)\rceil+1,\hat{n}).
\]
Obviously we have $v+p\in \cM(\log,n)$ for any $v\in \cM(\log,n-\hat{n})$. Now, it follows from the coercivity and continuity of the homogeneous functional $\mathcal{L}\big(v;{\bf 0}\big)$ in \cref{equiv}, problems \cref{pde}, \cref{L-NN}, \cref{decop}, and the triangle inequality that
\begin{eqnarray*}
\alpha\, \vertiii{u-u_{_N}}^2_{\bm\beta}
&\le &  \mathcal{L}\big(u-u_{_N};{\bf 0}\big)
= \mathcal{L}\big(u_{_N};{\bf f}\big)
\le  \mathcal{L}\big(v+p;{\bf f}\big) = \mathcal{L}\big(u-v-p;{\bf 0}\big)\\[2mm]
&=&\mathcal{L}\big((\hat{u}-v) + (\chi -p);{\bf 0}\big)
\le M \,\vertiii{(\hat{u}-v) + (\chi -p)}^2_{\bm\beta} \\[2mm]
&\leq & 2M\left(\vertiii{(\hat{u}-v) }^2_{\bm\beta} + \vertiii{(\chi -p)}^2_{\bm\beta}\right),
\end{eqnarray*}
which, together with \cref{chi-curve lem}, implies the validity of \cref{tau_1-error-2}. This completes the proof of the theorem.
\end{proof}

\begin{remark}
\cref{curve-theorem} mainly focused on the depth of NNs. A ReLU NN architecture with less than $\lceil\log_2(d+1)\rceil+1$ layers, prescribed widths, exact weights and biases to approximate piecewise constant functions will be addressed in a forthcoming paper.
\end{remark}

\begin{remark}
The continuous part of the solution $\hat{u}$ can be approximated well by a shallow NN, i.e., with depth $L=2$ (see, e.g., \cite{DeVore1997,Petrushev1998,mao2023rates,oswald1990degree} in various norms).
\end{remark}

\begin{remark}\label{uat}
The estimate in \cref{tau_1-error} holds even for the shallow neural network ($L=2$). However, the second term of the upper bound, $\inf_{v\in \cM(2,n)} \vertiii{\hat{u}+p-v}_{\bm\beta}$, depends on the inverse of $\varepsilon$ (as well as the norm of the directional derivative of $\hat{u}+p-v$) because the $p$ has a sharp transition layer of width $\varepsilon$. For any fixed $n$, $\inf_{v\in \cM(2,n)} \|p-v\|_{\infty}$ could be large depending on the size of $\varepsilon$, even though the universal approximation theorem implies 
\[
\lim_{n\to \infty} \inf_{v\in \cM(2,n)} \|p-v\|_{\infty}=0.
\]
Moreover, as $\varepsilon$ approaches $0$, $p$ approaches $\chi$, which is discontinuous; hence, in practice, the universal approximation theorem does not guarantee the convergence of the problem. On the other hand, by \cref{cpwl=relu}, deeper networks (with depth at least $\lceil \log_2(d+1)\rceil+1$) are capable of approximating such functions. As an illustration of this, we define a CPWL function $p(x,y)$ with a sharp transition layer of width $4\varepsilon/\sqrt{5}$ on $[0,1]^2$ as follows (see \cref{target}),
 \[
 p(x,y)=\left\{\begin{array}{rl}
 -1+\frac{1}{\varepsilon}\left(y+\frac{1}{2}x-\frac{4}{5}+\varepsilon\right), &\text{if } y\ge x,\  y>-\frac{1}{2}x+\frac{4}{5}-\varepsilon,\ y\le-\frac{1}{2}x+\frac{4}{5}+\varepsilon,\\[2mm]
 -1+\frac{1}{\varepsilon}\left(\frac{1}{2}y+x-\frac{4}{5}+\varepsilon\right), &\text{if } y<x,\ y>-2(x+\varepsilon)+\frac{8}{5},\ y\le-2(x-\varepsilon)+\frac{8}{5},\\[2mm]
 -1, &\text{if } y\le -\frac{1}{2}x+\frac{4}{5}-\varepsilon,\ y\le -2(x+\varepsilon)+\frac{8}{5},\\[2mm]
 1, &\text{otherwise.}
 \end{array}
 \right.
 \]
In \cite{Cai2021linear}, we proved CPWL functions of the form $p(x,y)$ are 2--4--4--1 ReLU NN functions. Here and in what follows, by a $d$--$n_1$-- $\cdots$ --$n_{L-1}$--$c$ ReLU NN function, we mean a ReLU NN function from $\mathbb{R}^d$ to $\mathbb{R}^c$ with an $L$-layer representation with width $n_l$ of the $l^{\text{th}}$ hidden layer of the representation for $l=1,\ldots,L-1$. Therefore, we can expect 2--$n_1$--$n_2$--1 (depth $3=\lceil \log_2(d+1)\rceil+1$ for $d=2$) ReLU NN functions outperform 2--$n_3$--1 ReLU NN functions, and \cref{12 comparison figure} and \cref{12 comparison table} demonstrate this for approximating $p(x,y)$ with $\varepsilon=0.1,0.01,0.001$ by ReLU NN functions with depth 2, 3 and various numbers of neurons where we minimized the squared $L^2$ norm loss function (the midpoint rule on a
uniform mesh with mesh size $h=10^{-2}$) using the Adam optimization algorithm for 100000 iterations with the learning rate 0.004. Even though the relative errors for the one-hidden-layer ReLU NN function approximations decrease as the number of neurons increases (\cref{12 comparison table}), the graphs (\cref{158,158 trace}) exhibit oscillations near the location where the transition layer is formed (\cref{transition}). In contrast, the two-hidden-layer ReLU NN functions approximate the target function well (\cref{target,44,44 trace}. The breaking hyperplanes (\cref{158 breaking,44 breaking}) show where the transition layers are formed for both approximations. In particular, if we zoom in on \cref{44 breaking}, breaking hyperplanes are right around the location of the transition layer (two green dotted lines with width $4\varepsilon/\sqrt{5}$), which implies the two almost vertical planes in \cref{44}.

\begin{figure}[htbp]\label{12 comparison figure}
\centering
\subfigure[The location of the transition layer\label{transition}]{
\begin{minipage}[t]{0.4\linewidth}
\centering
\includegraphics[width=1.8in]{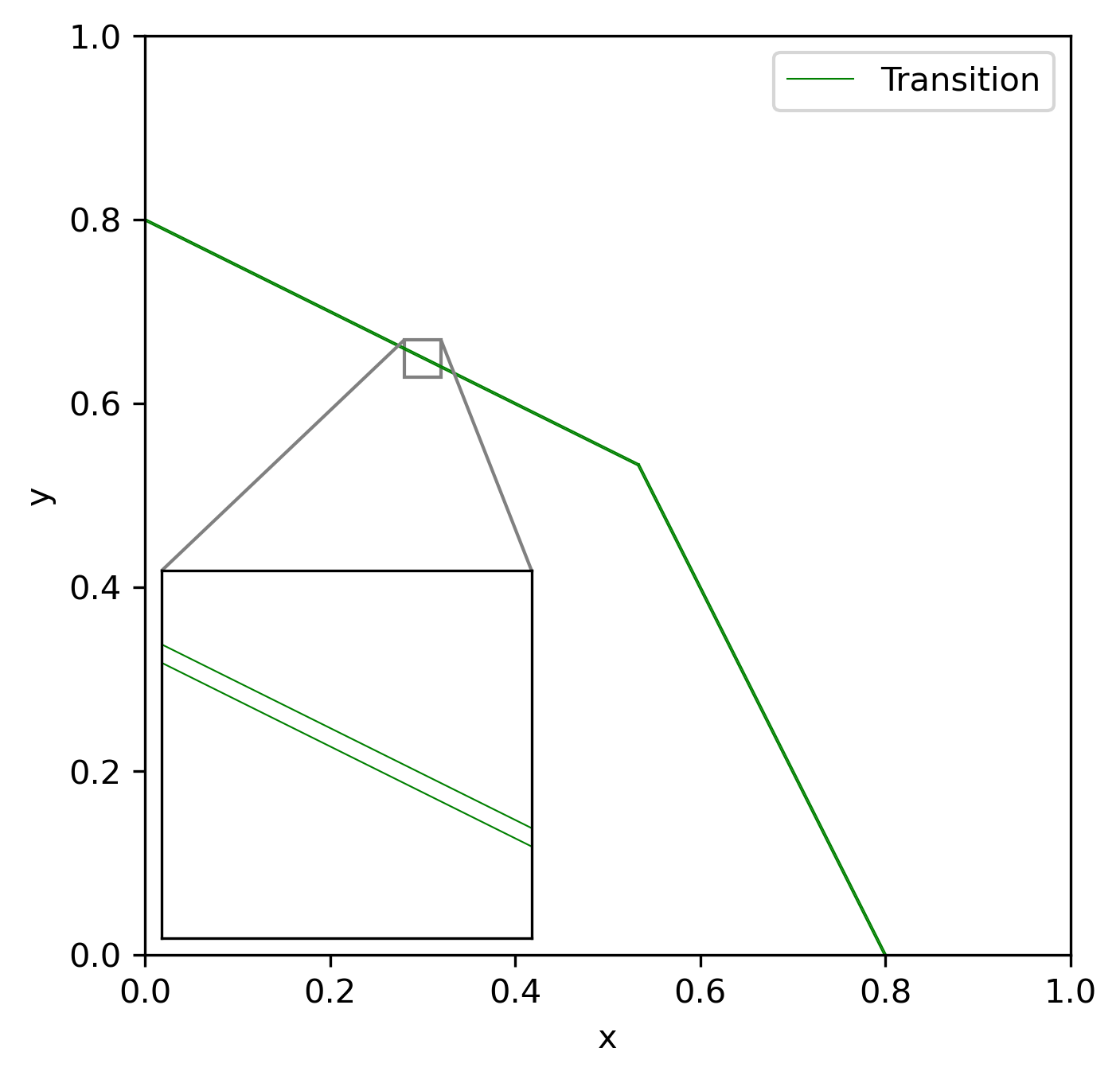}
\end{minipage}%
}%
\hspace{0.2in}
\subfigure[$p(x,y)$ with $\varepsilon=0.001$\label{target}]{
\begin{minipage}[t]{0.4\linewidth}
\centering
\includegraphics[width=1.8in]{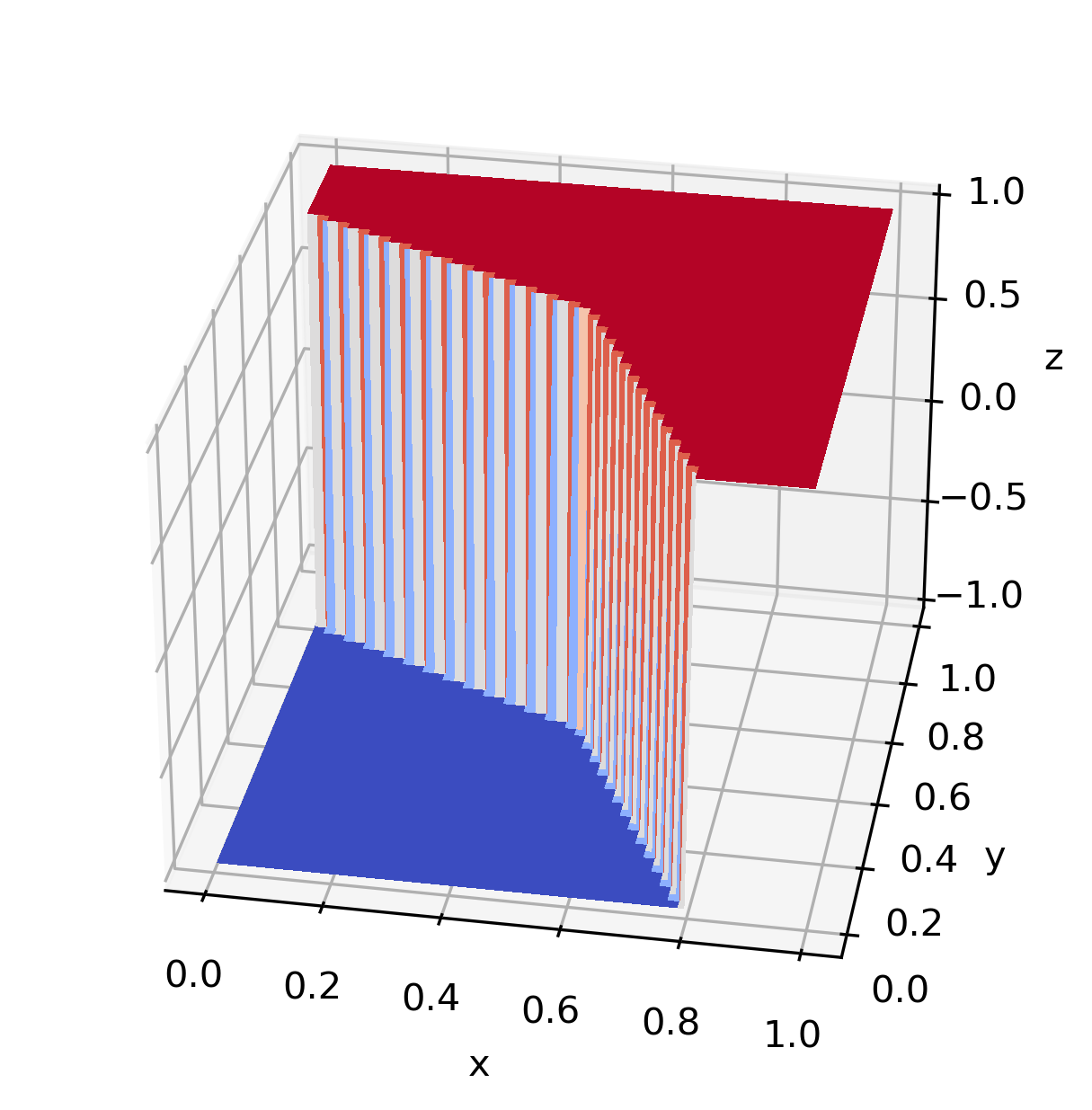}
\end{minipage}%
}%
\\
\subfigure[A 2--158--1 ReLU NN function approximation\label{158}]{
\begin{minipage}[t]{0.4\linewidth}
\centering
\includegraphics[width=1.8in]{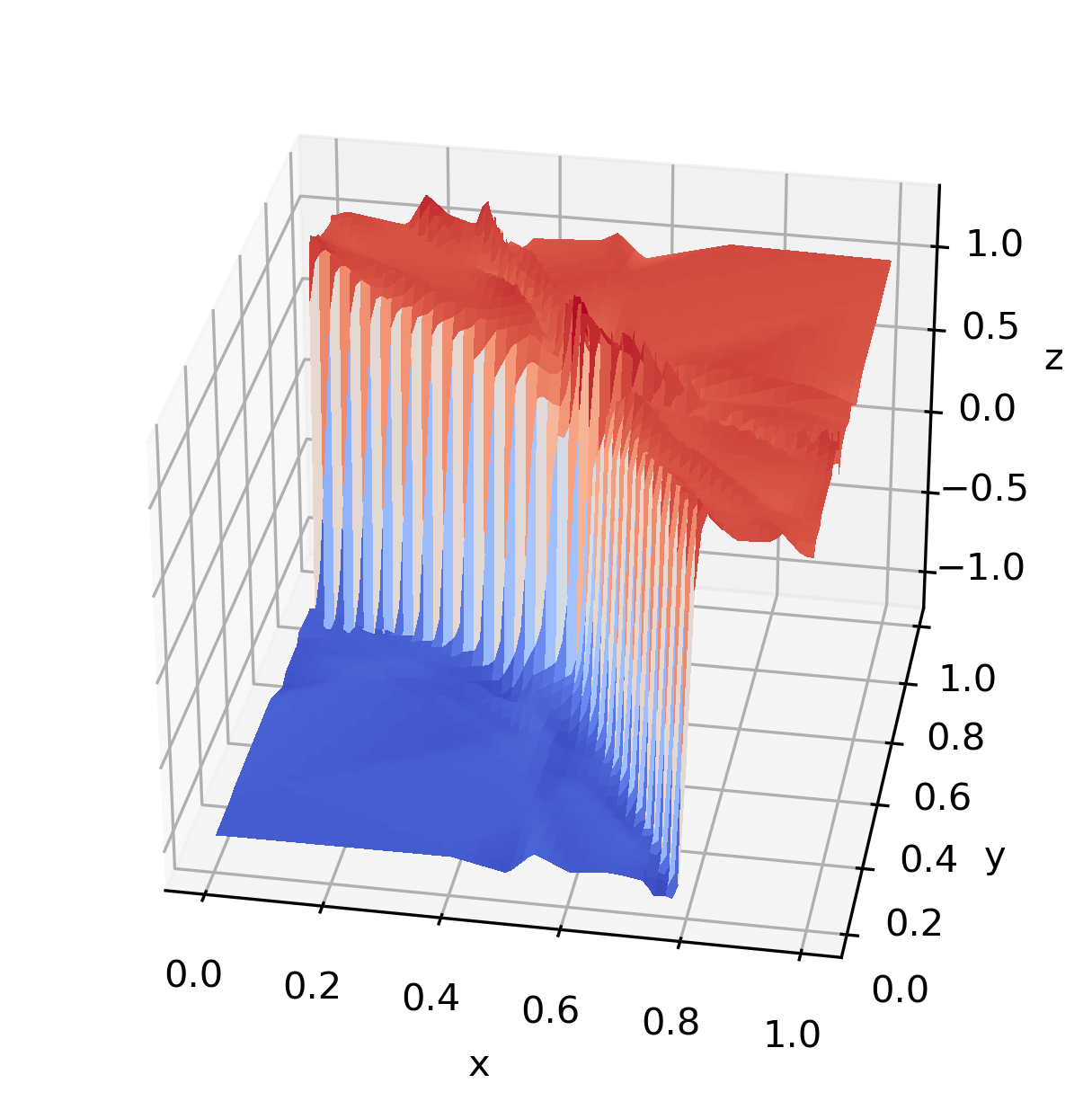}
\end{minipage}%
}%
\hspace{0.2in}
\subfigure[A 2--4--4--1 ReLU NN function approximation\label{44}]{
\begin{minipage}[t]{0.4\linewidth}
\centering
\includegraphics[width=1.8in]{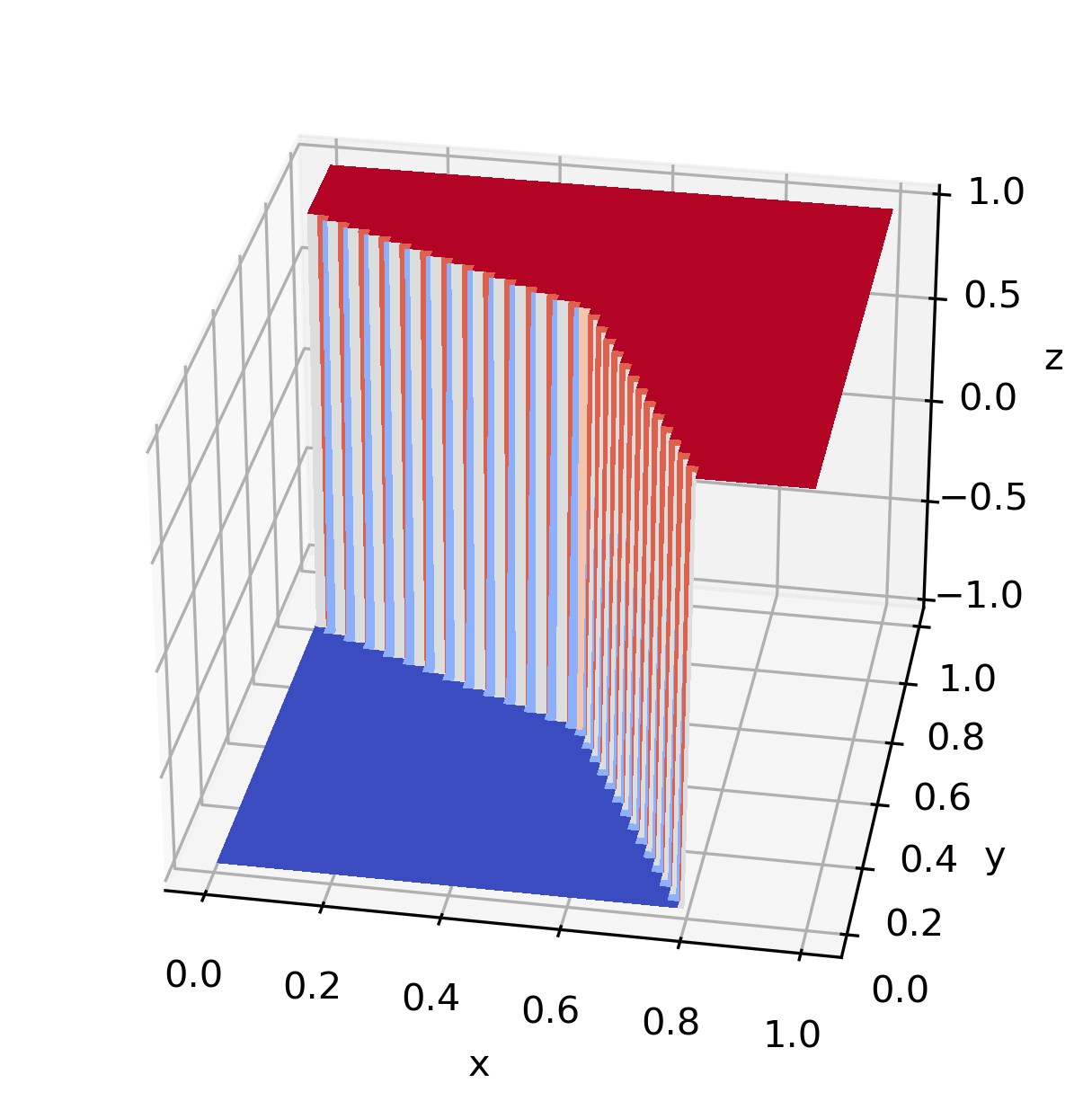}
\end{minipage}%
}%
\\
\subfigure[The trace of \cref{158} on $y=x$\label{158 trace}]{
\begin{minipage}[t]{0.4\linewidth}
\centering
\includegraphics[width=1.8in]{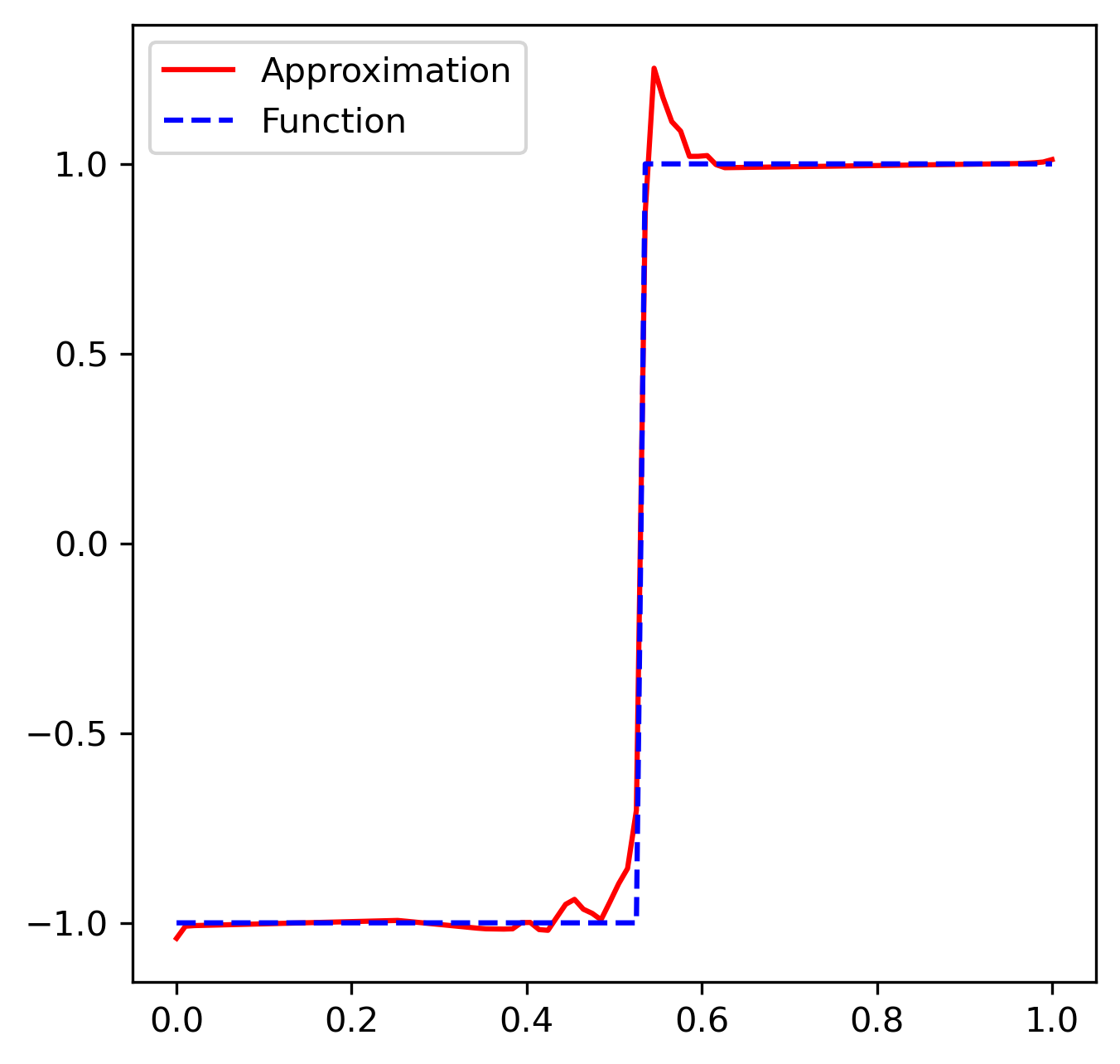}
\end{minipage}%
}%
\hspace{0.2in}
\subfigure[The trace of \cref{44} on $y=x$\label{44 trace}]{
\begin{minipage}[t]{0.4\linewidth}
\centering
\includegraphics[width=1.8in]{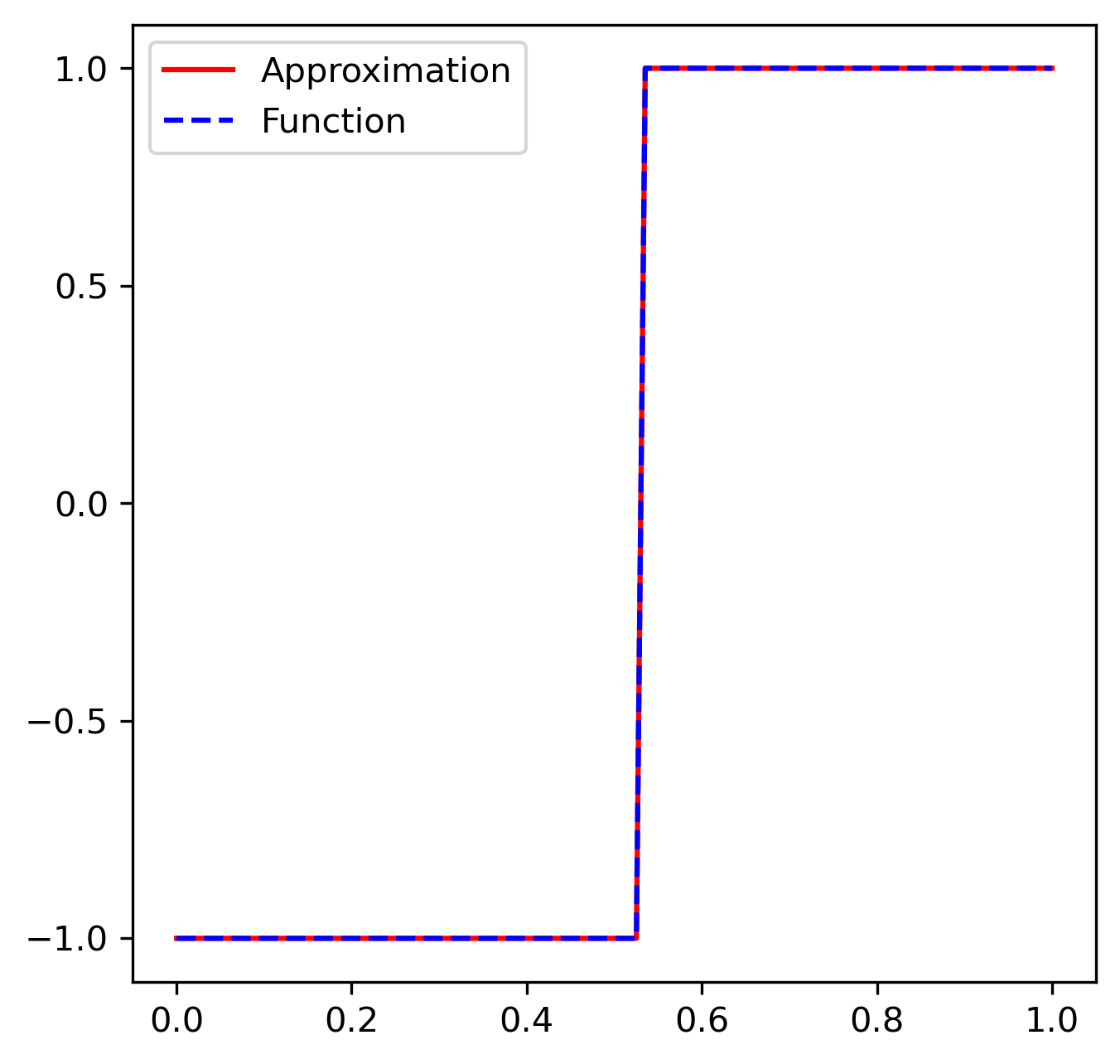}
\end{minipage}%
}%
\\
\subfigure[The breaking hyperplanes of the approximation in \cref{158}\label{158 breaking}]{
\begin{minipage}[t]{0.4\linewidth}
\centering
\includegraphics[width=1.8in]{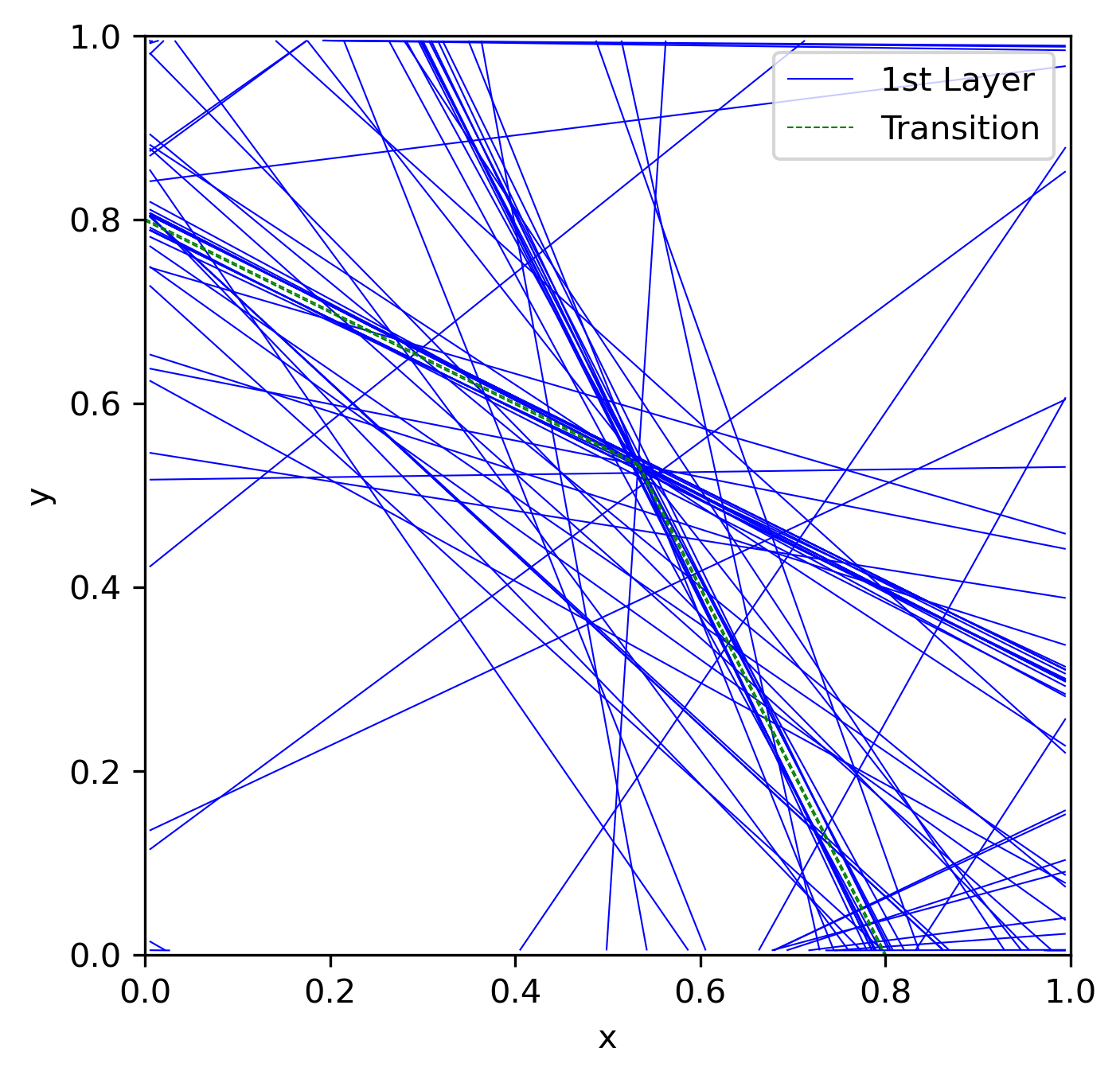}
\end{minipage}%
}%
\hspace{0.2in}
\subfigure[The breaking hyperplanes of the approximation in \cref{44}\label{44 breaking}]{
\begin{minipage}[t]{0.4\linewidth}
\centering
\includegraphics[width=1.8in]{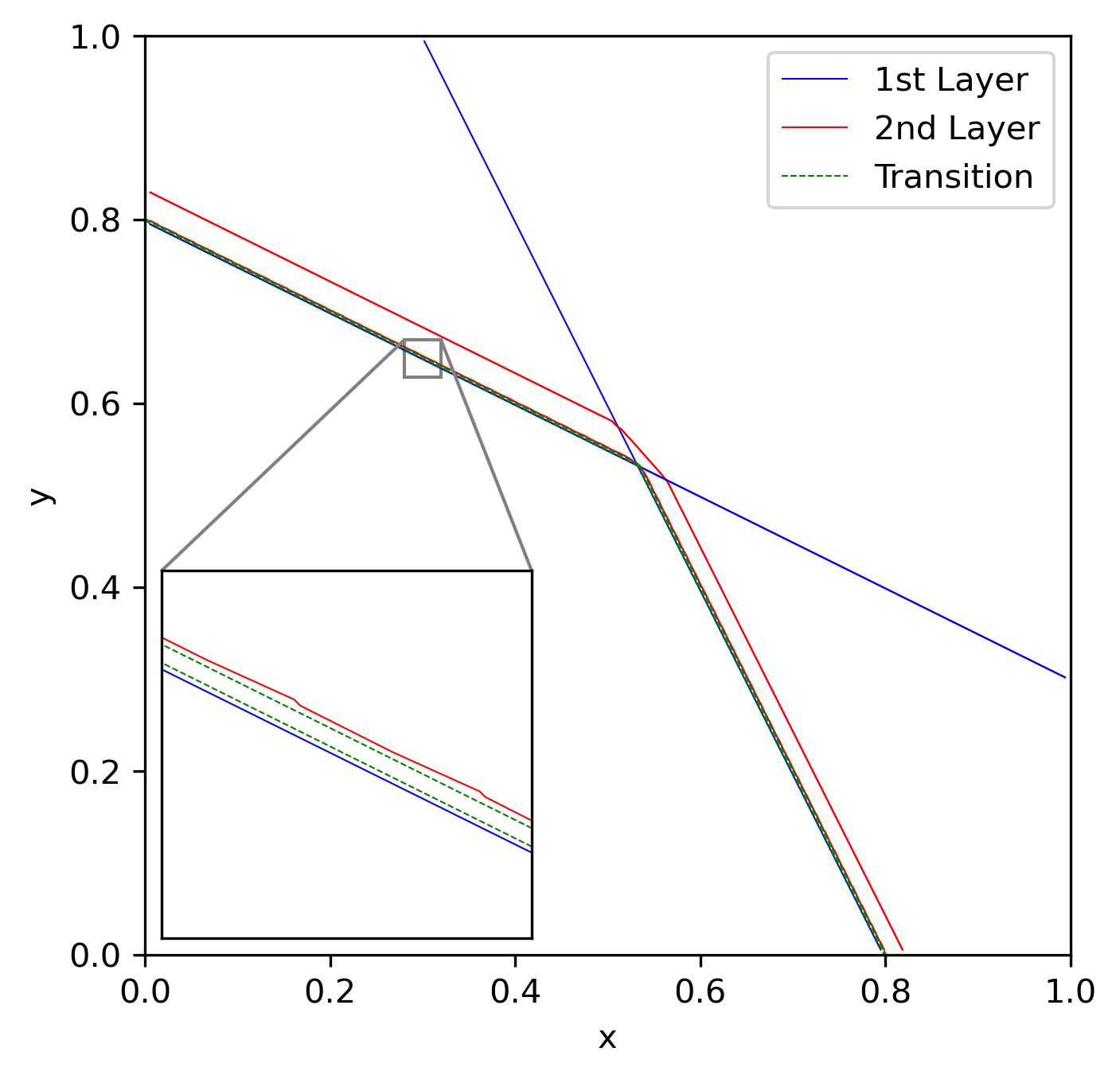}
\end{minipage}%
}%
\caption{$L^2$ norm approximation results of $p(x,y)$ with $\varepsilon=0.001$}
\end{figure}

\begin{table}[htbp]\label{12 comparison table}
\caption{Relative errors in the $L^2$ norm for approximating $p(x,y)$ with $\varepsilon=0.1,0.01,0.001$ by ReLU NN functions with depth $2$, $3$ and various numbers of neurons}
\centering
\begin{tabular}{|l|l|l|l|l|l|}
\hline
 & 2--8--1  &2--58--1& 2--108--1 & 2--158--1 & 2--4--4--1 \\ \hline
$\varepsilon=0.1$& $0.292446$ & $0.028176$ & $0.021200$ & $0.010259$ & $1.906427\times 10^{-7}$\\ \hline
$\varepsilon=0.01$& $0.254603$ &$0.078549$ &$0.065465$ & $0.030623$ & $7.536140\times 10^{-6}$\\ \hline
$\varepsilon=0.001$& $0.404299$ & $0.102757$ & $0.100136$ & $0.088885$ & $9.473783\times 10^{-7}$\\ \hline
\end{tabular}
\end{table}
Finally, we also note that $\mathcal{M}(2)\subsetneq\mathcal{M}(\lceil \log_2(d+1)\rceil+1)$ when $d\ge 2$ (see \cite{DeVore2021}).
\end{remark}

\begin{lemma}\label{uutn}
Let $u$, $u_{_N}$, and $u_{_\cT}^{_N}$ be the solutions of problems \cref{minimization1}, \cref{L-NN}, and \cref{discrete_minimization_functional}, respectively. Then there exist positive constants $C_{1}$ and $C_{2}$ such that
\begin{equation}\label{Cea-L-d}
\begin{split}
     \vertiii{u-u^{_N}_{_\cT}}_{\bm\beta}
    \le\  &C_{1}\,\left(\big|(\cL-\cL_{_\cT})(u_{_N}-u_{_\cT}^{_N}, {\bf 0})\big|
    + \big|(\cL-\cL_{_\cT})(u-u_{_N}, {\bf 0})\big|
    \right)^{1/2}\\
    &+C_{2}\,\left(\big|\alpha_1-\alpha_2\big|\, \sqrt{\varepsilon} + \inf_{v\in \cM(L,n)} \vertiii{\hat{u}+p-v}_{\bm\beta}
 \right).
\end{split}
\end{equation}
\end{lemma}

\begin{proof}
By the triangle inequality
\[\vertiii{u-u^{_N}_{_\cT}}_{\bm\beta}\le\vertiii{u-u_{_N}}_{\bm\beta}+\vertiii{u_{_N}-u^{_N}_{_\cT}}_{\bm\beta},\]
and the fact that
\[\vertiii{u_{_N}-u^{_N}_{_\cT}}_{\bm\beta}
    \le C_{1}\,\left(\left(\big|(\cL-\cL_{_\cT})(u_{_N}-u_{_\cT}^{_N}, {\bf 0})\big|
    + \big|(\cL-\cL_{_\cT})(u-u_{_N}, {\bf 0})\big|
    \right)^{1/2}+\vertiii{u-u_{_N}}_{\bm\beta}\right)\]
from the proof of Lemma~3.4 in \cite{Cai2021linear}, \cref{Cea-L-d} is a direct consequence of \cref{curve-theorem}.
\end{proof}

\section{Numerical experiments}\label{numerics}
In this section, we report numerical results for both two- and three-dimensional test problems with piecewise constant, or variable advection velocity fields. All numerical experiments have rectangular domains, and we used numerical integration (the midpoint rule) to implement the scheme in \cref{L-NN} (see \cite{Cai2021linear}). Numerical integration used a uniform mesh with mesh size $h=10^{-2}$. In \cref{finite_diff}, we set $\rho=h/4$ (except for the last test problem, which used $\rho=h/12$). We used the Adam optimization algorithm \cite{kingma2015} to iteratively solve the discrete minimization problem \eqref{discrete_minimization_functional}. For each numerical experiment, the learning rate started with 0.004 and was reduced by half for every 50000 iterations. Due to the possibility of the neural network getting trapped in a local minimum, we first trained the network with 5000 iterations 10 times, chose the weights and biases with the minimum loss function value, and trained further to get the results.

In \cref{2d test1 table,2d test2 table,2d test3 table,2d test4 table,2d test5 table,3d test1 table,3d test2 table,3d test3 table}, parameters indicate the total number of weights and biases, and $1/2$ in $\mathcal{L}^{1/2}$ for the relative error in the LS functional indicates the square root. We employed ReLU NN functions with width $n$ and depth $3=\lceil \log_2(d+1)\rceil+1$ for $d=2,\,3$. The basic principle for choosing the number of neurons is to start with a small number and increase the number to obtain a better approximation. For an automatic approach to design the architecture of deep neural networks (DNNs) for a given problem with a prescribed accuracy, see the recent work on the adaptive neural network method in \cite{LiuCai1, Cai2021DeepAdaptive}.

\subsection{Two-dimensional problems}
We present numerical results for five two-dimensional test problems with piecewise constant or variable advection velocity fields. The fifth test problem compares the LSNN method to other relevant methods. All five test problems are defined on the domain $\Omega=(0,1)^2$ with $\gamma=1$ (except for the fifth test problem with $\gamma=0.1$), and the exact solutions are the same as the right-hand side functions, $u(x,y)=f(x,y)$, which are step functions (except for the fourth and fifth test problems) along 3-line segment, 4-line segment, or curved interfaces. By \cref{curve-theorem}, the LSNN method with 3-layer ReLU NN functions leads to 
\[
 \vertiii{u-u_{_{N}}}_{\bm\beta}
 \leq C\,\big|\alpha_1-\alpha_2\big|\, \sqrt{\varepsilon},
\]
because the continuous part of the solution $\hat{u}$ in \eqref{decop} is zero (again, except for the fourth and fifth test problems).

\subsubsection{A problem with a 3-line segment interface}\label{2d test1 section}
\begin{figure}[htbp]
\centering
\subfigure[The interface\label{interface3}]{
\begin{minipage}[t]{0.4\linewidth}
\centering
\includegraphics[width=1.8in]{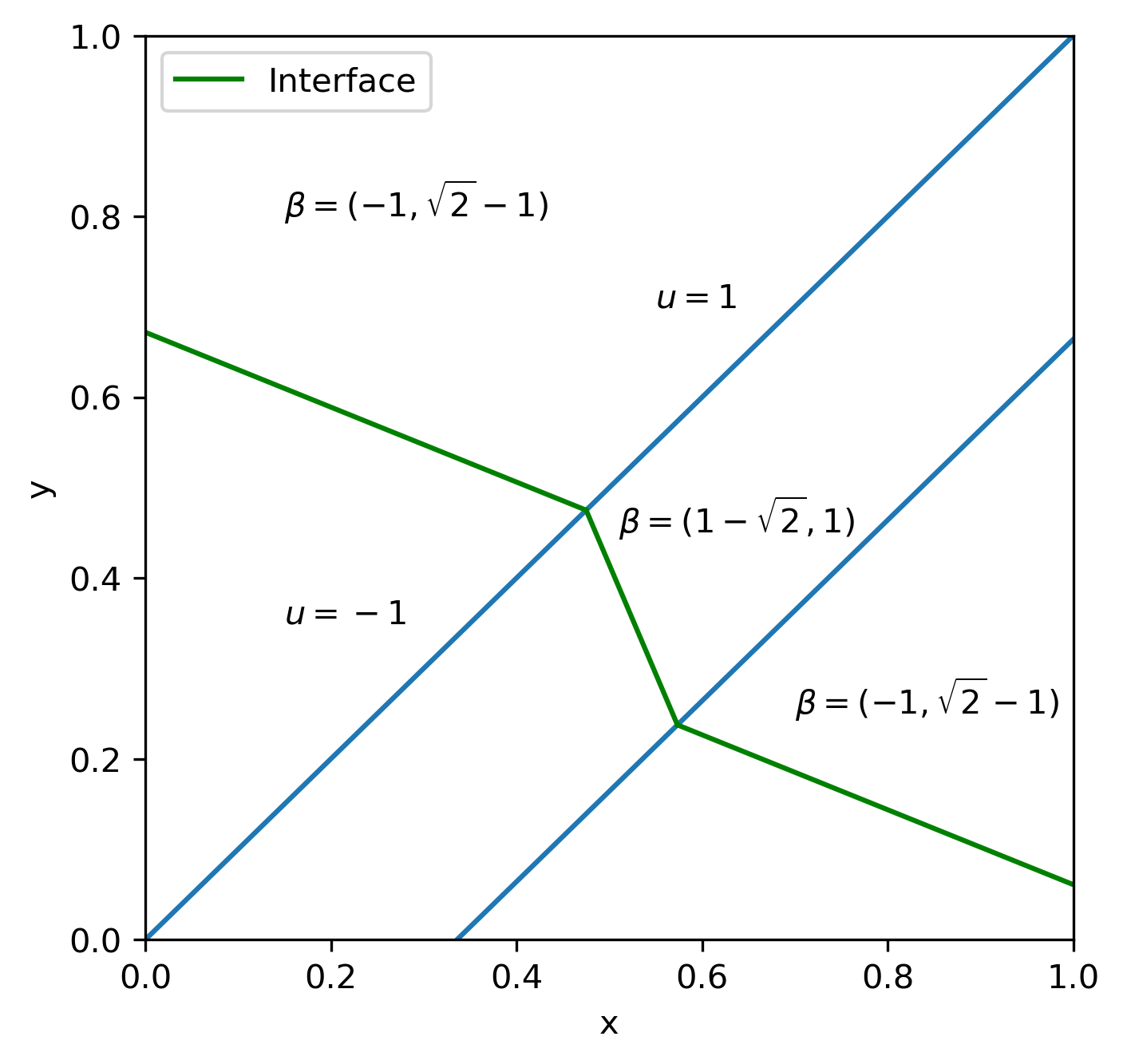}
\end{minipage}%
}%
\hspace{0.2in}
\subfigure[The exact solution\label{comparison31}]{
\begin{minipage}[t]{0.4\linewidth}
\centering
\includegraphics[width=1.8in]{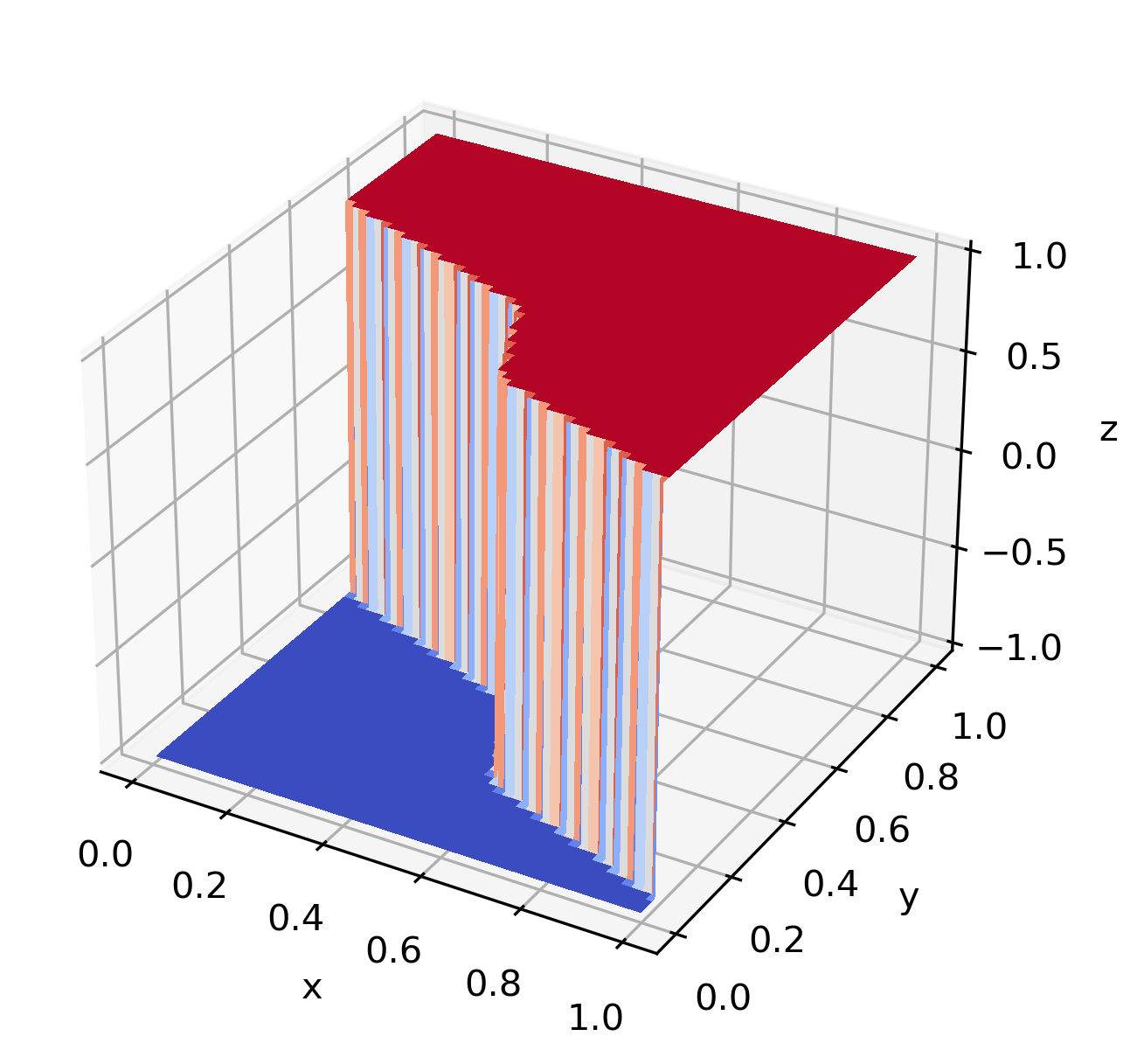}
\end{minipage}%
}%
\\
\subfigure[A 2--300--1 ReLU NN function approximation\label{comparison3_one}]{
\begin{minipage}[t]{0.4\linewidth}
\centering
\includegraphics[width=1.8in]{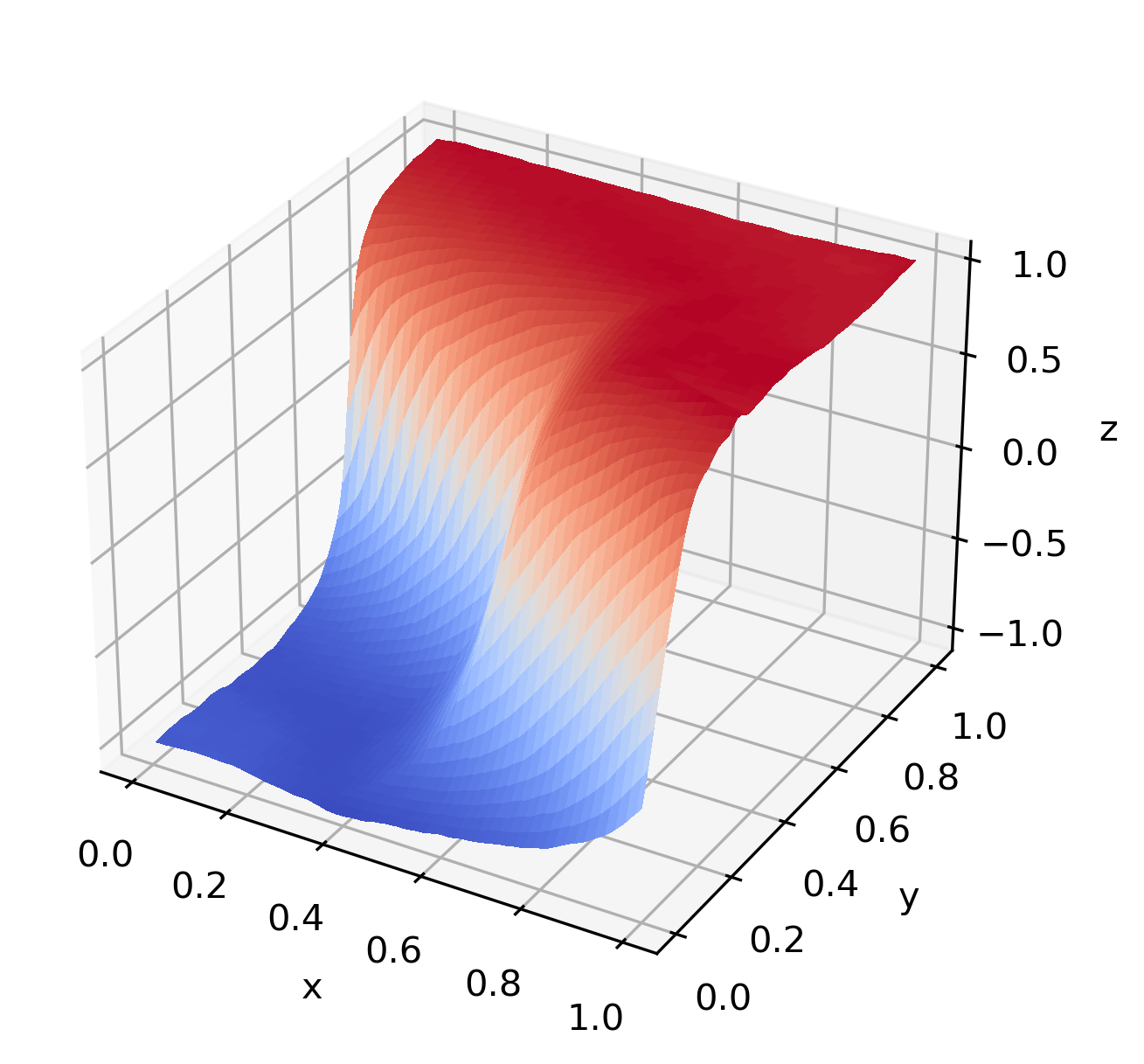}
\end{minipage}%
}%
\hspace{0.2in}
\subfigure[A 2--5--5--1 ReLU NN function approximation\label{comparison32}]{
\begin{minipage}[t]{0.4\linewidth}
\centering
\includegraphics[width=1.8in]{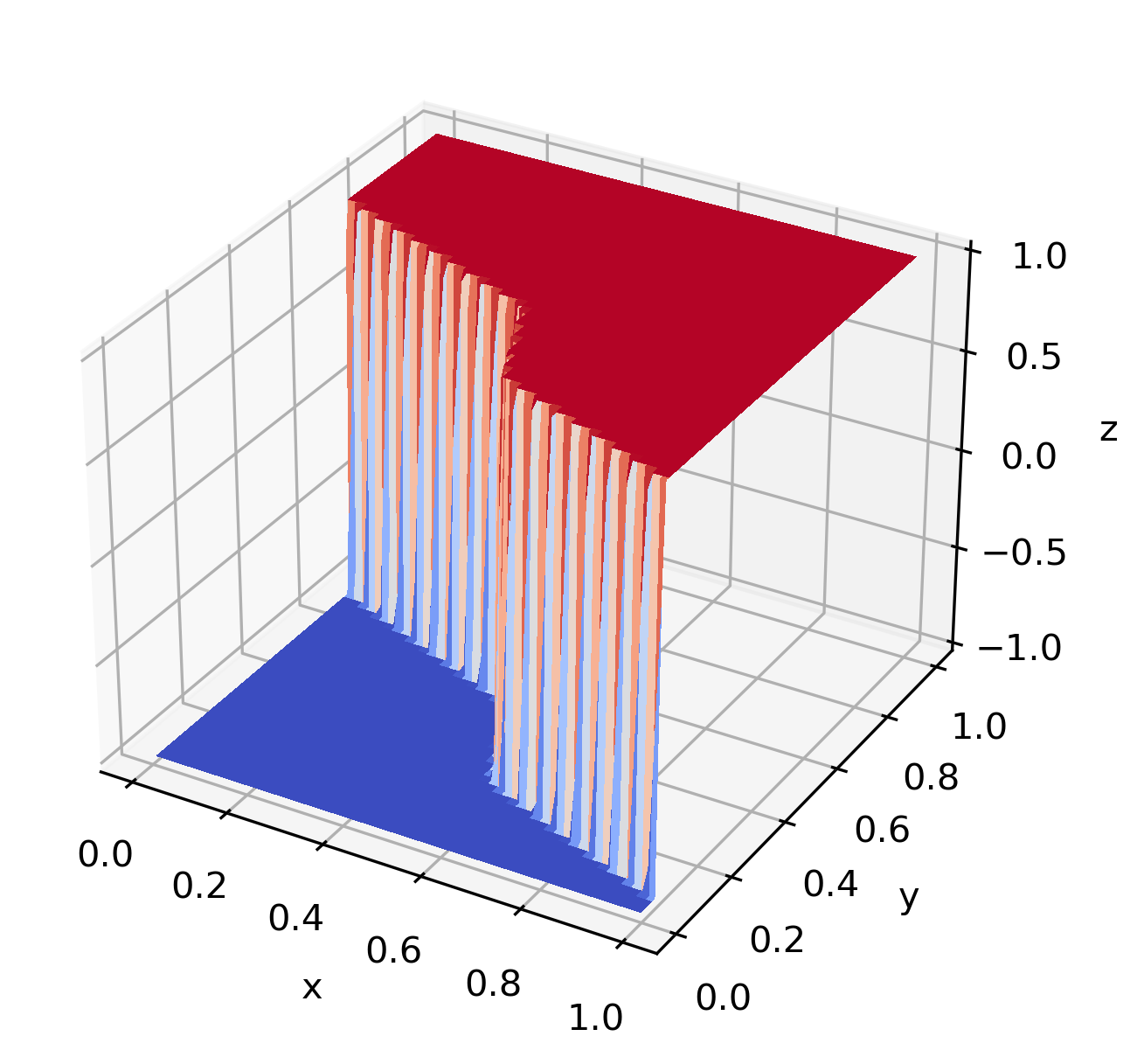}
\end{minipage}%
}%
\\
\subfigure[The trace of Figure \ref{comparison3_one} on $y=x$\label{vertical3_one}]{
\begin{minipage}[t]{0.4\linewidth}
\centering
\includegraphics[width=1.8in]{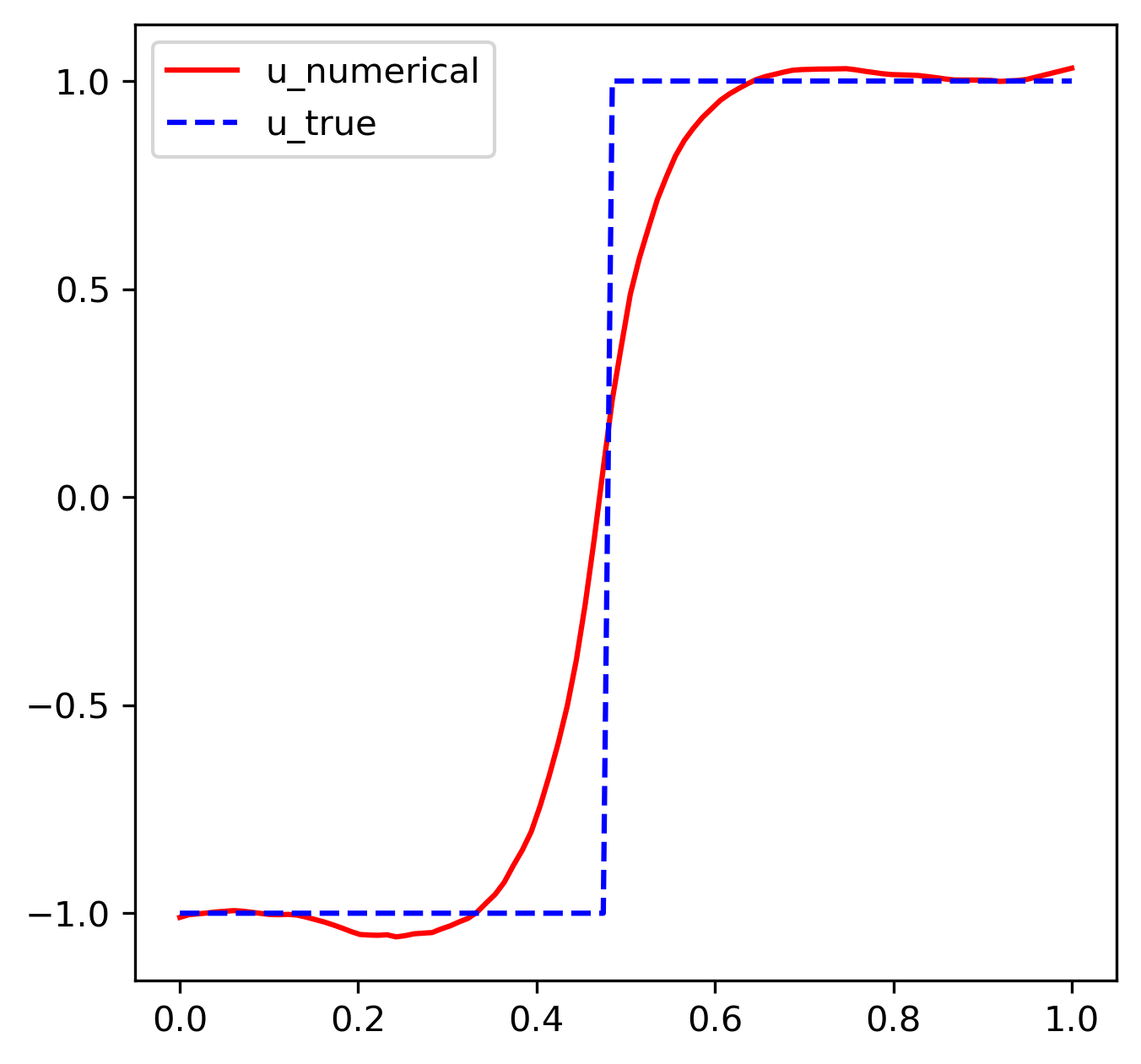}
\end{minipage}%
}%
\hspace{0.2in}
\subfigure[The trace of Figure \ref{comparison32} on $y=x$\label{vertical3}]{
\begin{minipage}[t]{0.4\linewidth}
\centering
\includegraphics[width=1.8in]{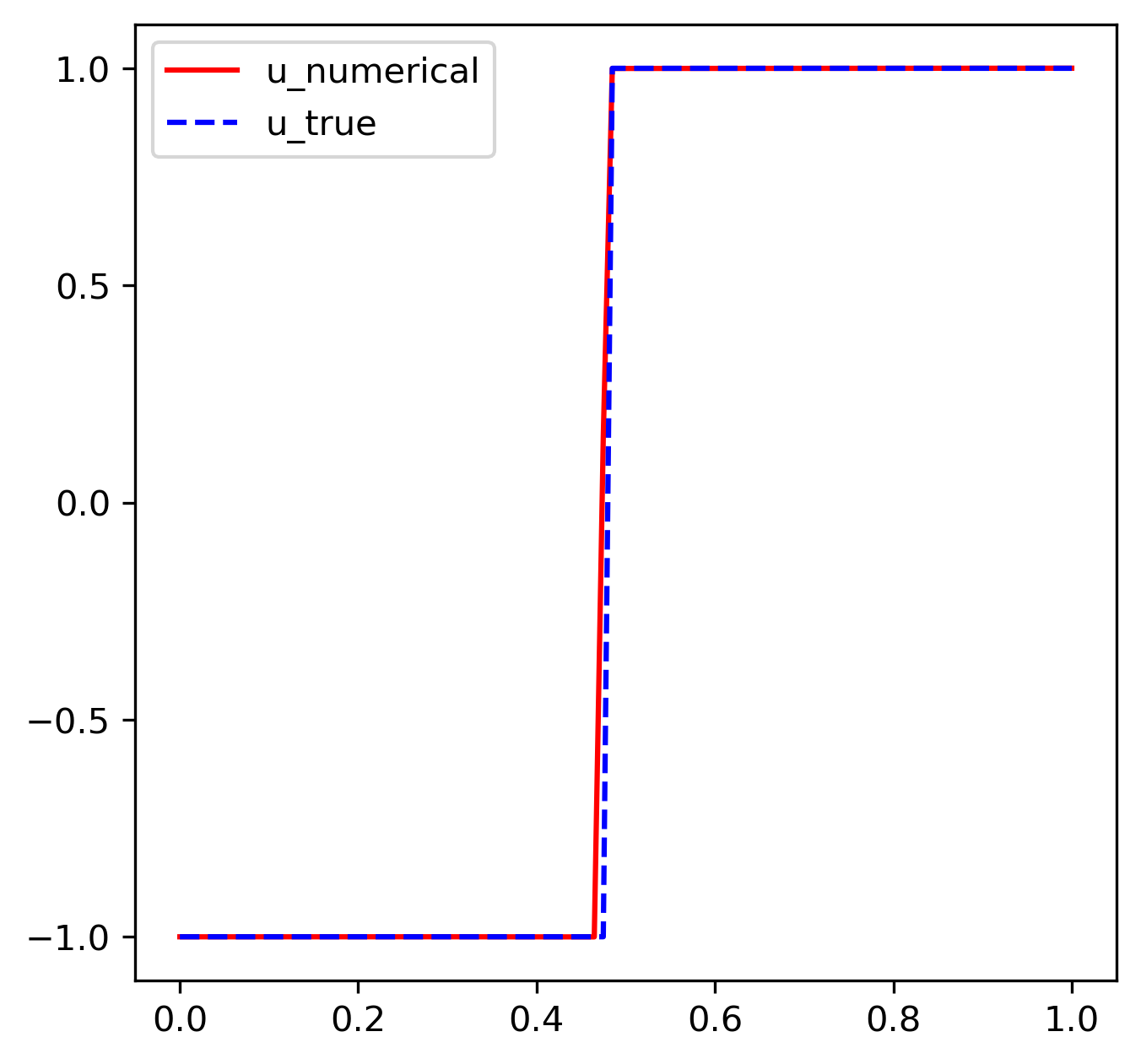}
\end{minipage}%
}%
\\
\subfigure[The breaking hyperplanes of the approximation in Figure \ref{comparison3_one}\label{breaking3_one}]{
\begin{minipage}[t]{0.4\linewidth}
\centering
\includegraphics[width=1.8in]{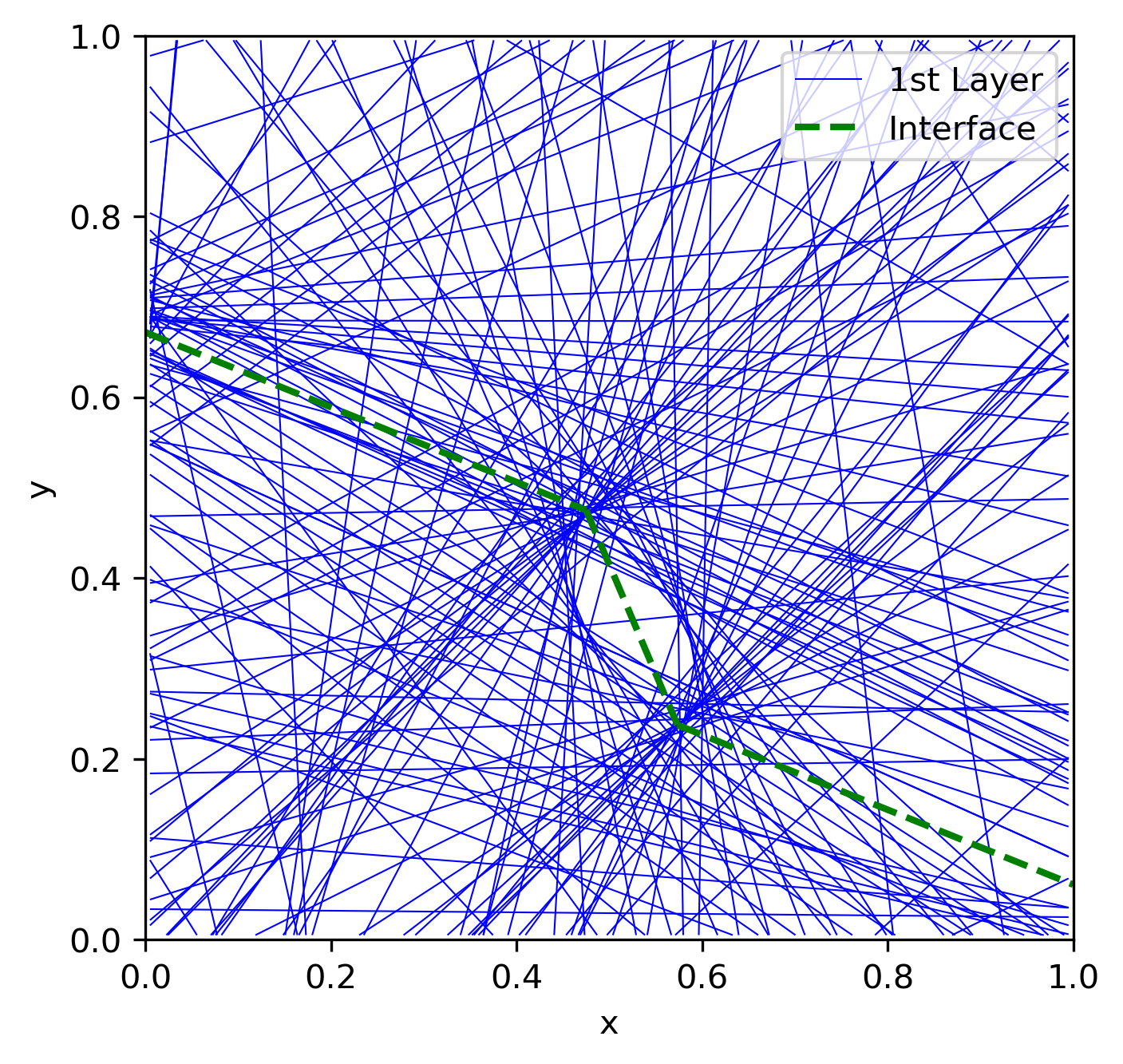}
\end{minipage}%
}%
\hspace{0.2in}
\subfigure[The breaking hyperplanes of the approximation in Figure \ref{comparison32}\label{breaking3}]{
\begin{minipage}[t]{0.4\linewidth}
\centering
\includegraphics[width=1.8in]{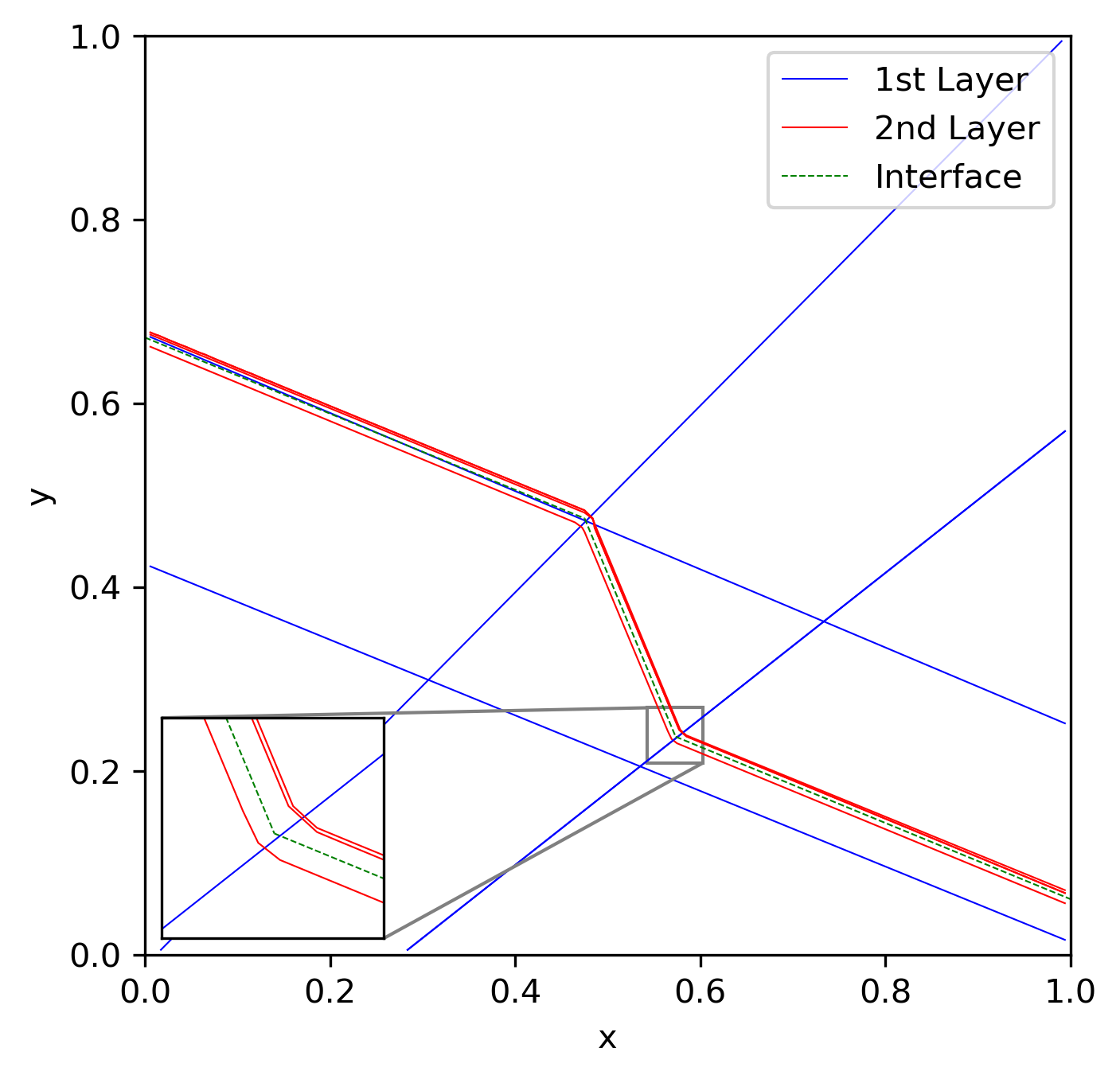}
\end{minipage}%
}%
\caption{Approximation results of the problem in \Cref{2d test1 section}}
\label{2d test1}
\end{figure}

This example is a modification of one from \cite{Cai2021linear}. Let $\bar{\Omega}=\bar{\Upsilon}_1\cup \bar{\Upsilon}_2 \cup \bar{\Upsilon}_3$ and 
\[
    \Upsilon_1=\{(x,y)\in\Omega:y\ge x\},\,\,
    \Upsilon_2=\{(x,y)\in\Omega: x-\tfrac{a}{2}\le y< x\},\,\mbox{ and }\,
    \Upsilon_3=\{(x,y)\in\Omega:y<x-\tfrac{a}{2}\}
\]
with $a=43/64$. The advective velocity field is a piecewise constant field given by \begin{equation}\label{2d test1 beta}
\bm{\beta}(x,y) =\left\{ \begin{array}{rl}
(-1,\sqrt{2}-1)^T,&(x,y)\in\Upsilon_1,\\[2mm]
(1-\sqrt{2},1)^T,&(x,y)\in\Upsilon_2,\\[2mm]
(-1,\sqrt{2}-1)^T,&(x,y)\in\Upsilon_3.
\end{array}\right.
\end{equation}
The inflow boundary and the inflow boundary condition are given by
\begin{eqnarray*}
\Gamma_{-}&=&\{(1,y):y\in(0,1)\}\cup\{(x,0):x\in(0,1)\}\\[2mm]
\mbox{and }\,\, g(x,y)&=&\left\{ \begin{array}{rl}
 1,& (x,y)\in \Gamma^1_-\equiv \{(1,y): y\in[1-\sqrt{2}+\frac{\sqrt{2}}{2}a,1)\}, \\[2mm]
 -1, &(x,y)\in \Gamma^2_-=\Gamma_-\setminus \Gamma_-^1,
\end{array}\right.
\end{eqnarray*} 
respectively. Let 
\begin{eqnarray*}
   && \hat{\Upsilon}_1=\{(x,y)\in\Upsilon_1: y<(1-\sqrt{2})x+a\},\,\,
   \hat{\Upsilon}_2=\{(x,y)\in\Upsilon_2: y<\tfrac{1}{1-\sqrt{2}}(x-\tfrac{a}{\sqrt{2}})+\tfrac{a}{\sqrt{2}}\},\\[2mm]
   \mbox{and }\,\, &&  \hat{\Upsilon}_3=\{(x,y)\in\Upsilon_3:y<({1-\sqrt{2}})x+\tfrac{\sqrt{2}}{2}a\}.
\end{eqnarray*}
The following right-hand side function is (see \cref{comparison31})
\begin{equation}
f(x,y)=\left\{ \begin{array}{rl}
 -1,& (x,y)\in \Omega_1\equiv\hat{\Upsilon}_1\cup \hat{\Upsilon}_2\cup \hat{\Upsilon}_3, \\[2mm]
 1, & (x,y)\in\Omega_2=\Omega\setminus\Omega_1.
\end{array}\right.
\end{equation}

200000 iterations were implemented with 2--300--1 and 2--5--5--1 ReLU NN functions. The numerical results are presented in \cref{2d test1,2d test1 table}. The numerical errors (\cref{2d test1 table}), trace (\cref{vertical3_one}), and approximation graph (\cref{comparison3_one}) of the 2-layer ReLU NN function approximation imply that the 2-layer network structure failed to approximate the solution (\cref{comparison31}) especially around the discontinuity interface (\cref{interface3}), although the breaking hyperplanes (\cref{breaking3_one}) indicate that the approximation roughly formed the transition layer around the interface. This and the remaining examples suggest the 2-layer network structure may not be able to approximate discontinuous solutions well as we expected in \cref{uat}. On the other hand, the 3-layer ReLU NN function approximation with the 2-5-5-1 structure with 4\% of the number of parameters of the 2-layer one approximates the solution accurately. Again, this and the remaining examples suggest that 3-layer ReLU NN functions may be more efficient than 2-layer ones of even bigger sizes. In this example, because of the shape of the interface and $\hat{u}=0$, the CPWL function $p$ with small $\varepsilon$, which we constructed in \cref{chi-curve lem} is expected to be a good approximation of the solution, and \cref{comparison32} indicates that the approximation in $\cM(3,10)$ is indeed such a function. The second-layer breaking hyperplanes (\cref{breaking3}) also help us to verify that a sharp transition layer was generated along the discontinuity interface, which is again consistent with our convergence analysis. The trace (\cref{vertical3}) of the 3-layer ReLU NN function approximation exhibits no oscillation.

\begin{table}[htbp]\label{2d test1 table}
\caption{Relative errors of the problem in \Cref{2d test1 section}}
\centering
\begin{tabular}{|l|l|l|l|l|}
\hline
Network structure  &$\frac{\|u-{u}^{_N}_{_\cT}\|_0}{\|u\|_0}$ &$\frac{\vertiii{u-{u}^{_N}_{_\cT}}_{\bm\beta}}{\vertiii{u}_{\bm\beta}}$ & $\frac{\mathcal{L}^{1/2}({u}^{_N}_{_\cT},\bf f)}{\mathcal{L}^{1/2}({u}^{_N}_{_\cT},\bf 0)}$ & Parameters \\ \hline
2--300--1 & 0.279867 & 0.404376 & 0.300774 & 1201\\ \hline
2--5--5--1  & 0.074153 & 0.079193 & 0.044987   & 51\\ \hline
\end{tabular}
\end{table}

\subsubsection{A problem with a 4-line segment interface}\label{2d test2 section}
\begin{figure}[htbp]\label{2d test2}
\centering
\subfigure[The interface\label{interface4}]{
\begin{minipage}[t]{0.4\linewidth}
\centering
\includegraphics[width=1.8in]{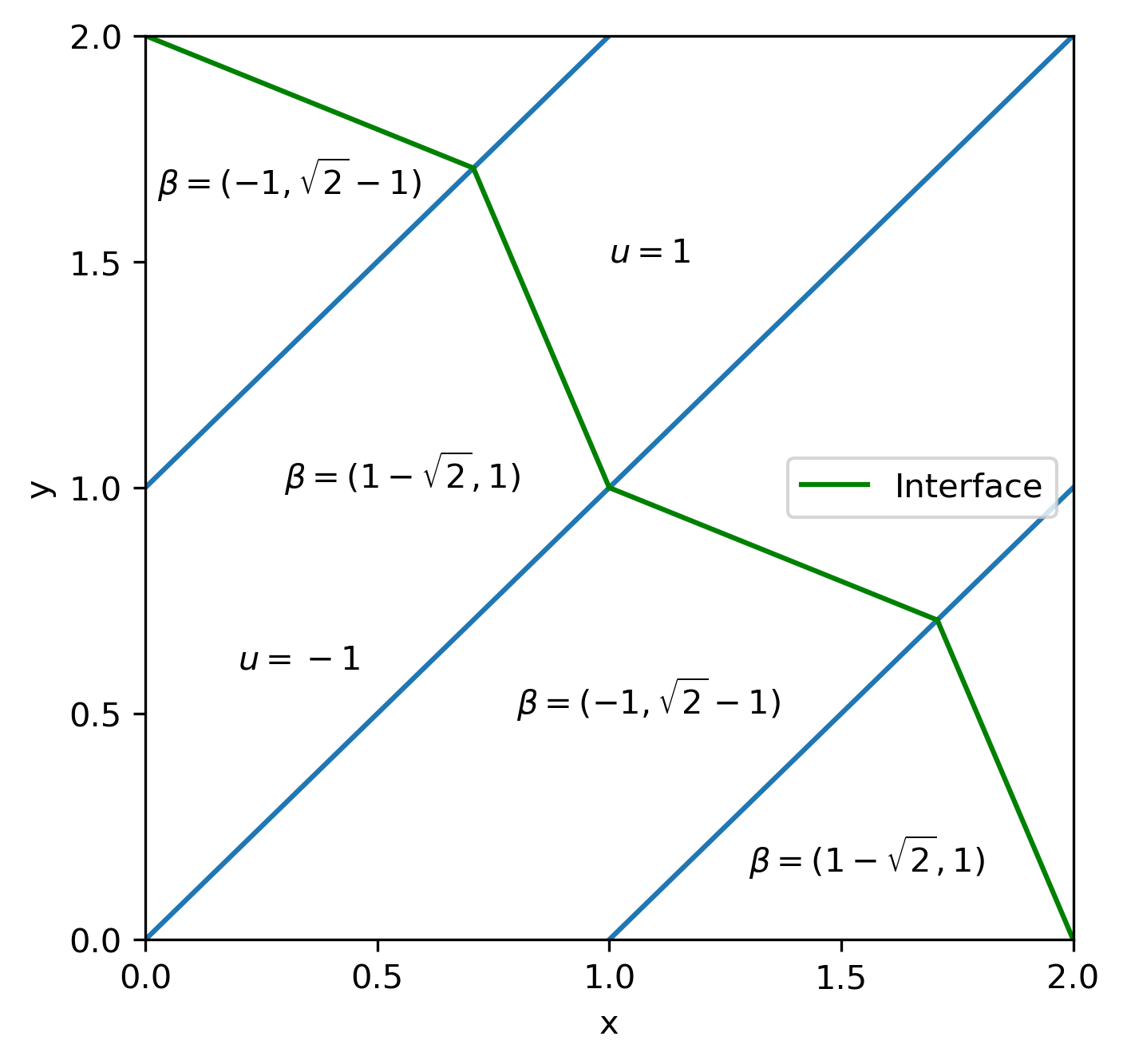}
\end{minipage}%
}%
\hspace{0.2in}
\subfigure[The exact solution\label{comparison41}]{
\begin{minipage}[t]{0.4\linewidth}
\centering
\includegraphics[width=1.8in]{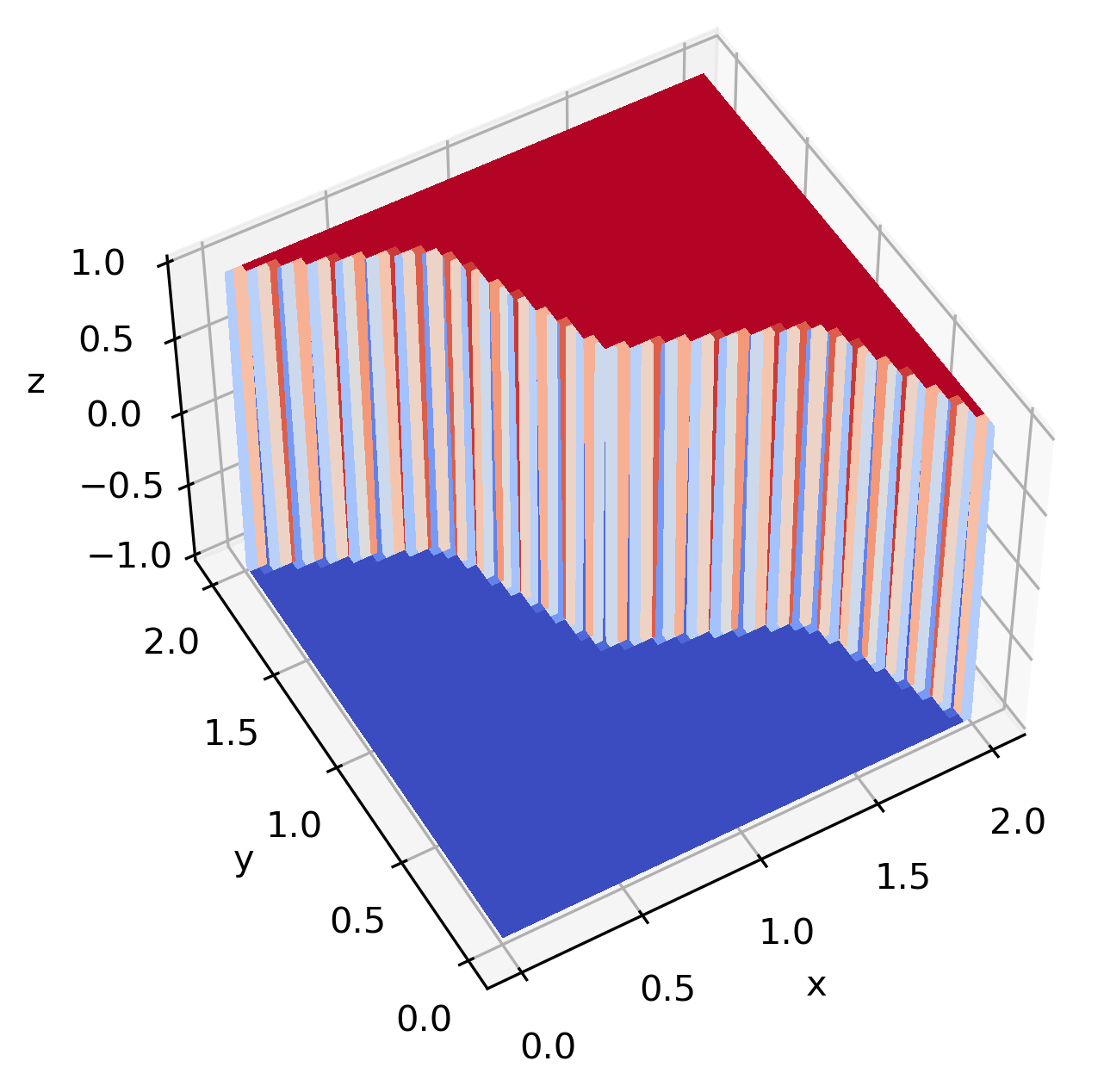}
\end{minipage}%
}%
\\
\subfigure[A 2--300--1 ReLU NN function approximation\label{comparison4_one}]{
\begin{minipage}[t]{0.4\linewidth}
\centering
\includegraphics[width=1.8in]{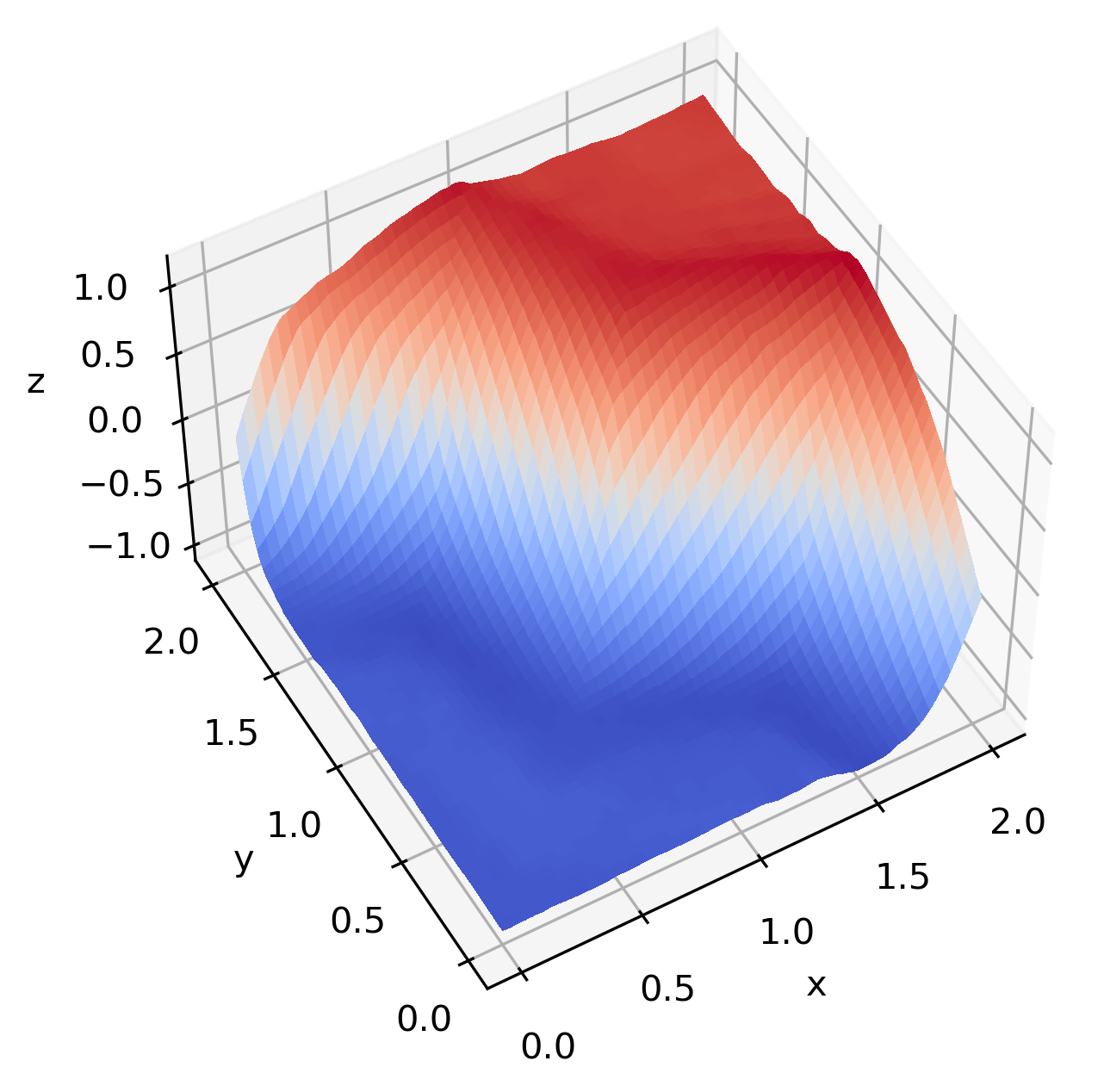}
\end{minipage}%
}%
\hspace{0.2in}
\subfigure[A 2--6--6--1 ReLU NN function approximation\label{comparison42}]{
\begin{minipage}[t]{0.4\linewidth}
\centering
\includegraphics[width=1.8in]{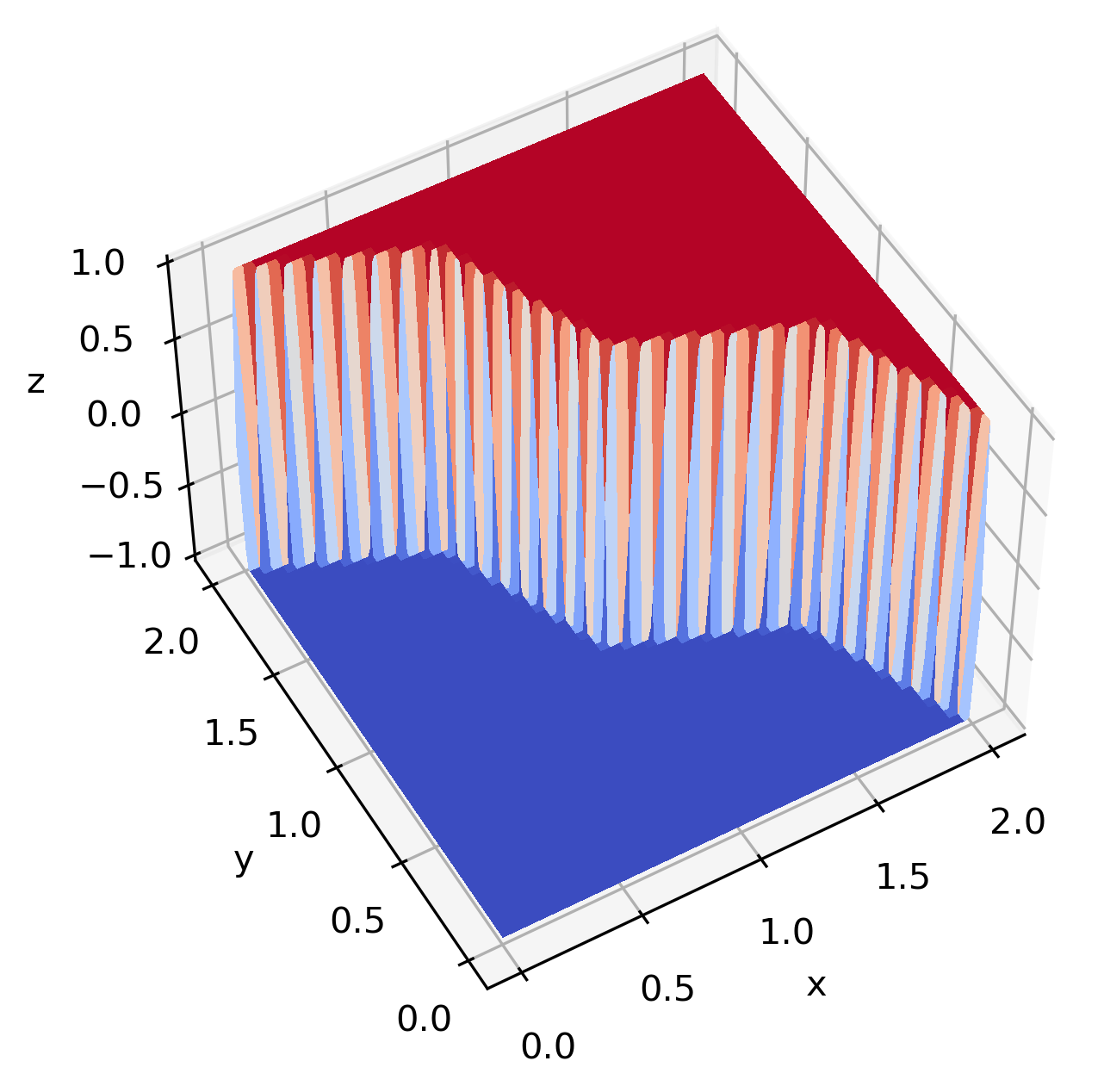}
\end{minipage}%
}%
\\
\subfigure[The trace of Figure \ref{comparison4_one} on $y=x$\label{vertical4_one}]{
\begin{minipage}[t]{0.4\linewidth}
\centering
\includegraphics[width=1.8in]{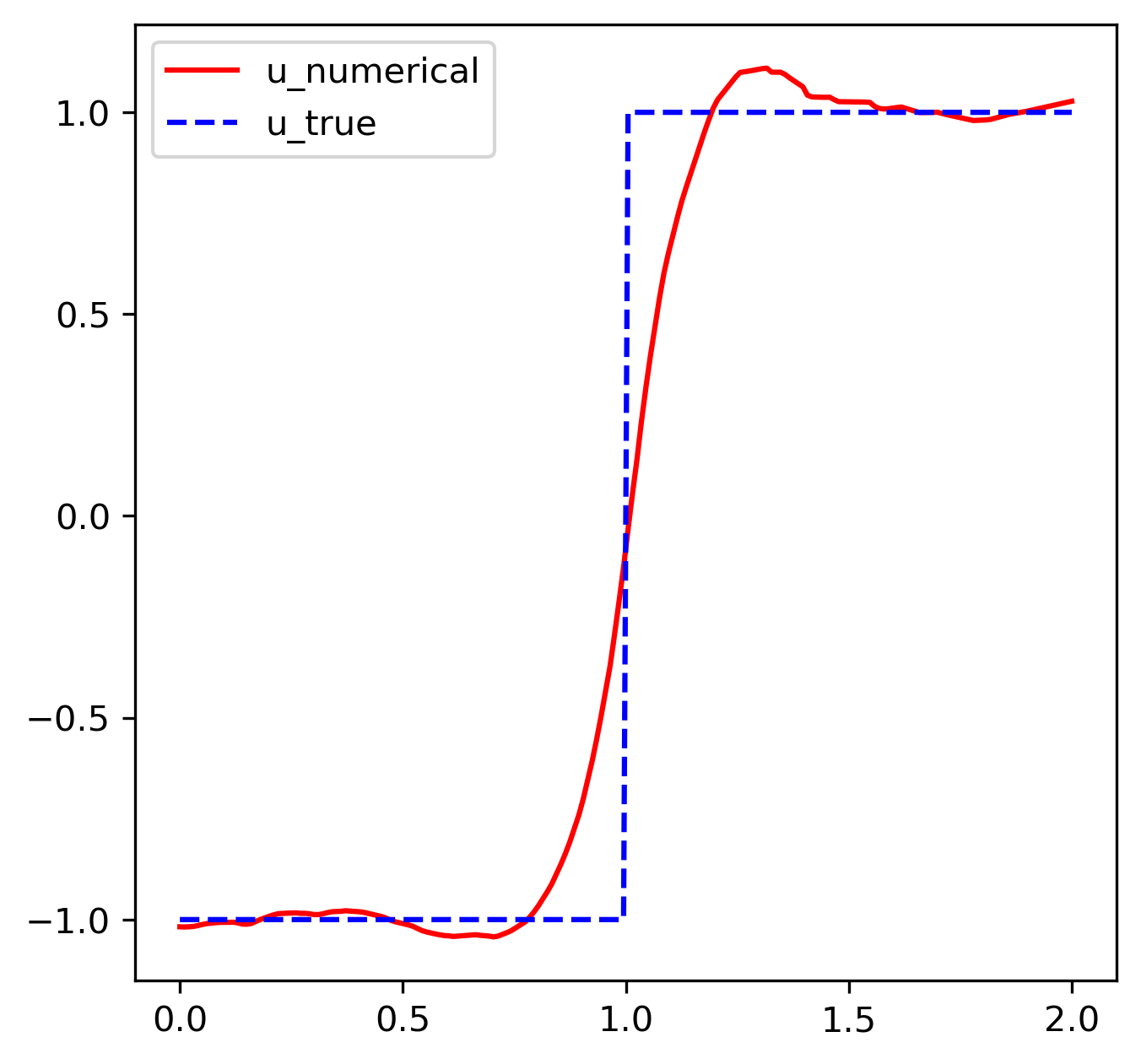}
\end{minipage}%
}%
\hspace{0.2in}
\subfigure[The trace of Figure \ref{comparison42} on $y=x$\label{vertical4}]{
\begin{minipage}[t]{0.4\linewidth}
\centering
\includegraphics[width=1.8in]{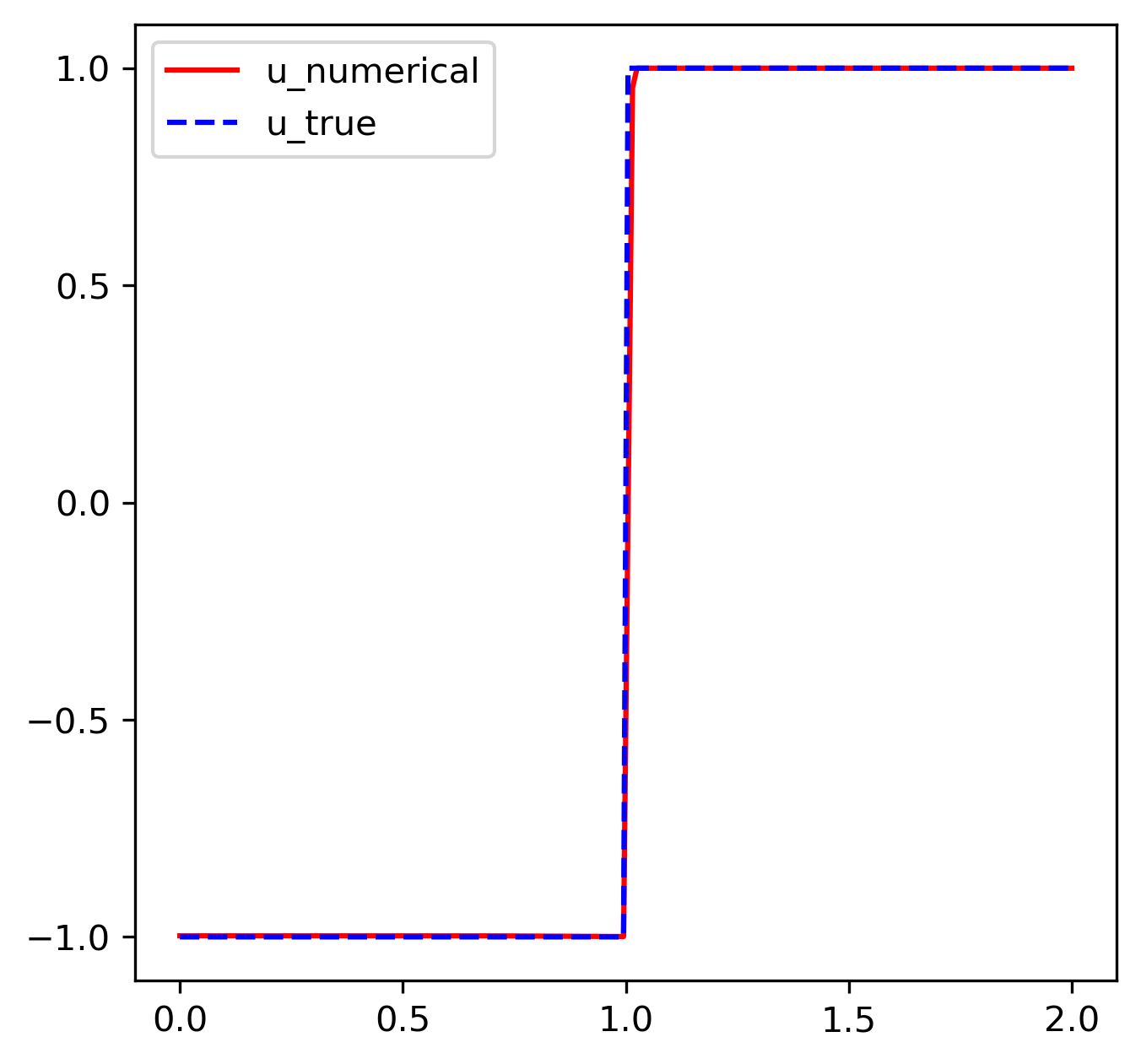}
\end{minipage}%
}%
\\
\subfigure[The breaking hyperplanes of the approximation in Figure \ref{comparison4_one}\label{breaking4_one}]{
\begin{minipage}[t]{0.4\linewidth}
\centering
\includegraphics[width=1.8in]{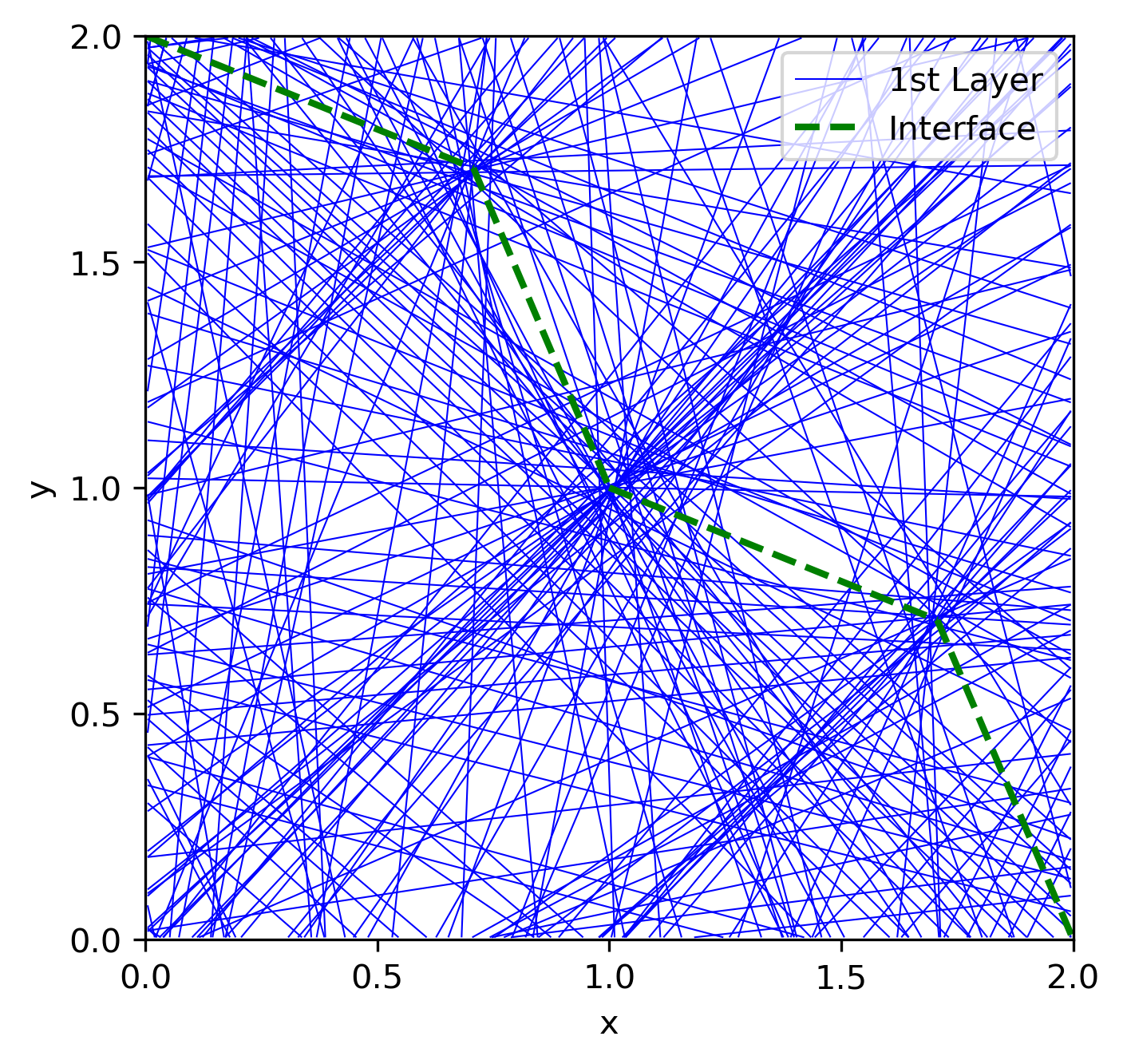}
\end{minipage}%
}%
\hspace{0.2in}
\subfigure[The breaking hyperplanes of the approximation in Figure \ref{comparison42}\label{breaking4}]{
\begin{minipage}[t]{0.4\linewidth}
\centering
\includegraphics[width=1.8in]{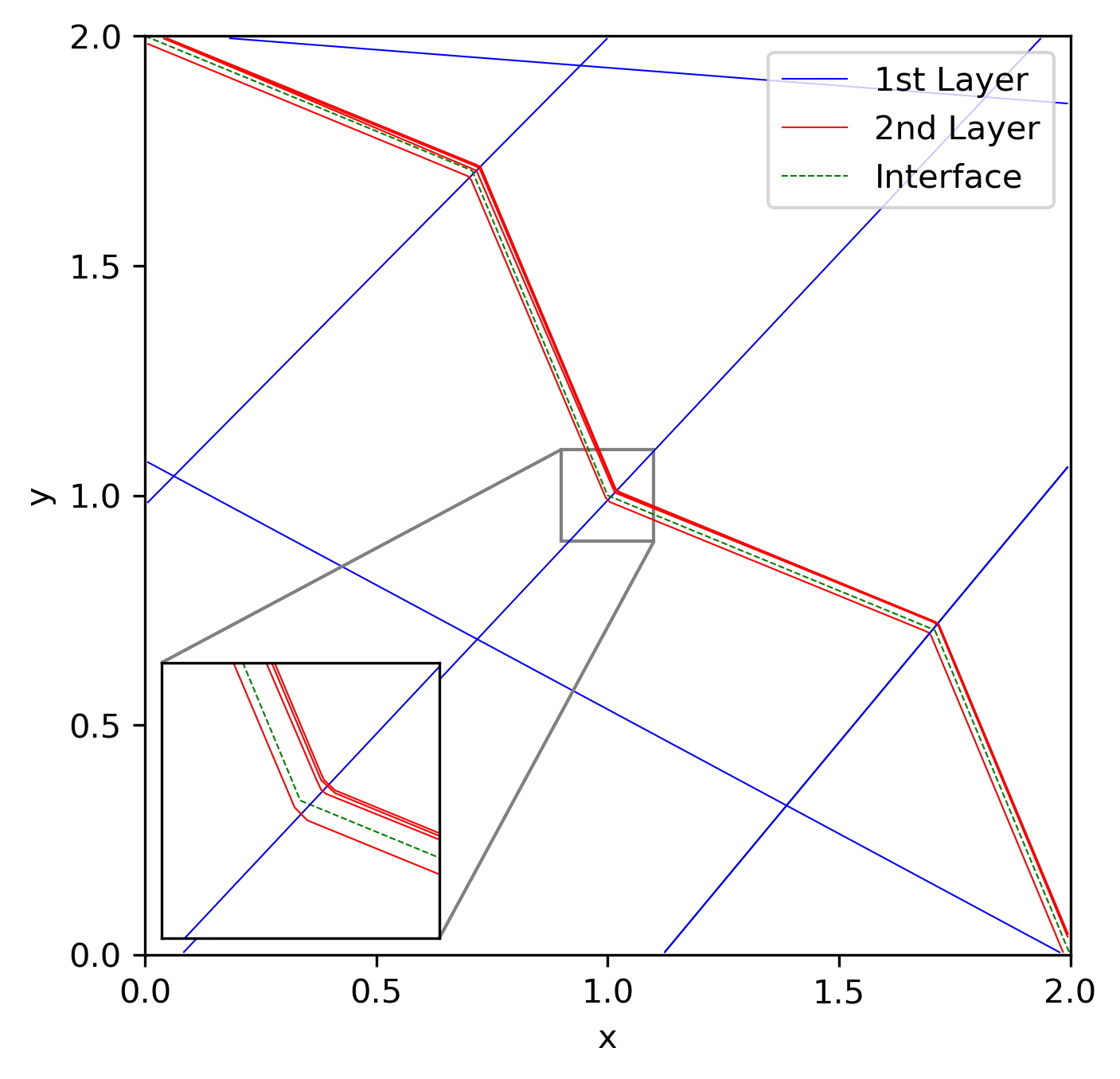}
\end{minipage}%
}%
\caption{Approximation results of the problem in \Cref{2d test2 section}}
\end{figure}

Let $\bar{\Omega}=\bar{\Upsilon}_1\cup \bar{\Upsilon}_2 \cup \bar{\Upsilon}_3\cup \bar{\Upsilon}_4$ and
\begin{multline*}
\Upsilon_1=\{(x,y)\in \Omega :\ y\ge x+1\},\,\, \Upsilon_2=\{(x,y)\in \Omega,\ x\le y<x+1\},\\[2mm] \Upsilon_3=\{(x,y)\in \Omega,\  x-1\le y<x\},\,\, \text{and } \Upsilon_4=\{(x,y)\in \Omega,\ y<x-1\}.
\end{multline*}
The advective velocity field is a piecewise constant field given by \begin{equation}
\bm{\beta}(x,y) =\left\{ \begin{array}{rl}
(-1,\sqrt{2}-1)^T,&(x,y)\in\Upsilon_1,\\[2mm]
(1-\sqrt{2},1)^T,&(x,y)\in\Upsilon_2,\\[2mm]
(-1,\sqrt{2}-1)^T,&(x,y)\in \Upsilon_3,\\[2mm]
(1-\sqrt{2},1)^T,&(x,y)\in\Upsilon_4.
\end{array}\right.
\end{equation}
The inflow boundary and the inflow boundary condition are given by
\begin{eqnarray*}
\Gamma_{-}&=&\{(2,y):y\in(0,2)\}\cup\{(x,0):x\in(0,2)\}\\[2mm]
\mbox{and }\,\, g(x,y)&=&\left\{ \begin{array}{rl}
 1,& (x,y)\in \Gamma^1_-\equiv \{(2,y): y\in(0,2)\}, \\[2mm]
 -1, & (x,y)\in \Gamma^2_-=\Gamma_-\setminus \Gamma_-^1,
\end{array}\right.
\end{eqnarray*} 
respectively. Let
\begin{multline*}
\hat{\Upsilon}_1=\{(x,y)\in\Upsilon_1: y<(1-\sqrt{2})x+2\},\,\, \hat{\Upsilon}_2=\{(x,y)\in\Upsilon_2: y<\tfrac{1}{1-\sqrt{2}}(x-1)+1\},\\[2mm]
\hat{\Upsilon}_3=\{(x,y)\in\Upsilon_3:y<({1-\sqrt{2}})(x-1)+1\},\,\,\text{and }\hat{\Upsilon}_4=\{(x,y)\in\Upsilon_3:y<\tfrac{1}{1-\sqrt{2}}x+\tfrac{2}{\sqrt{2}-1}\}.
\end{multline*}
The following right-hand side function is (see \cref{comparison41})
\begin{equation}
f(x,y)=\left\{ \begin{array}{rl}
 -1,& (x,y)\in \Omega_1\equiv\hat{\Upsilon}_1\cup \hat{\Upsilon}_2\cup \hat{\Upsilon}_3\cup\hat{\Upsilon}_4, \\[2mm]
 1, & (x,y)\in\Omega_2=\Omega\setminus\Omega_1.
\end{array}\right.
\end{equation}

200000 iterations were implemented with 2--300--1 and 2--6--6--1 ReLU NN functions. The numerical results are presented in \cref{2d test2,2d test2 table}. Since the interface has one more line segment than that of Example \ref{2d test1 section}, we increased the number of hidden neurons to have higher expresiveness. The 2-6-6-1 structure with 6\% of the number of parameters of the 2-300-1 structure approximated the solution (\cref{comparison41}) accurately and \cref{comparison42} indicates that the approximation in $\mathcal{M}(3,12)$ is the CPWL function $p$ with small $\varepsilon$ in \cref{chi-curve lem}. The trace (\cref{vertical4}) shows no oscillation and the second-layer breaking hyperplanes (\cref{breaking4}) along the discontinuity interface (\cref{interface4}) show where a sharp transition layer was generated. On the other hand, the 2-300-1 ReLU NN function approximation roughly found the location of the interface (\cref{breaking4_one}) but did not approximate the solution well (\cref{comparison4_one,vertical4_one,2d test2 table}).

\begin{table}[htbp]\label{2d test2 table}
\caption{Relative errors of the problem in \Cref{2d test2 section}}
\centering
\begin{tabular}{|l|l|l|l|l|}
\hline
Network structure  &$\frac{\|u-{u}^{_N}_{_\cT}\|_0}{\|u\|_0}$ &$\frac{\vertiii{u-{u}^{_N}_{_\cT}}_{\bm\beta}}{\vertiii{u}_{\bm\beta}}$ & $\frac{\mathcal{L}^{1/2}({u}^{_N}_{_\cT},\bf f)}{\mathcal{L}^{1/2}({u}^{_N}_{_\cT},\bf 0)}$ & Parameters \\ \hline
2--300--1 & 0.288282 & 0.358756 & 0.306695 & 1201\\ \hline
2--6--6--1  & 0.085817 & 0.091800 & 0.069808   & 67\\ \hline
\end{tabular}
\end{table}

\subsubsection{A problem with a curved interface}\label{2d test3 section}
\begin{figure}[htbp]\label{2d test3}
\centering
\subfigure[The interface\label{interfacec}]{
\begin{minipage}[t]{0.4\linewidth}
\centering
\includegraphics[width=1.8in]{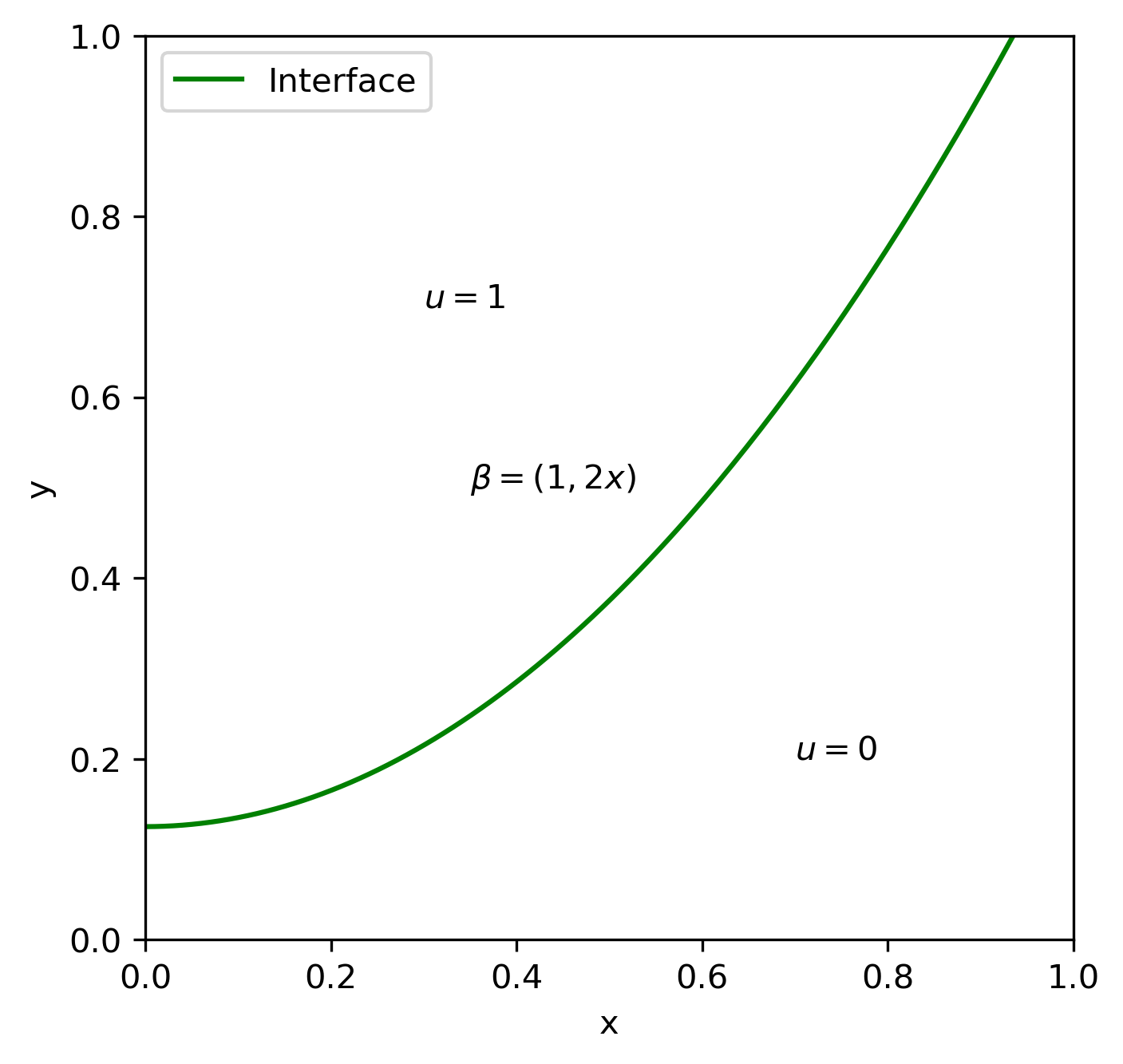}
\end{minipage}%
}%
\hspace{0.2in}
\subfigure[The exact solution\label{comparisonc1}]{
\begin{minipage}[t]{0.4\linewidth}
\centering
\includegraphics[width=1.8in]{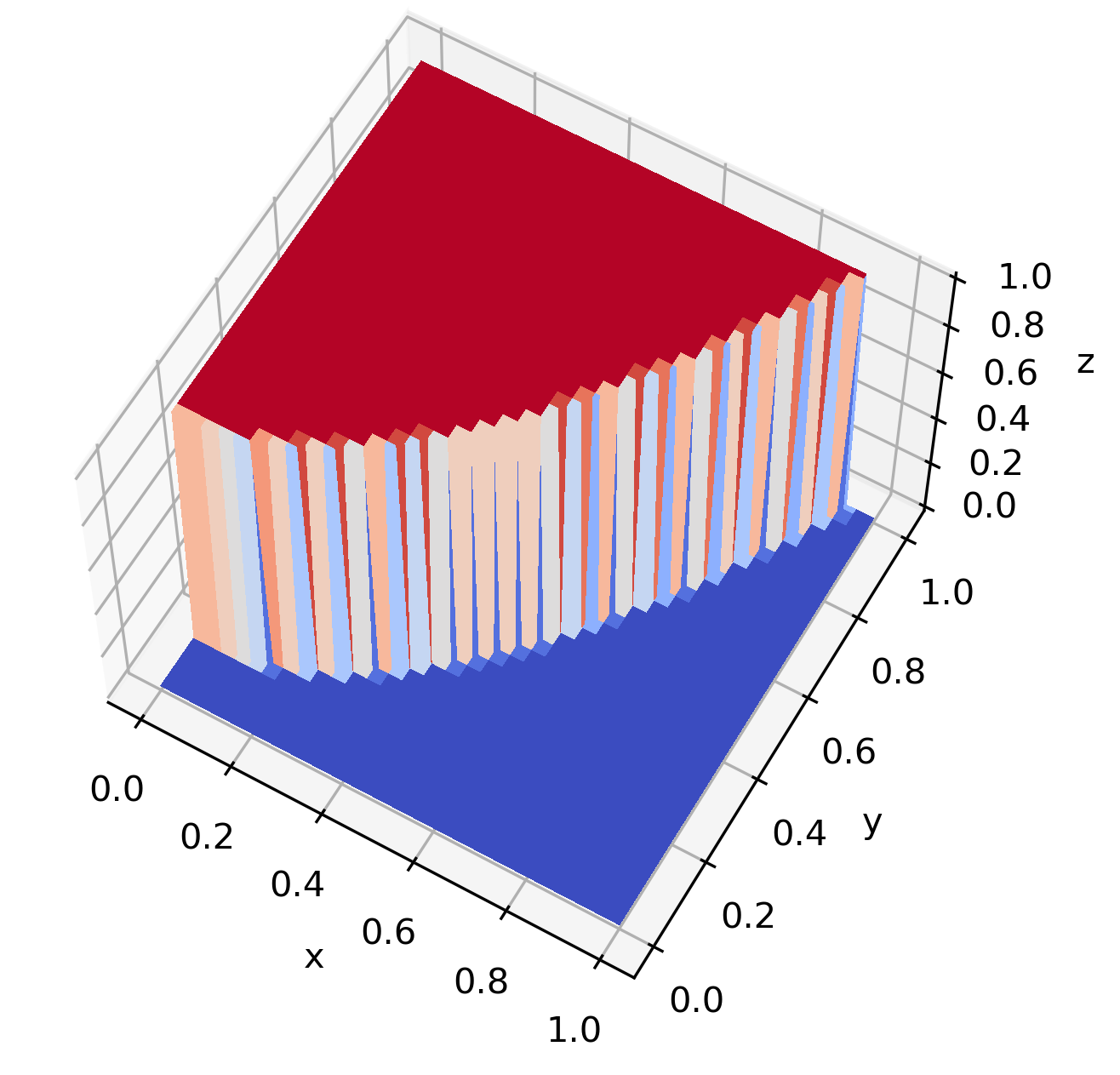}
\end{minipage}%
}%
\\
\subfigure[A 2--3000--1 ReLU NN function approximation\label{comparisonc_one}]{
\begin{minipage}[t]{0.4\linewidth}
\centering
\includegraphics[width=1.8in]{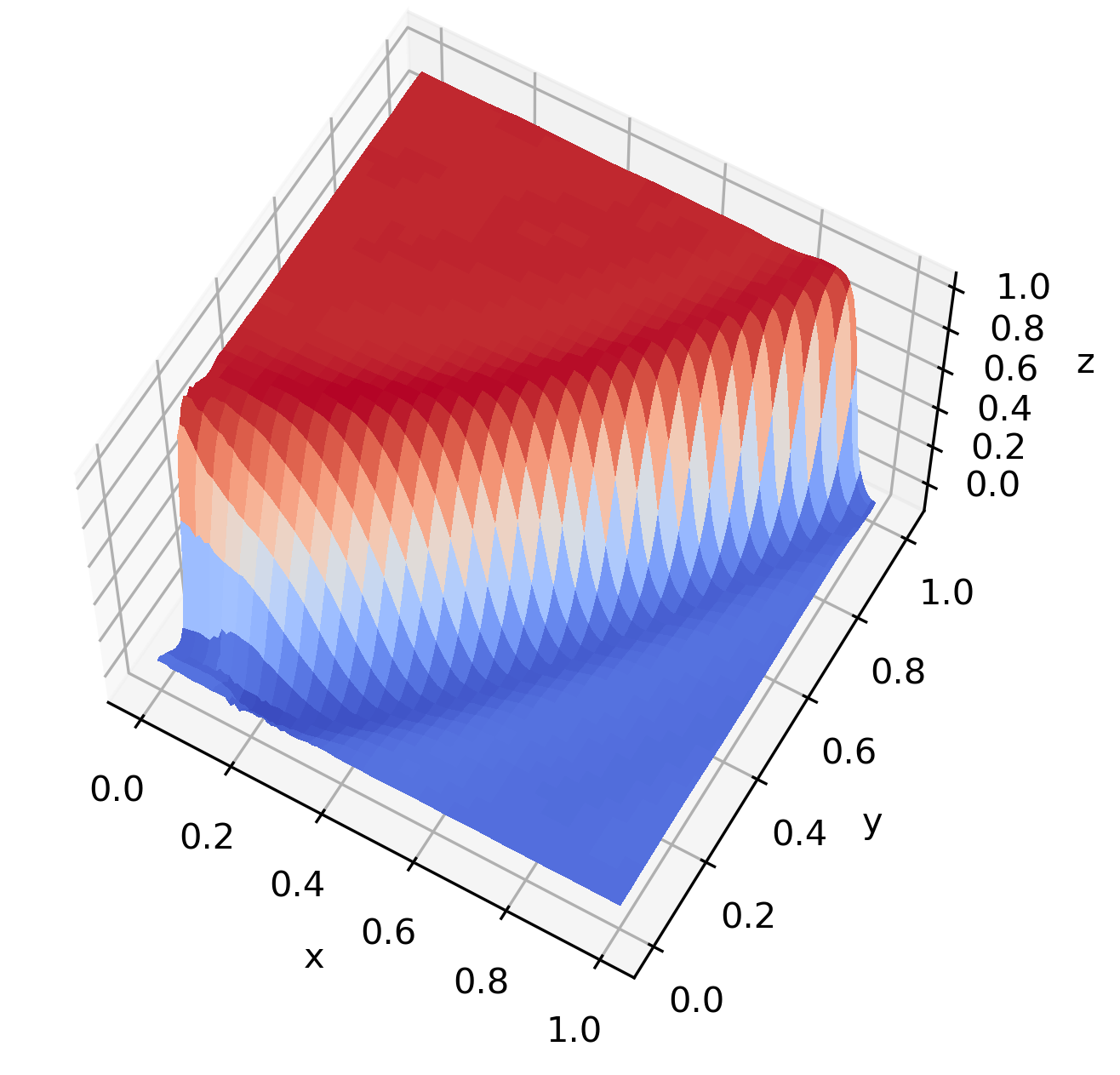}
\end{minipage}%
}%
\hspace{0.2in}
\subfigure[A 2--60--60--1 ReLU NN function approximation\label{comparisonc2}]{
\begin{minipage}[t]{0.4\linewidth}
\centering
\includegraphics[width=1.8in]{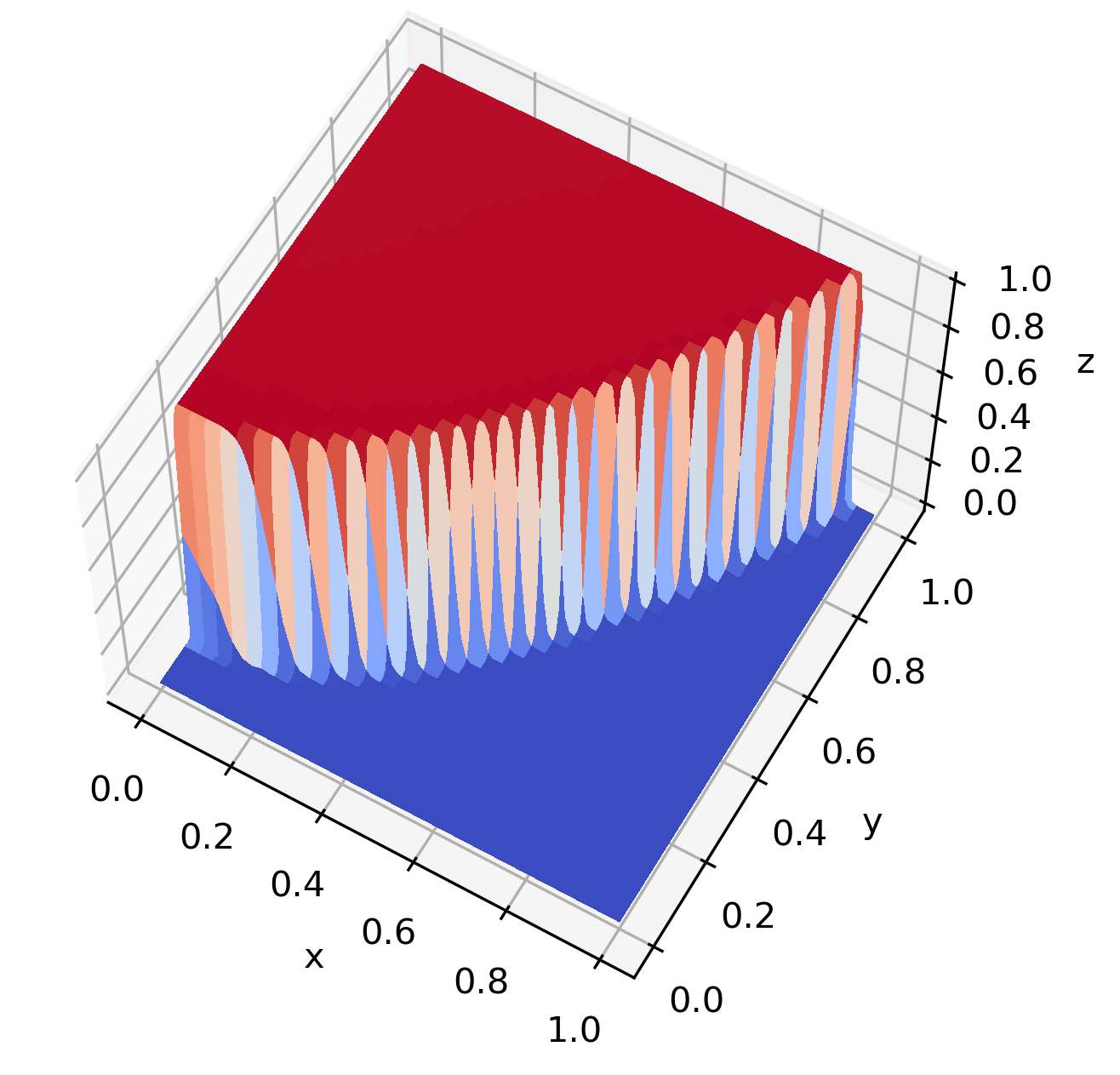}
\end{minipage}%
}%
\\
\subfigure[The trace of Figure \ref{comparisonc_one} on $y=1-x$\label{verticalc_one}]{
\begin{minipage}[t]{0.4\linewidth}
\centering
\includegraphics[width=1.8in]{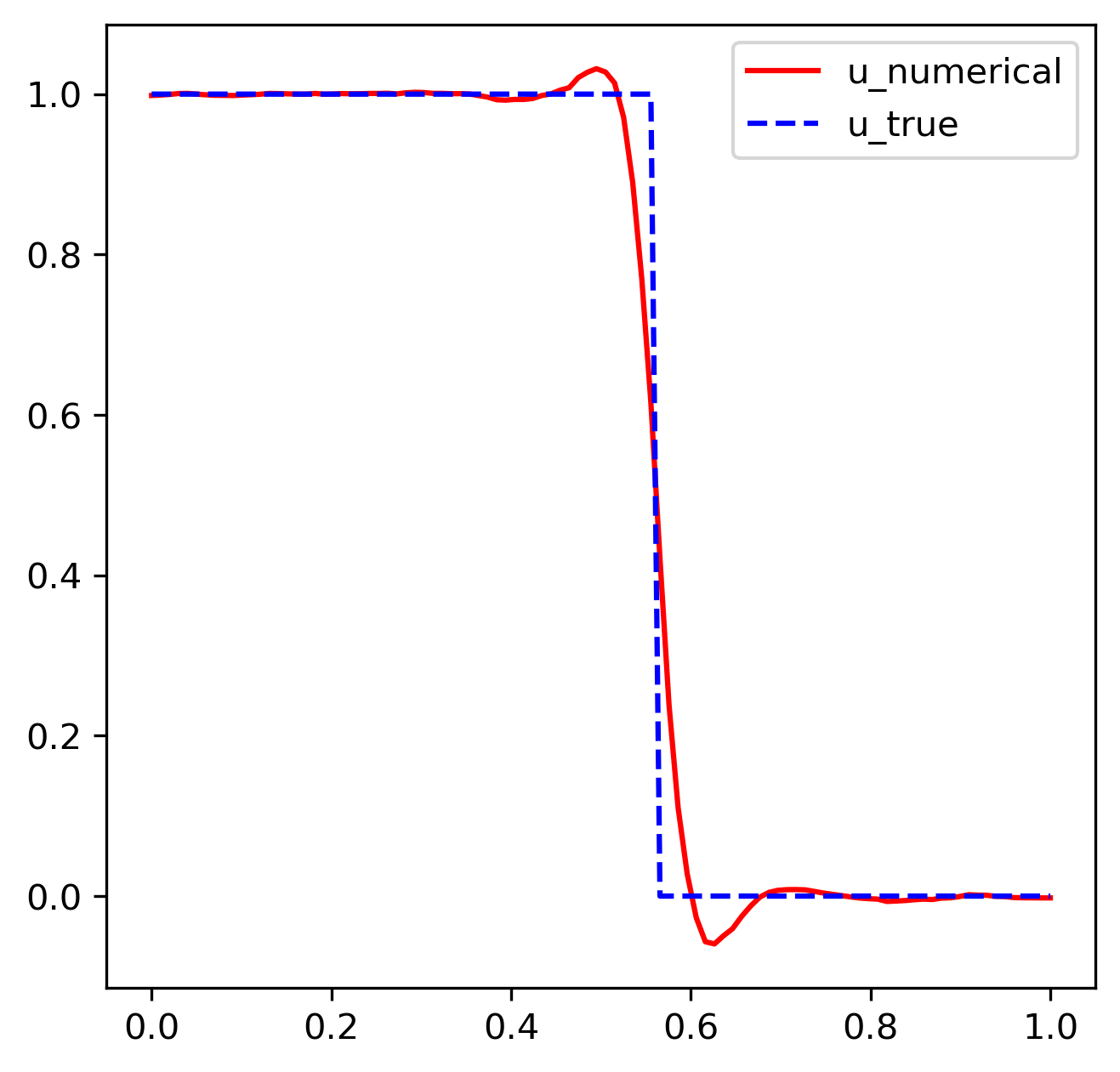}
\end{minipage}%
}%
\hspace{0.2in}
\subfigure[The trace of Figure \ref{comparisonc2} on $y=1-x$\label{verticalc}]{
\begin{minipage}[t]{0.4\linewidth}
\centering
\includegraphics[width=1.8in]{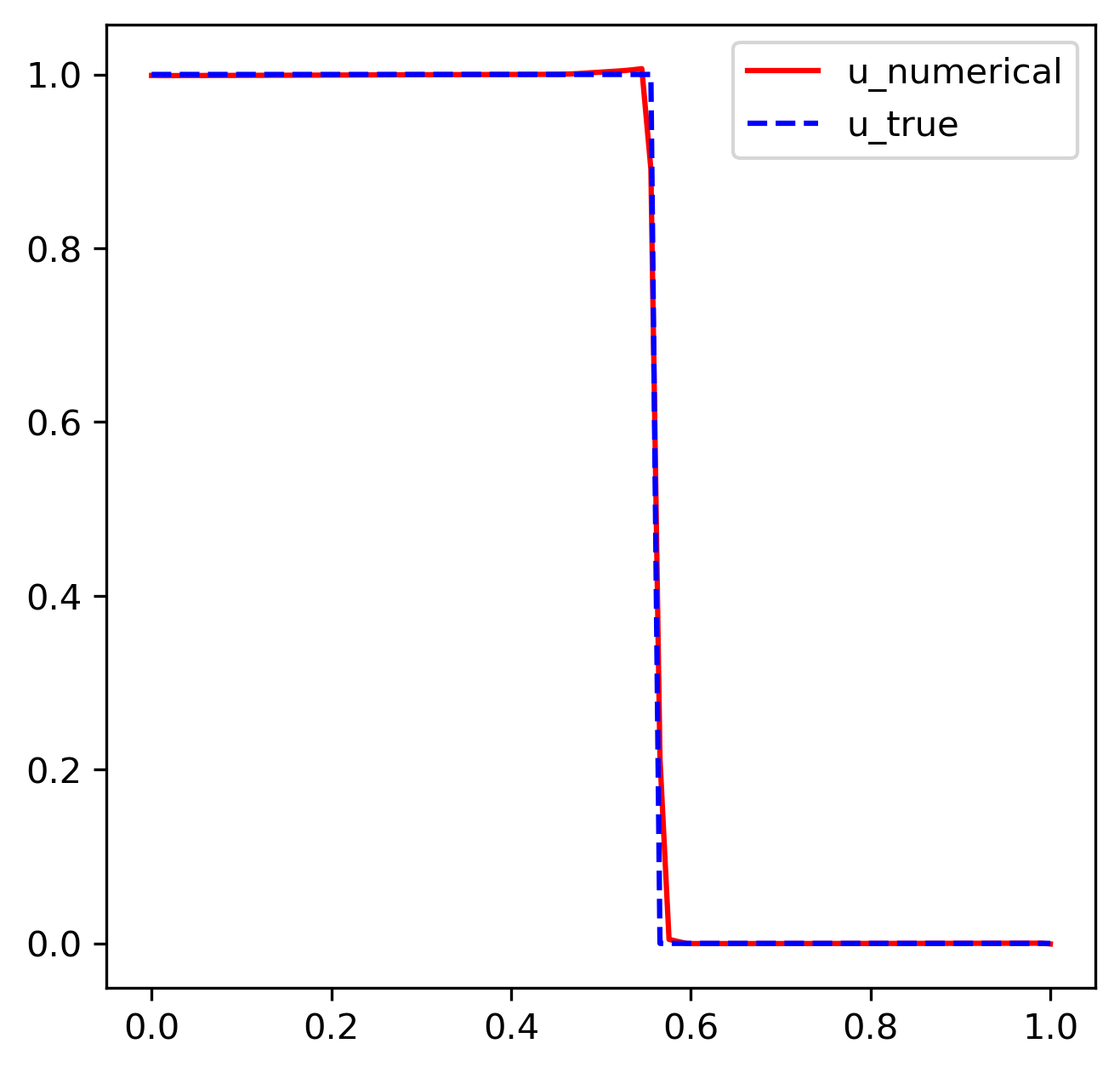}
\end{minipage}%
}%
\\
\subfigure[The breaking hyperplanes of the approximation in Figure \ref{comparisonc_one}\label{breakingc_one}]{
\begin{minipage}[t]{0.4\linewidth}
\centering
\includegraphics[width=1.8in]{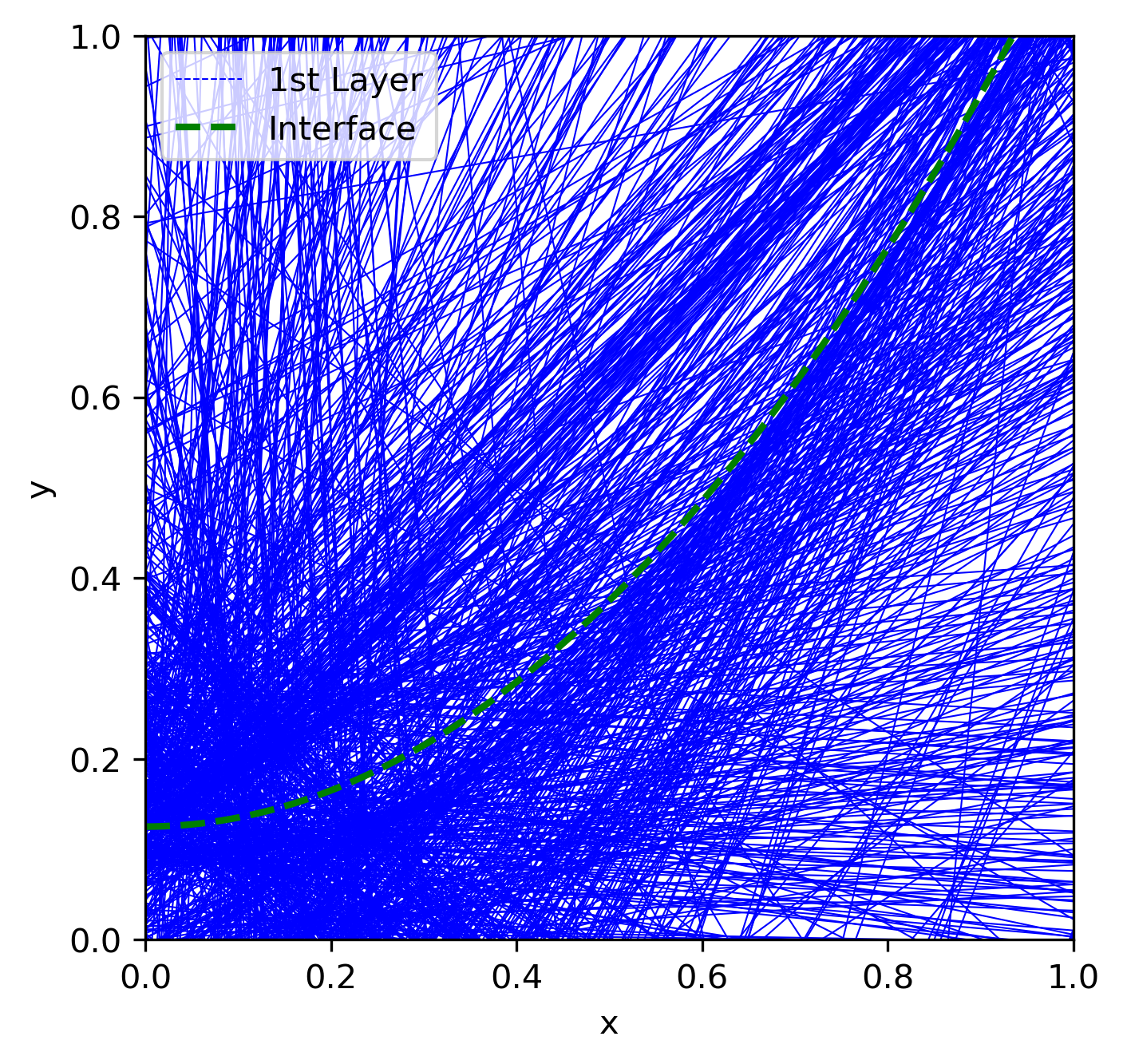}
\end{minipage}%
}%
\hspace{0.2in}
\subfigure[The breaking hyperplanes of the approximation in Figure \ref{comparisonc2}\label{breakingc}]{
\begin{minipage}[t]{0.4\linewidth}
\centering
\includegraphics[width=1.8in]{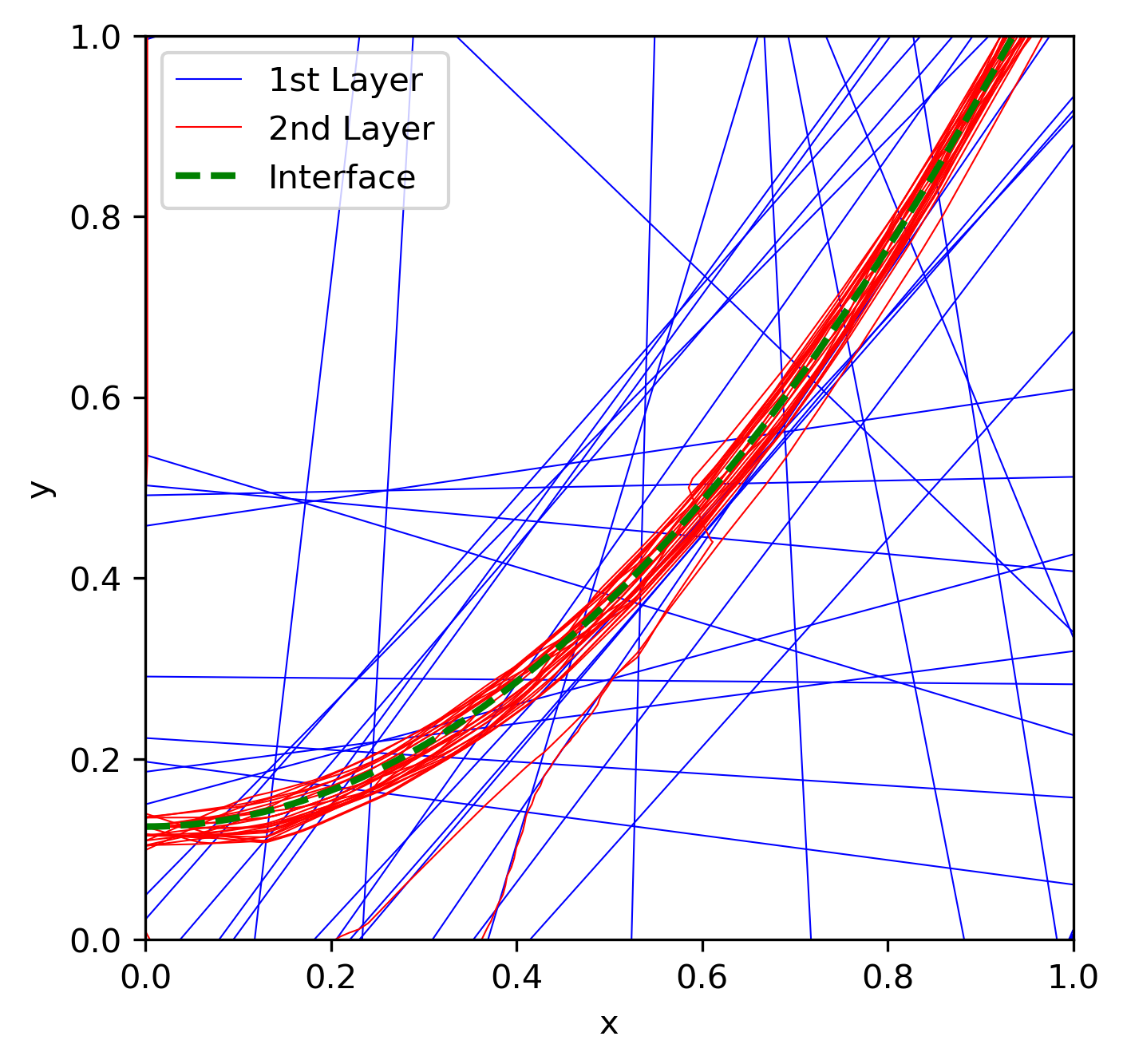}
\end{minipage}%
}%
\caption{Approximation results of the problem in \Cref{2d test3 section}}
\end{figure}

The advective velocity field is the variable field given by
\begin{equation}
\bm{\beta}(x,y) = (1,2x),\,\, (x,y)\in\Omega.
\end{equation}
The inflow boundary and the inflow boundary condition are given by
\begin{eqnarray*}
\Gamma_{-}&=&\{(0,y):y\in(0,1)\}\cup\{(x,0):x\in(0,1)\}\\[2mm]
\mbox{and }\,\, g(x,y)&=&\left\{ \begin{array}{rl}
 1,& (x,y)\in \Gamma^1_-\equiv\{(0,y): y\in[\frac{1}{8},1)\}, \\[2mm]
 0, & (x,y)\in \Gamma^2_-=\Gamma_-\setminus \Gamma_-^1,
\end{array}\right.
\end{eqnarray*} 
respectively. The following right-hand side function is (see \cref{comparisonc1})
\begin{equation}
f(x,y)=\left\{ \begin{array}{rl}
 0,& (x,y)\in \Omega_1\equiv\{(x,y)\in\Omega:y< x^2+\frac{1}{8}\}, \\[2mm]
 1, & (x,y)\in \Omega_2=\Omega\setminus\Omega_1.
\end{array}\right.
\end{equation}
300000 iterations were implemented with 2--3000--1 and 2--60--60--1 ReLU NN functions. The numerical results are presented in \cref{2d test3,2d test3 table}. We again increased the number of hidden neurons for the 3-layer network structure, assuming CPWL functions approximating a curved discontinuity interface (\cref{interfacec}) well would be a ReLU NN function with more hidden neurons. \cref{comparisonc_one,breakingc_one,verticalc_one} suggest that the 2-layer network structure failed to approximate the solution (\cref{comparisonc1}) around the discontinuity interface with more than three times the number of parameters of the 3-layer network structure. In contrast, the 3-layer network structure shows better numerical errors (\cref{2d test3 table}) and pointwise approximations (\cref{comparisonc2,verticalc}), locating the discontinuity interface (\cref{breakingc}).

\begin{table}[htbp]\label{2d test3 table}
\caption{Relative errors of the problem in \Cref{2d test3 section}}
\centering
\begin{tabular}{|l|l|l|l|l|}
\hline
Network structure  &$\frac{\|u-{u}^{_N}_{_\cT}\|_0}{\|u\|_0}$ &$\frac{\vertiii{u-{u}^{_N}_{_\cT}}_{\bm\beta}}{\vertiii{u}_{\bm\beta}}$ & $\frac{\mathcal{L}^{1/2}({u}^{_N}_{_\cT},\bf f)}{\mathcal{L}^{1/2}({u}^{_N}_{_\cT},\bf 0)}$ & Parameters \\ \hline
2--3000--1  & 0.134514 & 0.181499 & 0.078832   & 12001\\ \hline
2--60--60--1  & 0.066055 & 0.106095 & 0.030990   & 3901\\ \hline
\end{tabular}
\end{table}

\subsubsection{A problem with a curved interface and $\hat{u}\neq 0$ in \eqref{decop}}\label{2d test4 section}
\begin{figure}[htbp]\label{2d test4}
\centering
\subfigure[The interface\label{interfaceuhat}]{
\begin{minipage}[t]{0.4\linewidth}
\centering
\includegraphics[width=1.8in]{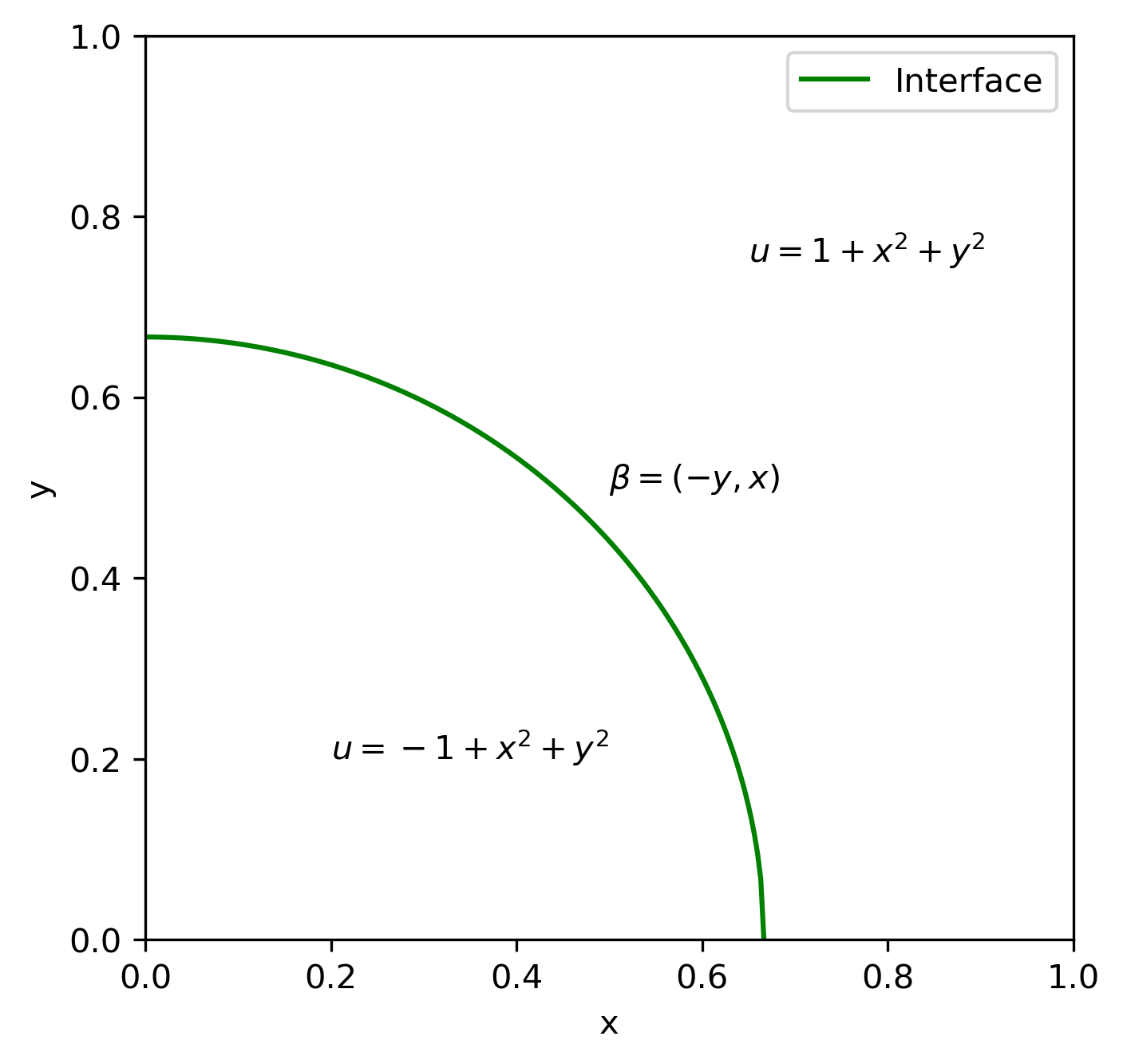}
\end{minipage}%
}%
\hspace{0.2in}
\subfigure[The exact solution\label{comparisonuhat1}]{
\begin{minipage}[t]{0.4\linewidth}
\centering
\includegraphics[width=1.8in]{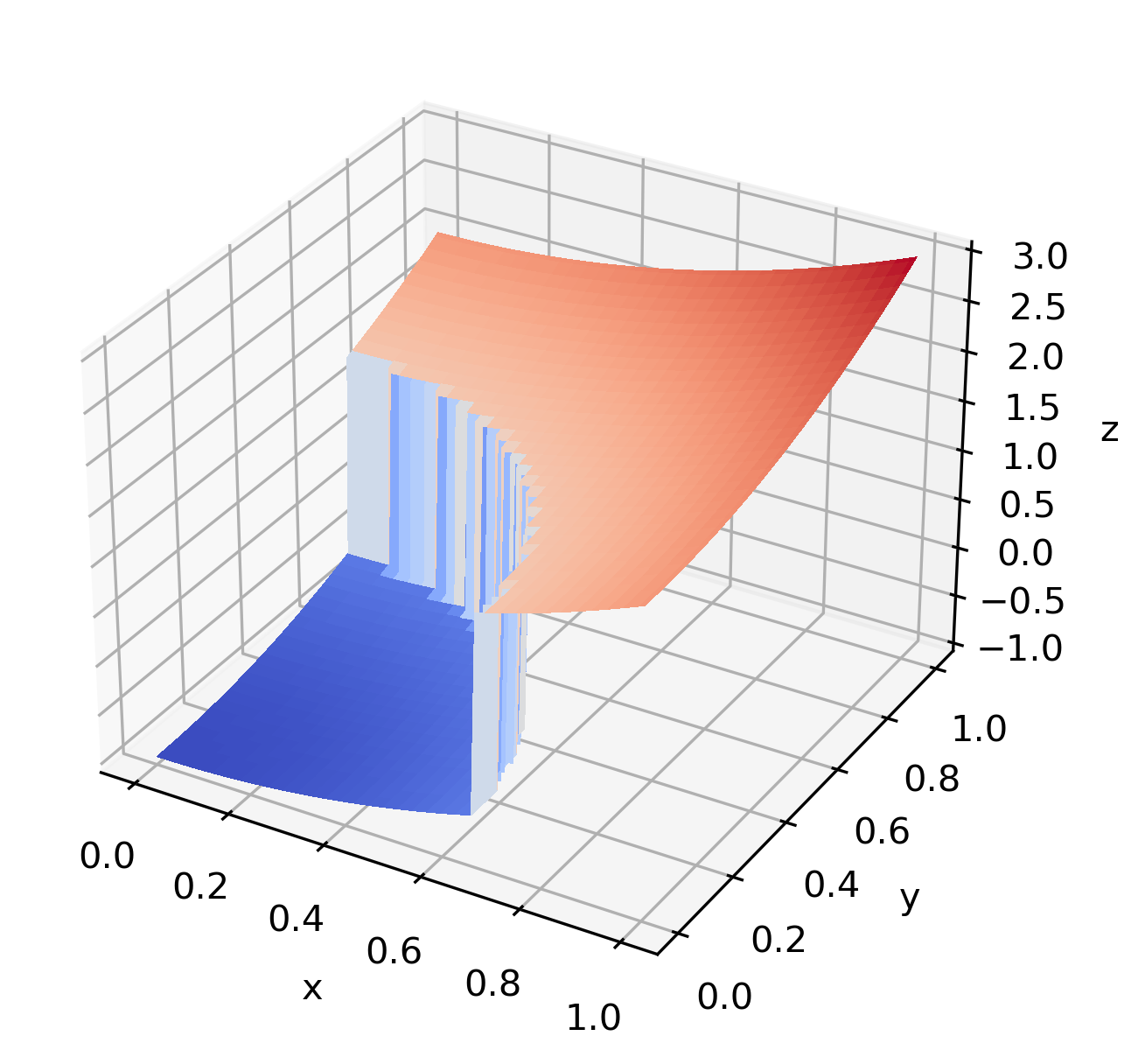}
\end{minipage}%
}%
\\
\subfigure[A 2--4000--1 ReLU NN function approximation\label{comparisonuhat_one}]{
\begin{minipage}[t]{0.4\linewidth}
\centering
\includegraphics[width=1.8in]{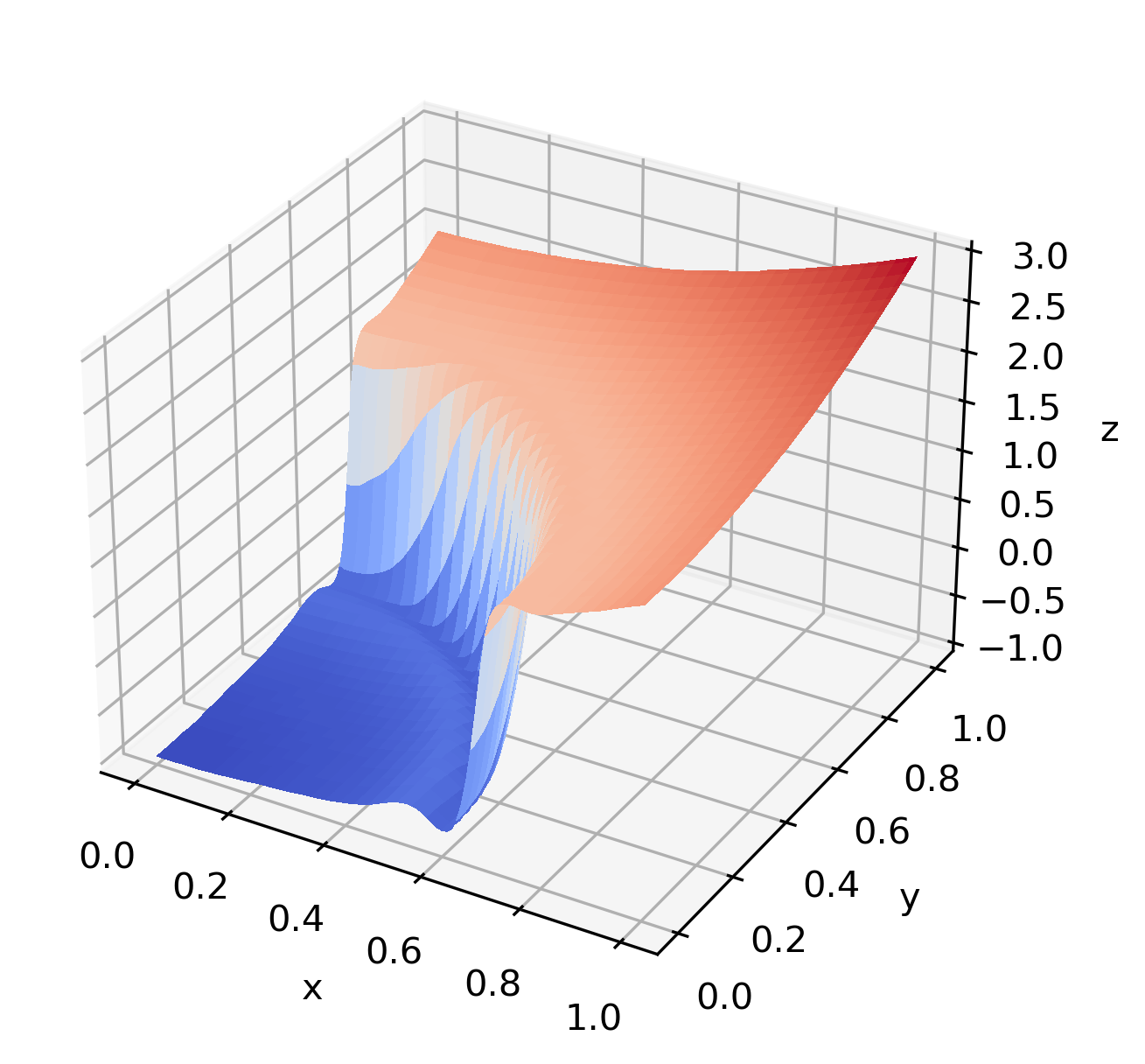}
\end{minipage}%
}%
\hspace{0.2in}
\subfigure[A 2--65--65--1 ReLU NN function approximation\label{comparisonuhat2}]{
\begin{minipage}[t]{0.4\linewidth}
\centering
\includegraphics[width=1.8in]{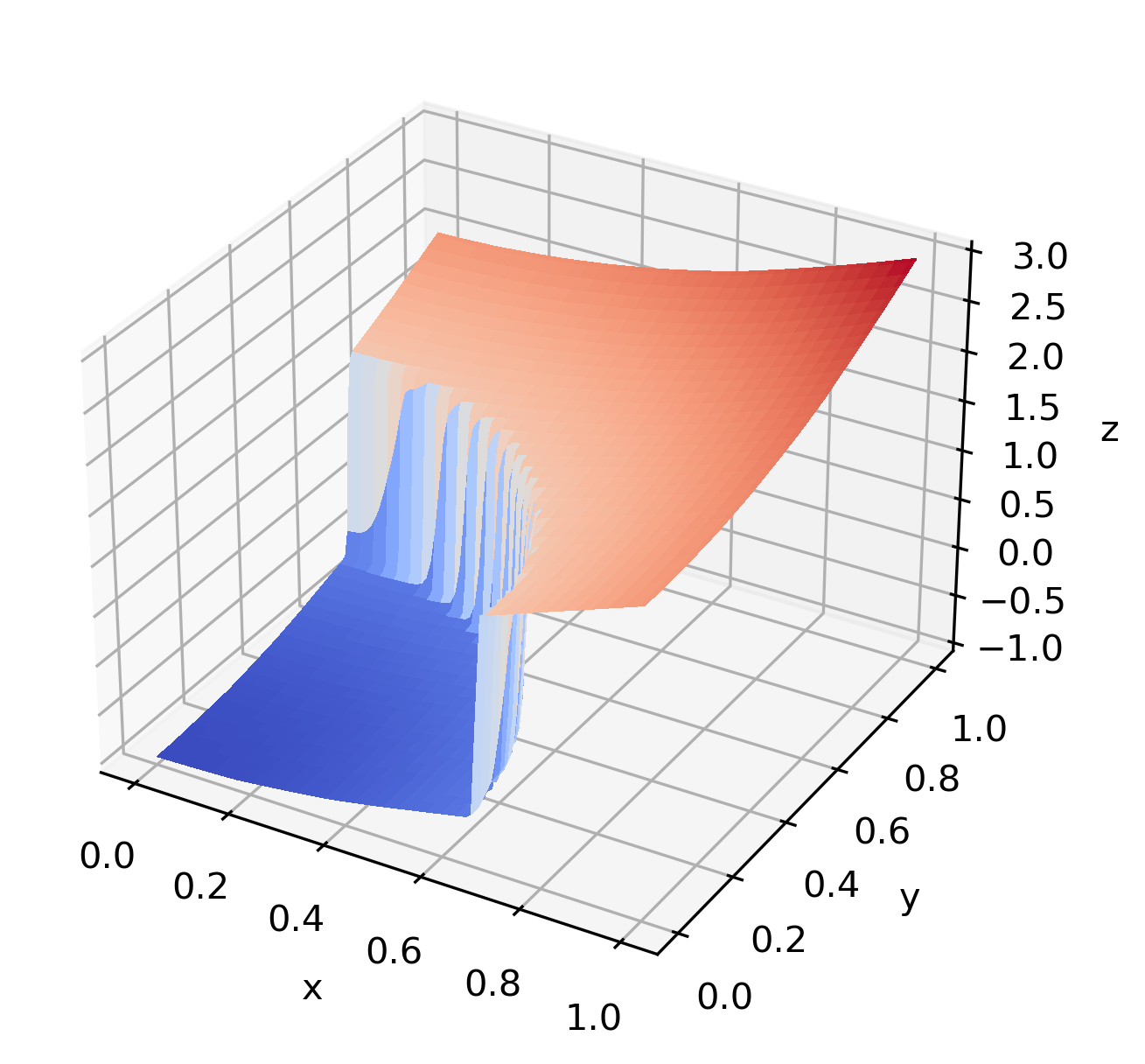}
\end{minipage}%
}%
\\
\subfigure[The trace of Figure \ref{comparisonuhat_one} on $y=x$\label{verticaluhat_one}]{
\begin{minipage}[t]{0.4\linewidth}
\centering
\includegraphics[width=1.8in]{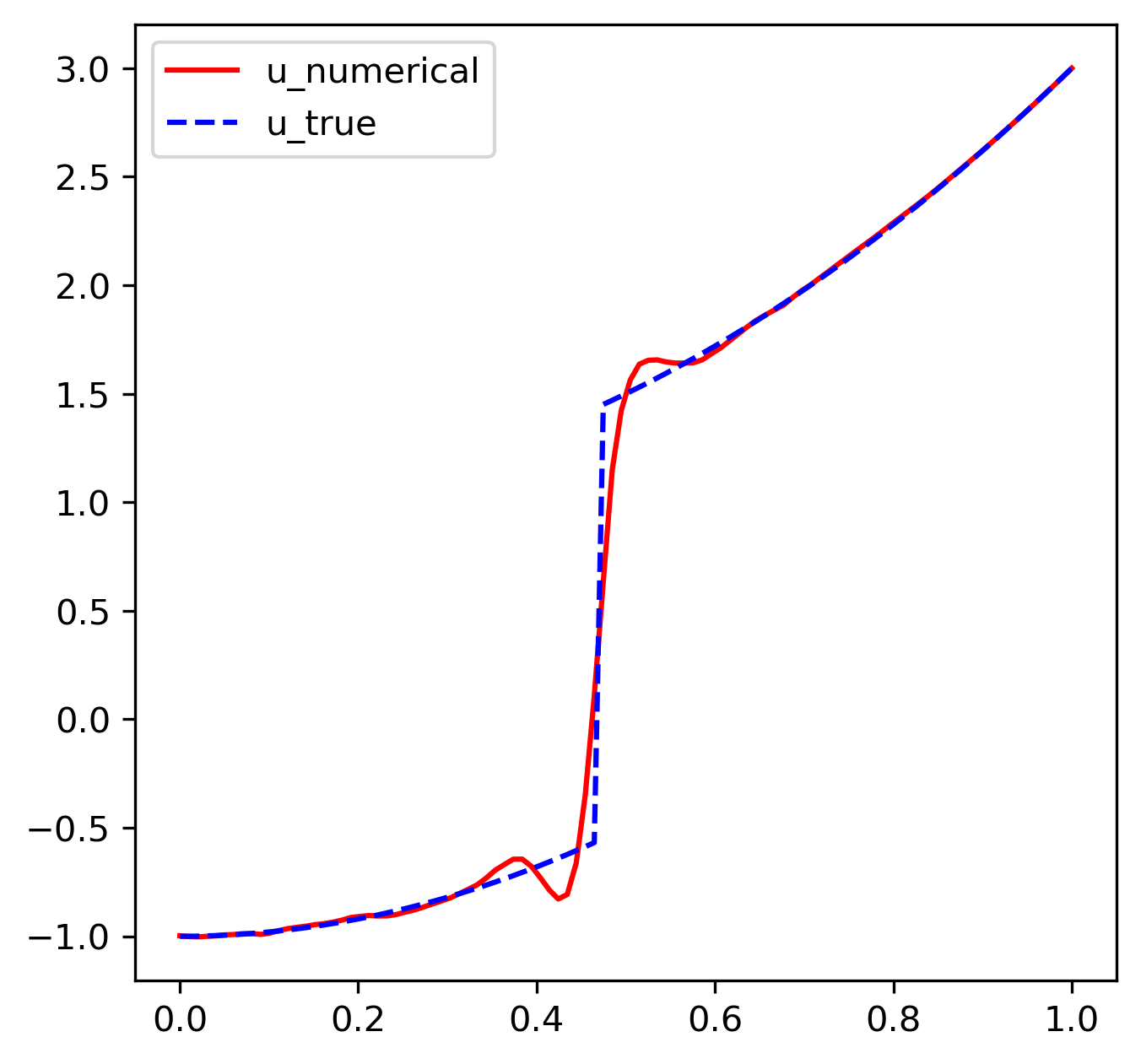}
\end{minipage}%
}%
\hspace{0.2in}
\subfigure[The trace of Figure \ref{comparisonuhat2} on $y=x$\label{verticaluhat}]{
\begin{minipage}[t]{0.4\linewidth}
\centering
\includegraphics[width=1.8in]{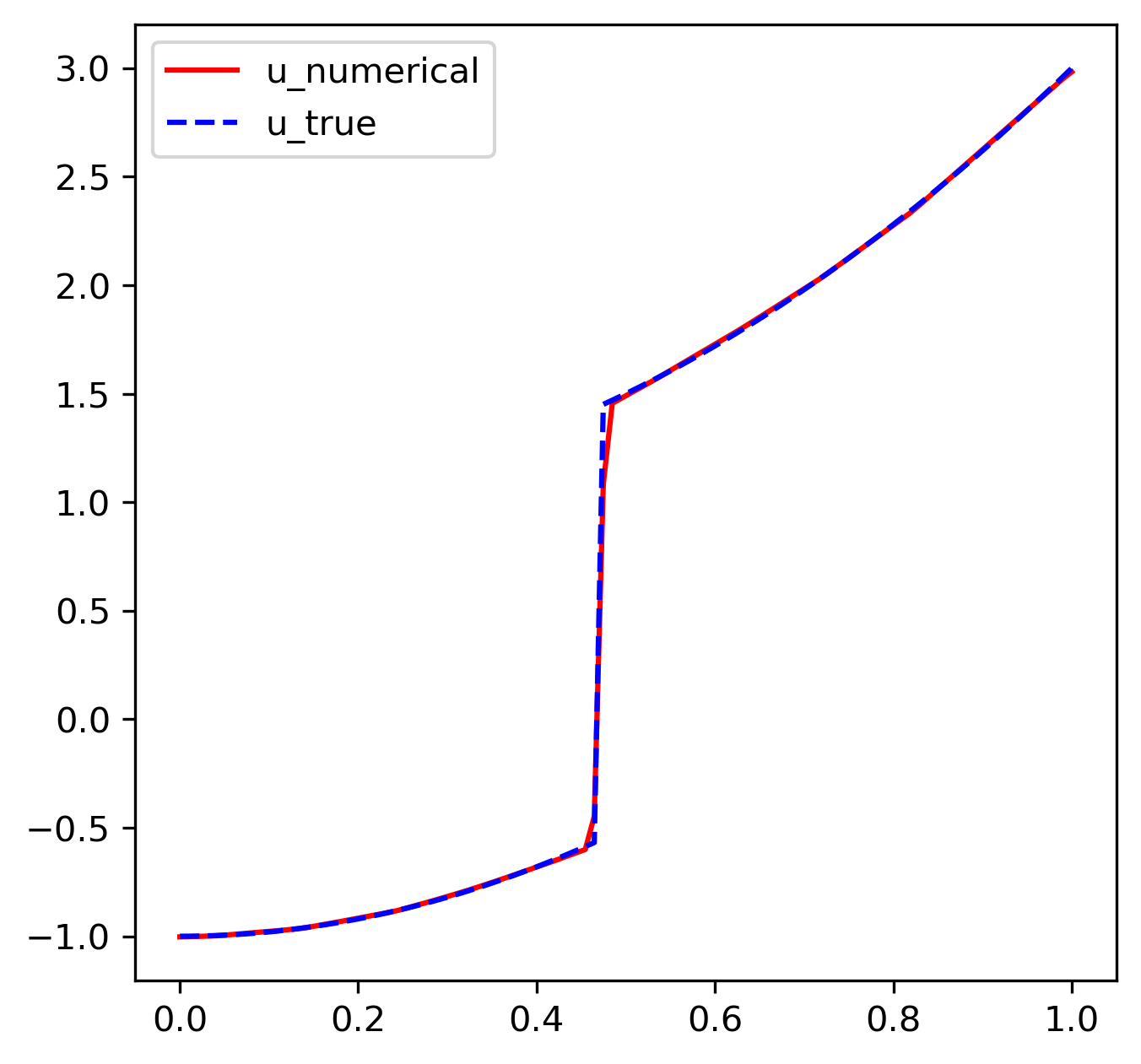}
\end{minipage}%
}%
\\
\subfigure[The breaking hyperplanes of the approximation in Figure \ref{comparisonuhat_one}\label{breakinguhat_one}]{
\begin{minipage}[t]{0.4\linewidth}
\centering
\includegraphics[width=1.8in]{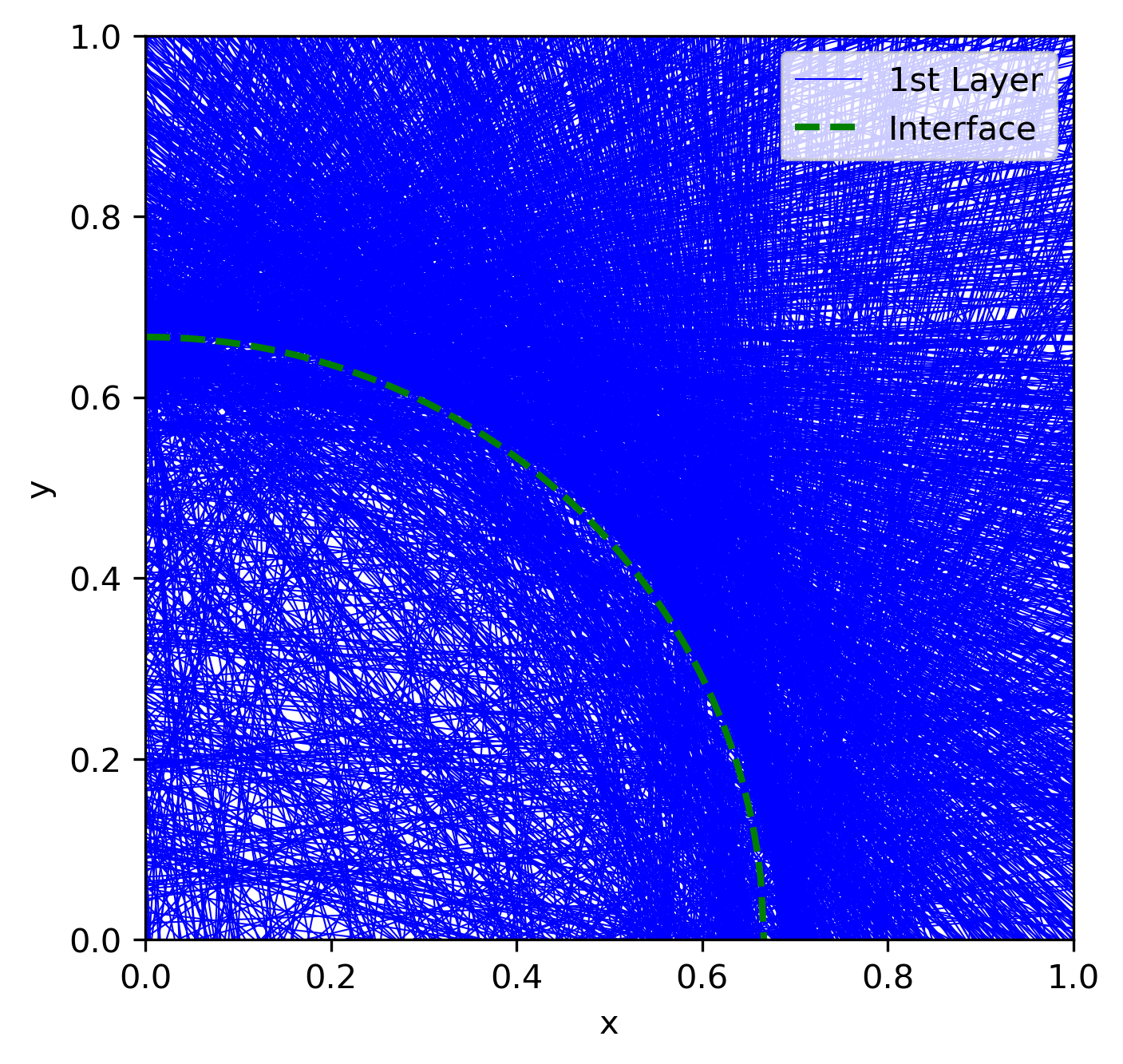}
\end{minipage}%
}%
\hspace{0.2in}
\subfigure[The breaking hyperplanes of the approximation in Figure \ref{comparisonuhat2}\label{breakinguhat}]{
\begin{minipage}[t]{0.4\linewidth}
\centering
\includegraphics[width=1.8in]{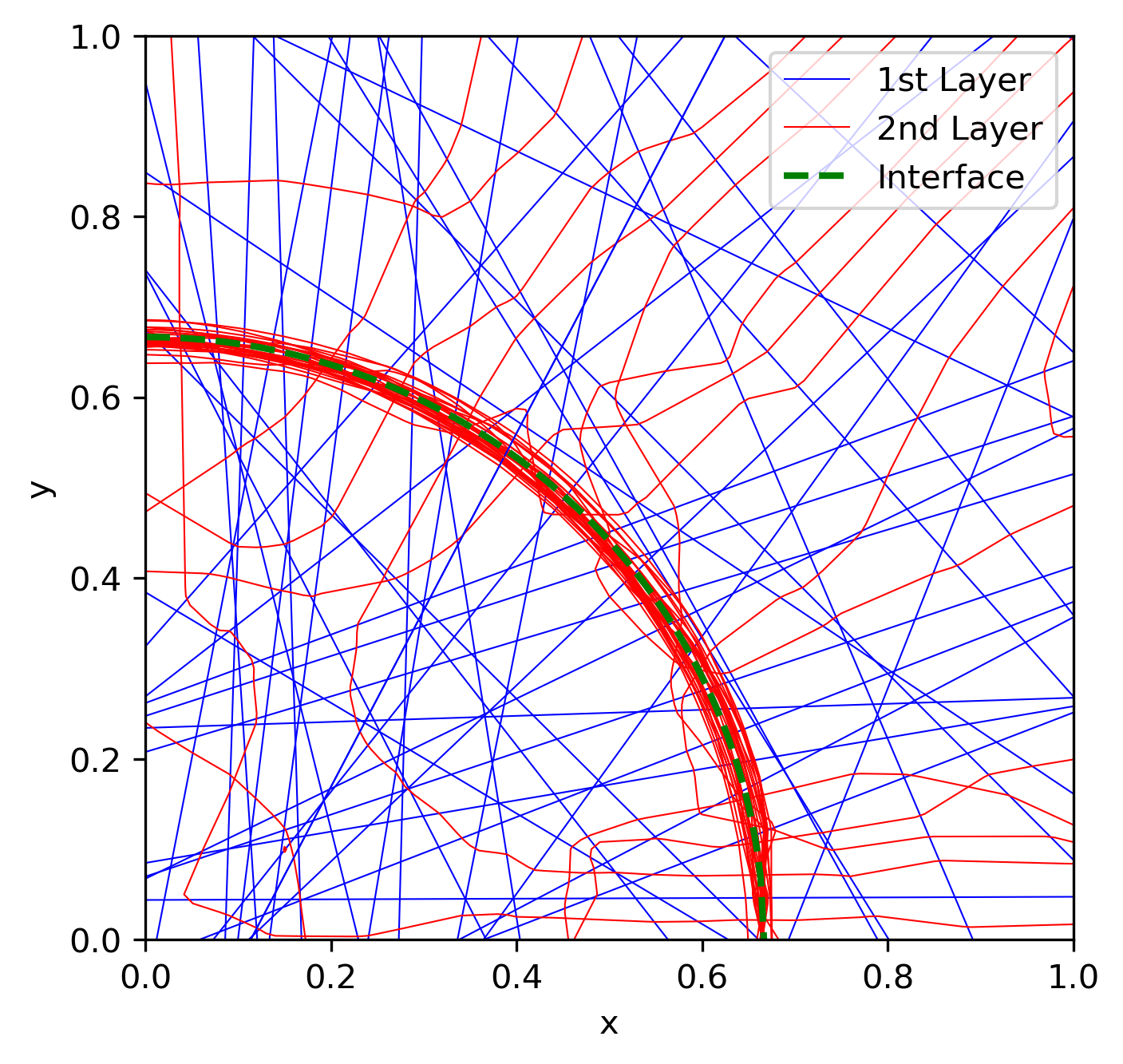}
\end{minipage}%
}%
\caption{Approximation results of the problem in \Cref{2d test4 section}}
\end{figure}

The advective velocity field is a variable field given by
\begin{equation}
\bm{\beta}(x,y) = (-y,x),\,\, (x,y)\in\Omega.
\end{equation}
The inflow boundary and the inflow boundary condition are given by
\begin{eqnarray*}
\Gamma_{-}&=&\{(1,y):y\in(0,1)\}\cup\{(x,0):x\in(0,1)\}\\[2mm]
\mbox{and }\,\, g(x,y)&=&\left\{ \begin{array}{rl}
 -1+x^2+y^2,& (x,y)\in \Gamma^1_-\equiv\{(x,0): x\in(0,\frac{2}{3})\}, \\[2mm]
 1+x^2+y^2, & (x,y)\in \Gamma^2_-=\Gamma_-\setminus \Gamma_-^1,
\end{array}\right.
\end{eqnarray*} 
respectively. The following right-hand side function is (see \cref{comparisonuhat1})
\begin{equation}
f(x,y)=\left\{ \begin{array}{rl}
 -1+x^2+y^2,& (x,y)\in \Omega_1\equiv\{(x,y)\in\Omega:y< \sqrt{\frac{4}{9}-x^2}\}, \\[2mm]
 1+x^2+y^2, & (x,y)\in \Omega_2=\Omega\setminus\Omega_1.
\end{array}\right.
\end{equation}
200000 iterations were implemented with 2--4000--1 and 2--65--65--1 ReLU NN functions. The numerical results are presented in \cref{2d test4,2d test4 table}. \cref{2d test4 table} indicates that the 2-layer network structure is capable of approximating the solution (\cref{comparisonuhat1}) on average, but \cref{comparisonuhat_one,verticaluhat_one,breakinguhat_one} show difficulty around the discontinuity interface (\cref{interfaceuhat}). Again, the 3-layer network structure with 28\% of the number of parameters of the 2-layer network structure presented better error results (\cref{2d test4 table}) and approximated the solution accurately pointwise (\cref{comparisonuhat2,verticaluhat}). Unlike other examples, some of the second-layer breaking hyperplanes (\cref{breakinguhat}) of the approximation spread out on the whole domain in addition to those around the interface, which implies that they are necessary for approximating the solution with $\hat{u}\neq0$ ($\hat{u}=x^2+y^2$ in this example).

\begin{table}[htbp]\label{2d test4 table}
\caption{Relative errors of the problem in \Cref{2d test4 section}}
\centering
\begin{tabular}{|l|l|l|l|l|}
\hline
Network structure  &$\frac{\|u-{u}^{_N}_{_\cT}\|_0}{\|u\|_0}$ &$\frac{\vertiii{u-{u}^{_N}_{_\cT}}_{\bm\beta}}{\vertiii{u}_{\bm\beta}}$ & $\frac{\mathcal{L}^{1/2}({u}^{_N}_{_\cT},\bf f)}{\mathcal{L}^{1/2}({u}^{_N}_{_\cT},\bf 0)}$ & Parameters \\ \hline
2--4000--1  & 0.088349 & 0.108430 & 0.058213   & 16001\\ \hline
2--65--65--1  & 0.048278 & 0.073095 & 0.015012   & 4551\\ \hline
\end{tabular}
\end{table}

\subsubsection{A problem with a sharp transition layer\label{2d test5 section}}
\begin{figure}[htbp]\label{2d test5}
\centering
\subfigure[A 2--70--70--1 ReLU NN function approximation\label{comparison2dc22}]{
\begin{minipage}[t]{0.4\linewidth}
\centering
\includegraphics[width=1.8in]{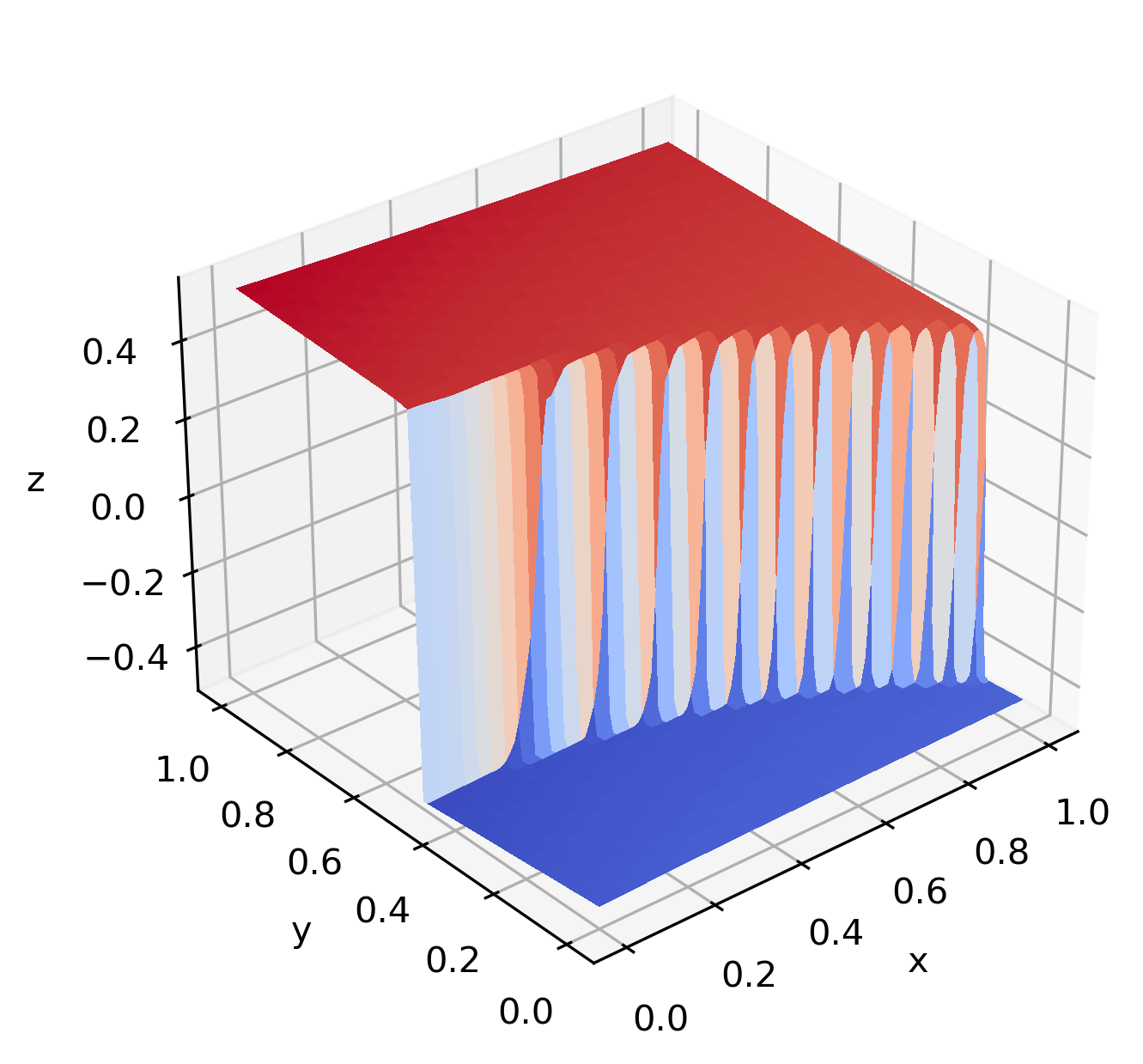}
\end{minipage}%
}%
\hspace{0.2in}
\subfigure[The exact solution with $\varepsilon=10^{-10}$ in \cref{2d test5 exact}\label{comparison2dc21}]{
\begin{minipage}[t]{0.4\linewidth}
\centering
\includegraphics[width=1.8in]{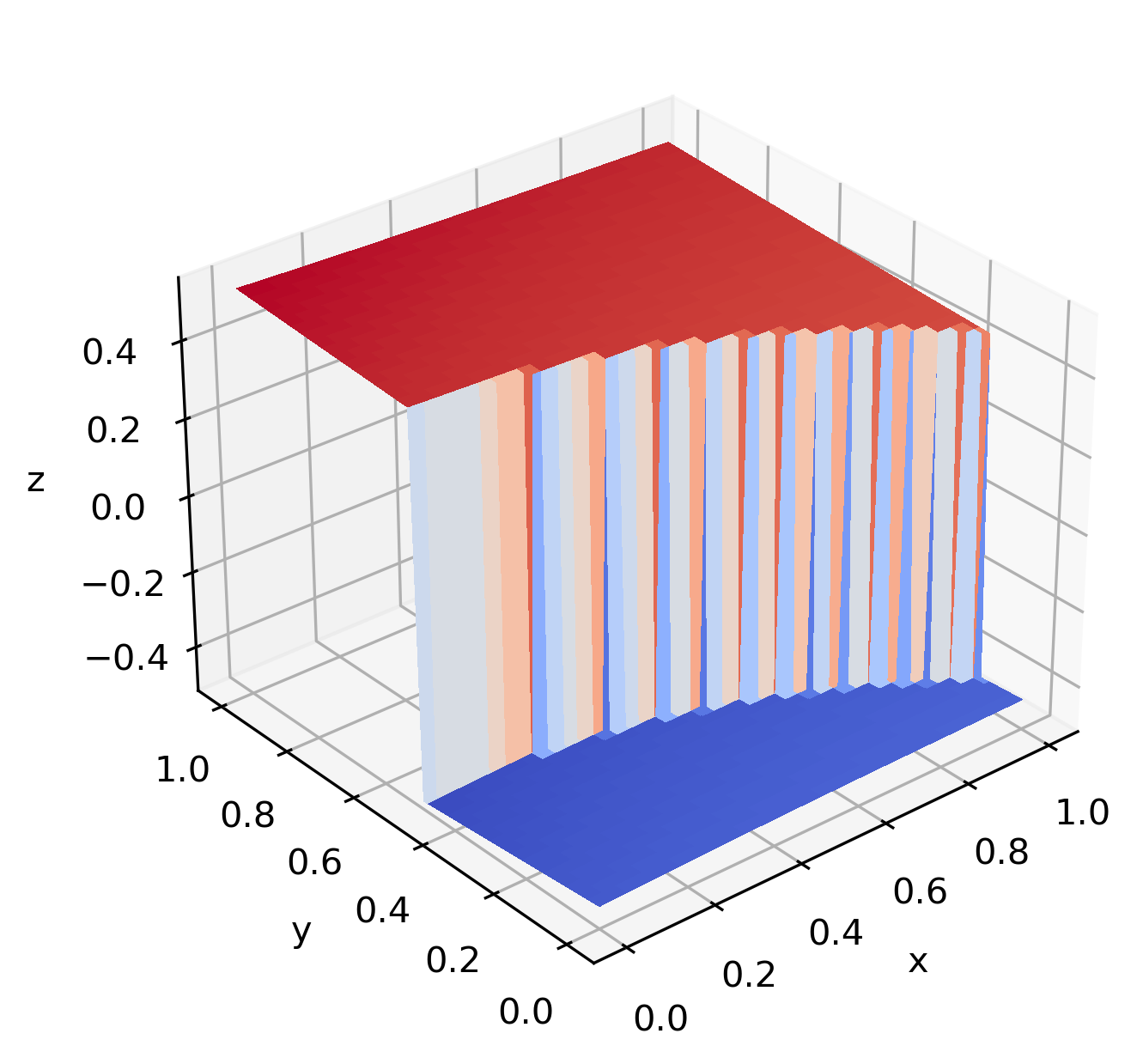}
\end{minipage}%
}%
\\
\subfigure[The trace of Figure \ref{comparison2dc22} on $y=x$\label{vertical2dc2}]{
\begin{minipage}[t]{0.4\linewidth}
\centering
\includegraphics[width=1.8in]{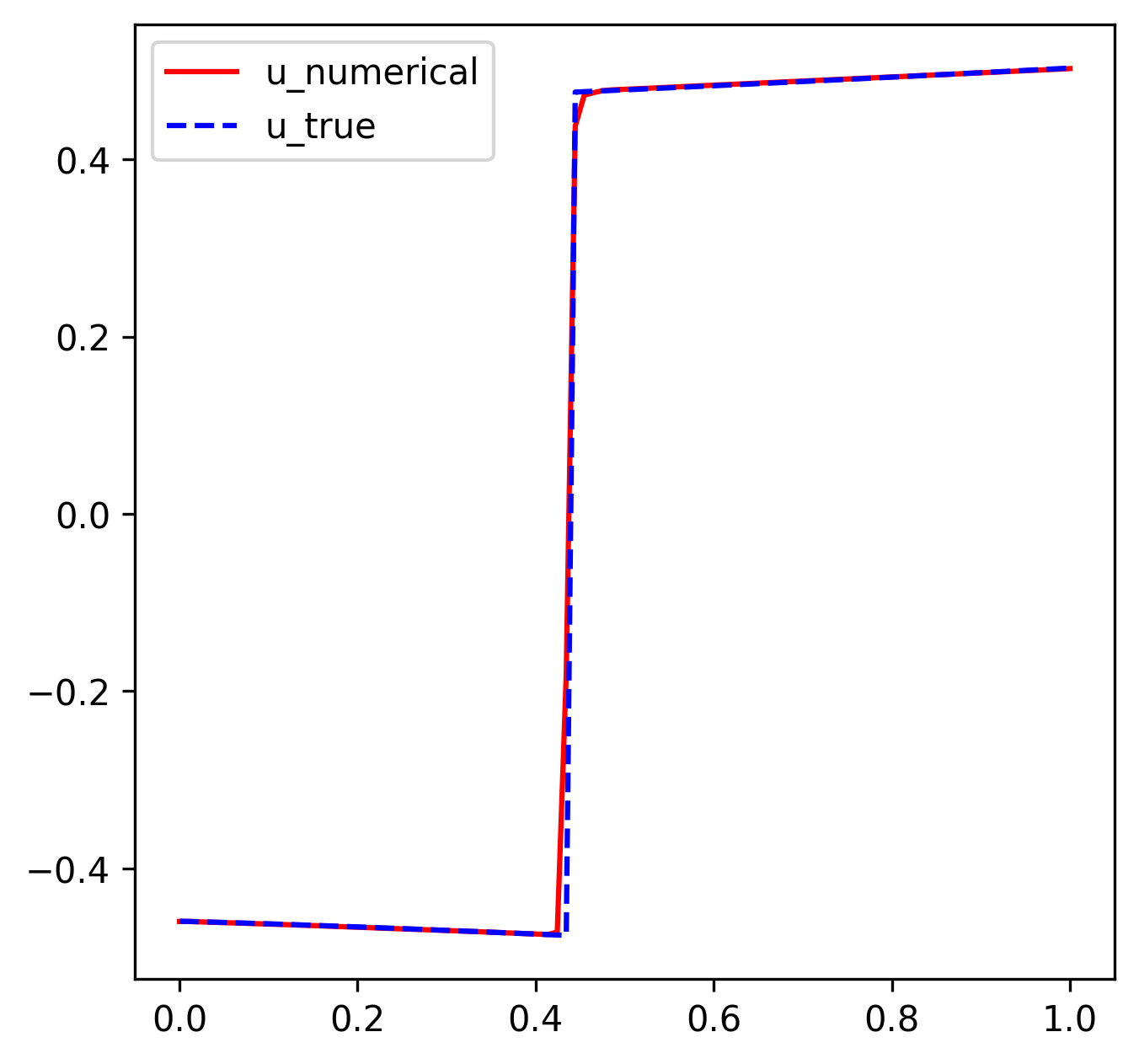}
\end{minipage}%
}%
\hspace{0.2in}
\subfigure[The breaking hyperplanes of the approximation in Figure \ref{comparison2dc22}\label{breaking2dc2}]{
\begin{minipage}[t]{0.4\linewidth}
\centering
\includegraphics[width=1.8in]{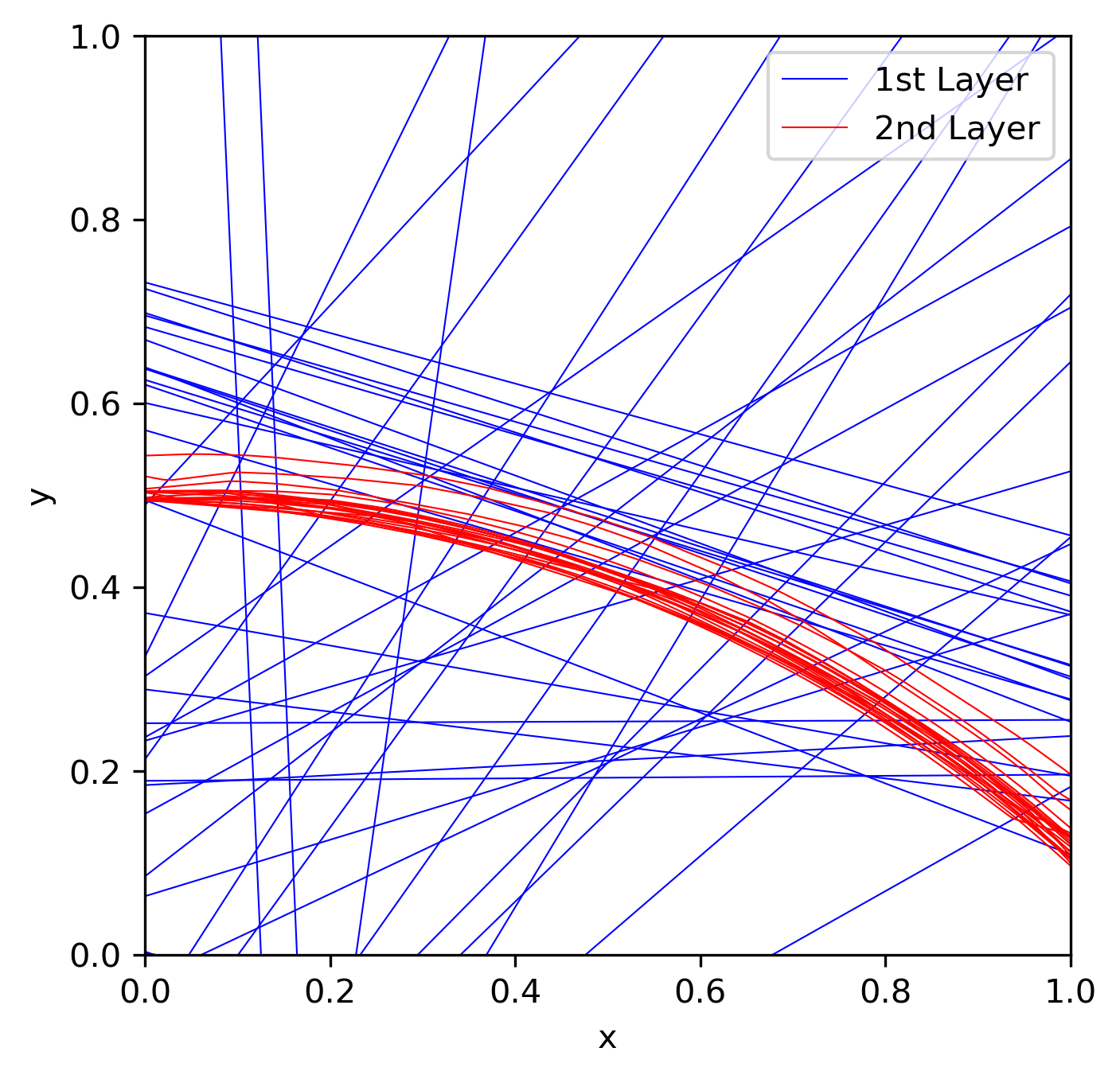}
\end{minipage}%
}%
\caption{Approximation results of the problem in \Cref{2d test5 section}}
\end{figure}

The advective velocity field is a variable field given by
\begin{equation}
\bm{\beta}(x,y) = (y+1,-x)/\sqrt{x^2+(y+1)^2},\,\, (x,y)\in\Omega,
\end{equation}
the inflow boundary is given by
\begin{equation*}
\Gamma_{-}=\{(0,y):y\in(0,1)\}\cup\{(x,1):x\in(0,1)\},
\end{equation*} 
and $f=0$. We choose an inflow boundary condition $g$ such that the exact solution is (see \cref{comparison2dc21})
\begin{equation}\label{2d test5 exact}
u(x,y)=\frac{1}{4}\exp\left(\gamma r\arcsin\left(\frac{y+1}{r}\right)\right)\arctan\left(\frac{r-1.5}{\varepsilon}\right) \text{ with }r=\sqrt{x^2+(y+1)^2}
\end{equation}

300000 iterations were implemented with 2--70--70--1 ReLU NN functions to approximate $u$ with $\varepsilon=10^{-10}$ in \cref{2d test5 exact}. The numerical results are presented in \cref{2d test5,2d test5 table}. The approximation results are similar to those of the previous examples. The same PDE was solved in \cite{liu2020adaptive,liu2020adaptive2,Zhang2022}, and in the experiments, the layer cannot be fully resolved and should be viewed as discontinuous (see \cref{comparison2dc21}). The $L^2$ errors of the least-squares finite element methods in \cite{liu2020adaptive,liu2020adaptive2} are approximately between $4\times10^{-2}$ and $6\times10^{-2}$ with $10^4$ to $10^6$ degrees of freedom, whereas the error by the LSNN method is approximately $3\times10^{-2}$ with 5251 parameters (\cref{2d test5 table}). The discontinuous Galerkin finite element methods (DGFEMs) in \cite{Zhang2022} give similar results with a P0-DGFEM solution having no overshoot and a P1-DGFEM solution having a non-trivial overshoot. There is no overshoot from the LSNN method (\cref{comparison2dc22,vertical2dc2}).

\begin{table}[htbp]\label{2d test5 table}
\caption{Relative errors of the problem in \Cref{2d test5 section}}
\centering
\begin{tabular}{|l|l|l|l|l|l|}
\hline
Network structure & $\|u-{u}^{_N}_{_\cT}\|_0$ &$\frac{\|u-{u}^{_N}_{_\cT}\|_0}{\|u\|_0}$ &$\frac{\vertiii{u-{u}^{_N}_{_\cT}}_{\bm\beta}}{\vertiii{u}_{\bm\beta}}$ & $\frac{\mathcal{L}^{1/2}({u}^{_N}_{_\cT},\bf f)}{\mathcal{L}^{1/2}({u}^{_N}_{_\cT},\bf 0)}$ & Parameters \\ \hline
2--70--70--1  & 0.032229 & 0.066337 & 0.069926 &  0.023664 & 5251\\ \hline
\end{tabular}
\end{table}

\subsection{Three-dimensional problems}
We present numerical results for three three-dimensional test problems with piecewise constant or variable advection velocity fields whose solutions are piecewise constant along a connected series of planes or a surface. All three test problems are defined on the domain $\Omega=(0,1)^3$, approximation results are depicted on $z=0.505$ (except for the last test problem on $z=0.205$), and again, as in the experiments for $d=2$, we have
\[
 \vertiii{u-u_{_{N}}}_{\bm\beta}
 \leq C\,\big|\alpha_1-\alpha_2\big|\, \sqrt{\varepsilon}.
\]

\subsubsection{A problem with a 2-plane segment interface}\label{3d test1 section}
\begin{figure}[htbp]\label{3d test1}
\centering
\subfigure[The interface\label{interface3d2l}]{
\begin{minipage}[t]{0.4\linewidth}
\centering
\includegraphics[width=1.8in]{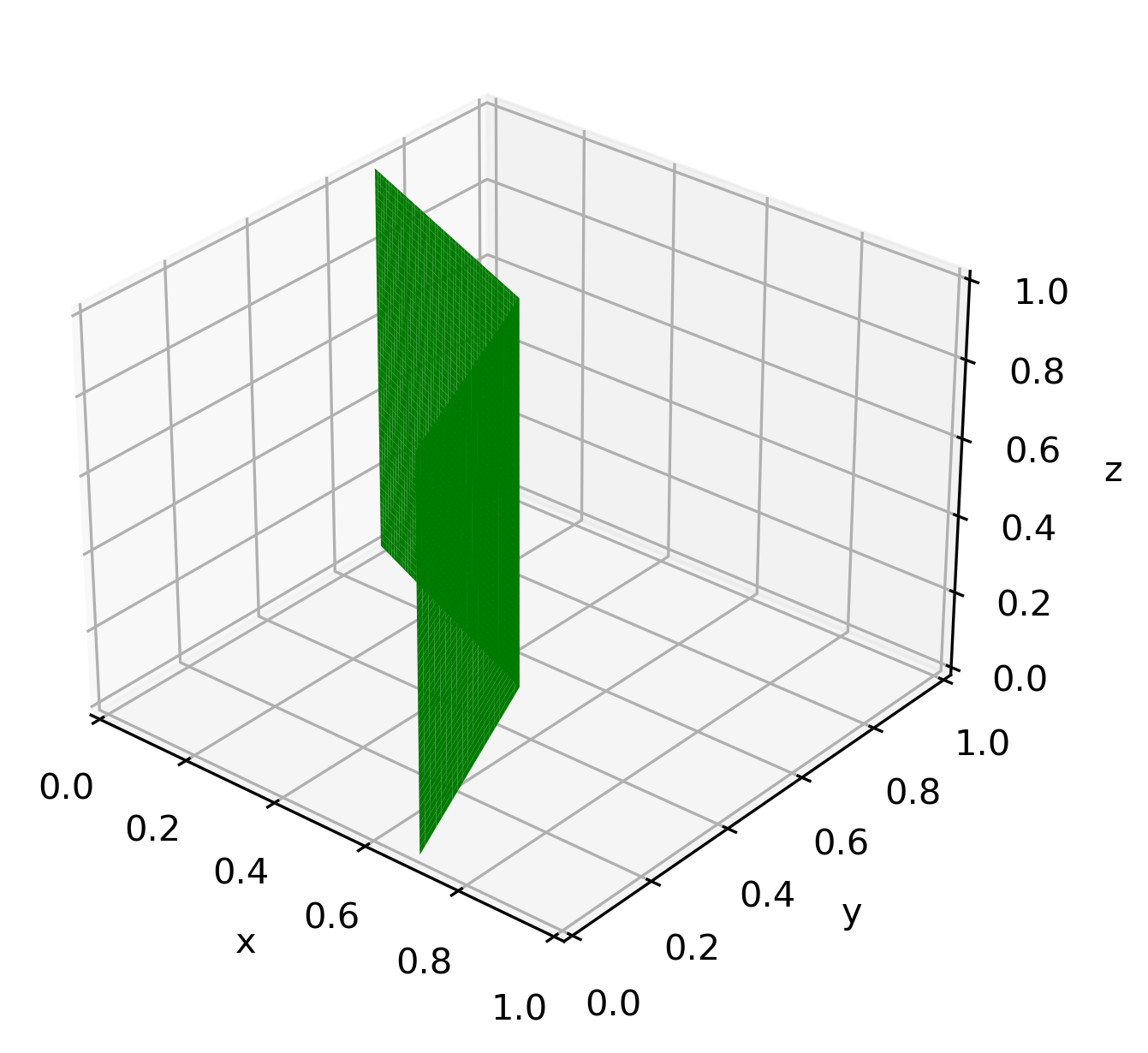}
\end{minipage}%
}%
\hspace{0.2in}
\subfigure[The exact solution\label{comparison3d2l1}]{
\begin{minipage}[t]{0.4\linewidth}
\centering
\includegraphics[width=1.8in]{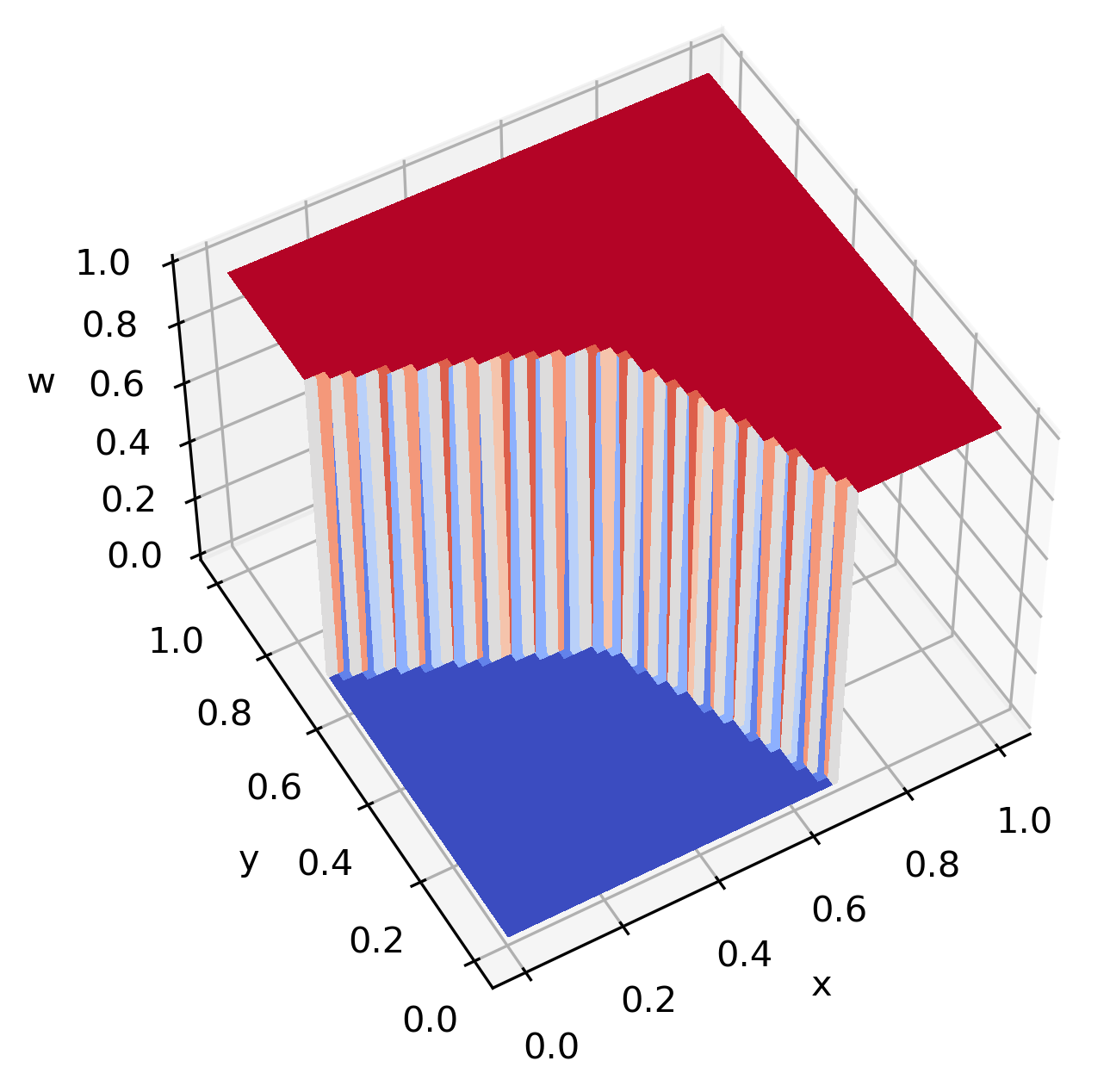}
\end{minipage}%
}%
\\
\subfigure[A 3--300--1 ReLU NN function approximation\label{comparison3d2l_one}]{
\begin{minipage}[t]{0.4\linewidth}
\centering
\includegraphics[width=1.8in]{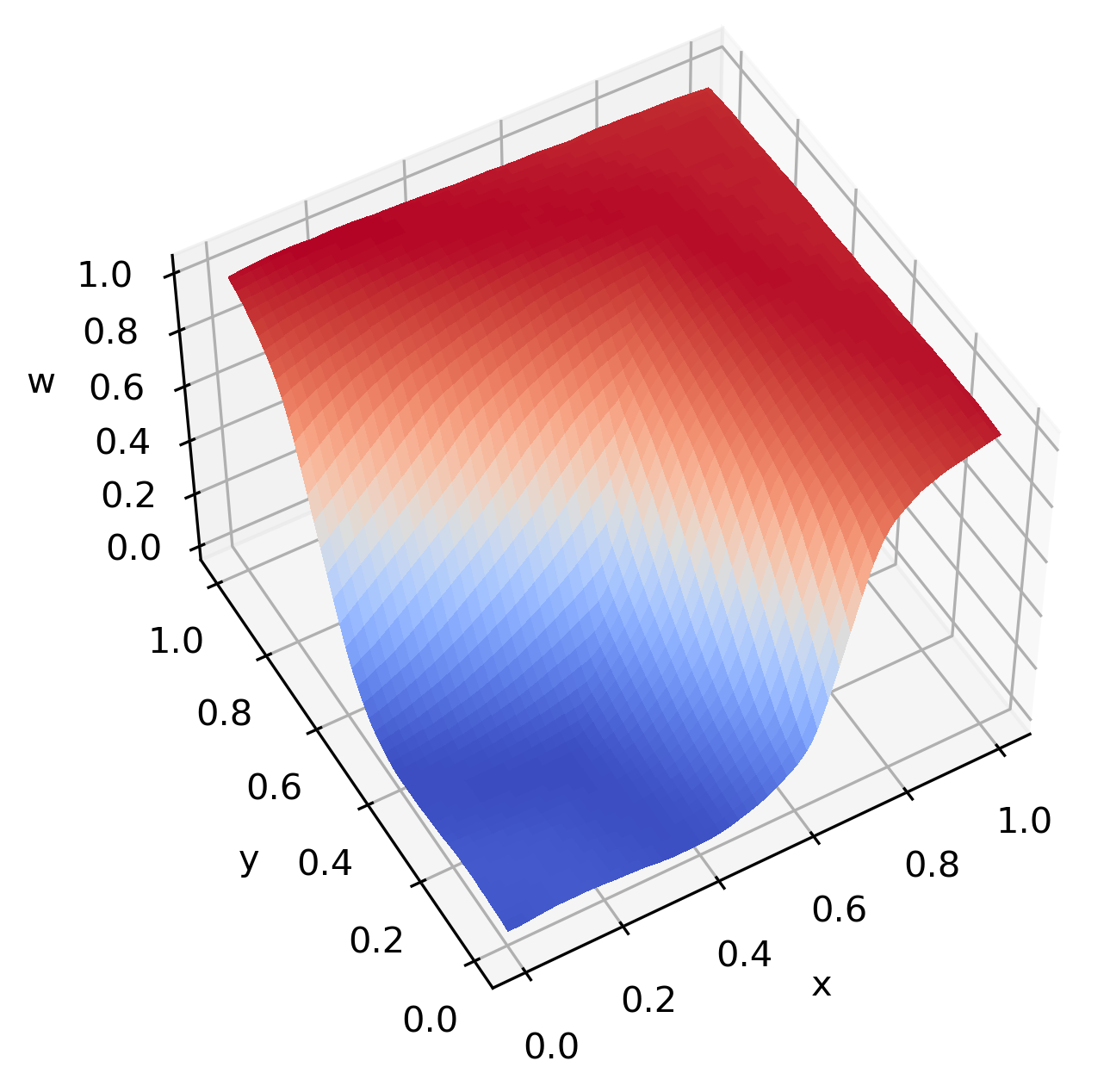}
\end{minipage}%
}%
\hspace{0.2in}
\subfigure[A 3--5--5--1 ReLU NN function approximation\label{comparison3d2l2}]{
\begin{minipage}[t]{0.4\linewidth}
\centering
\includegraphics[width=1.8in]{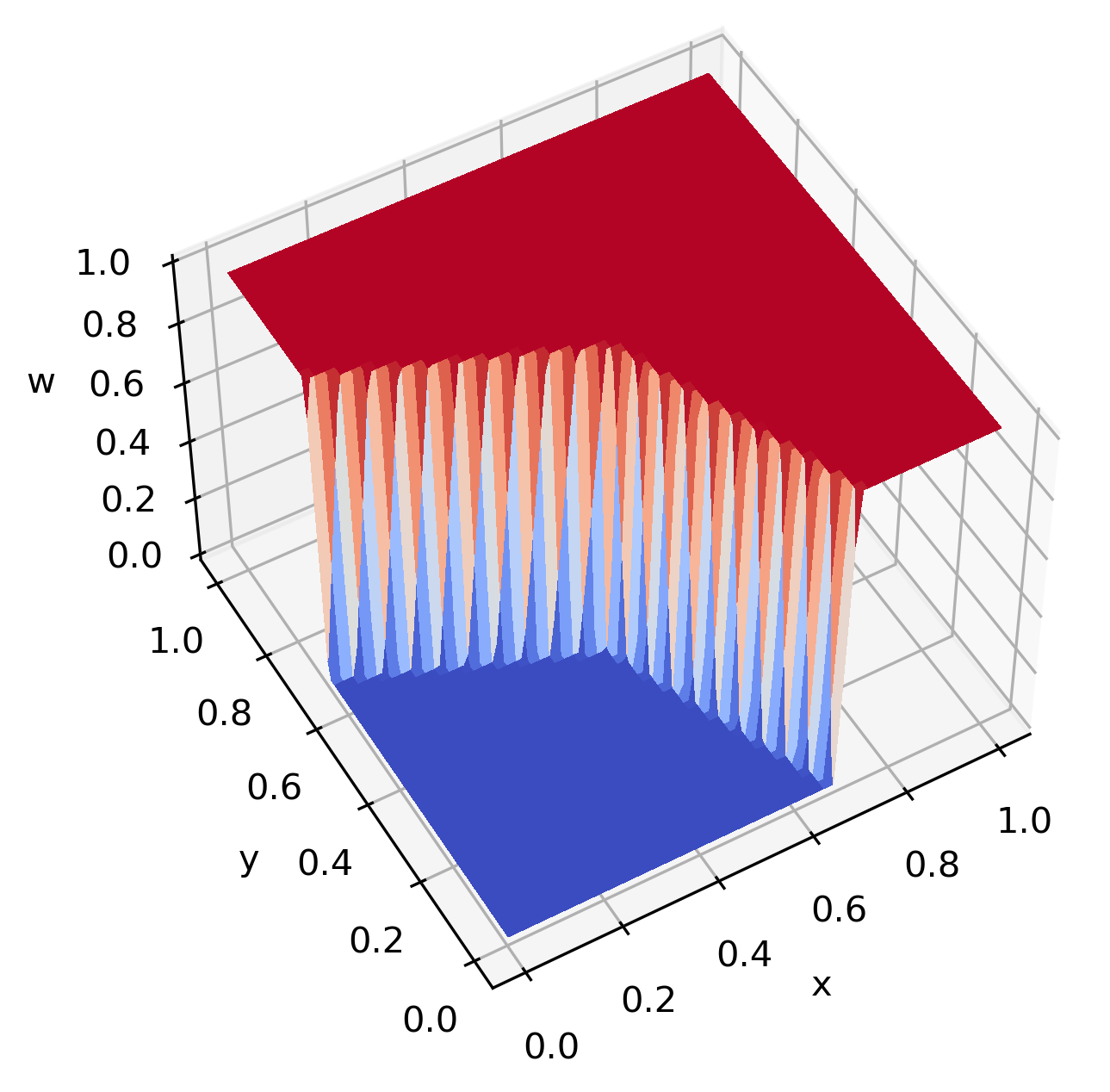}
\end{minipage}%
}%
\\
\subfigure[The trace of Figure \ref{comparison3d2l_one} on $y=x$\label{vertical3d2l_one}]{
\begin{minipage}[t]{0.4\linewidth}
\centering
\includegraphics[width=1.8in]{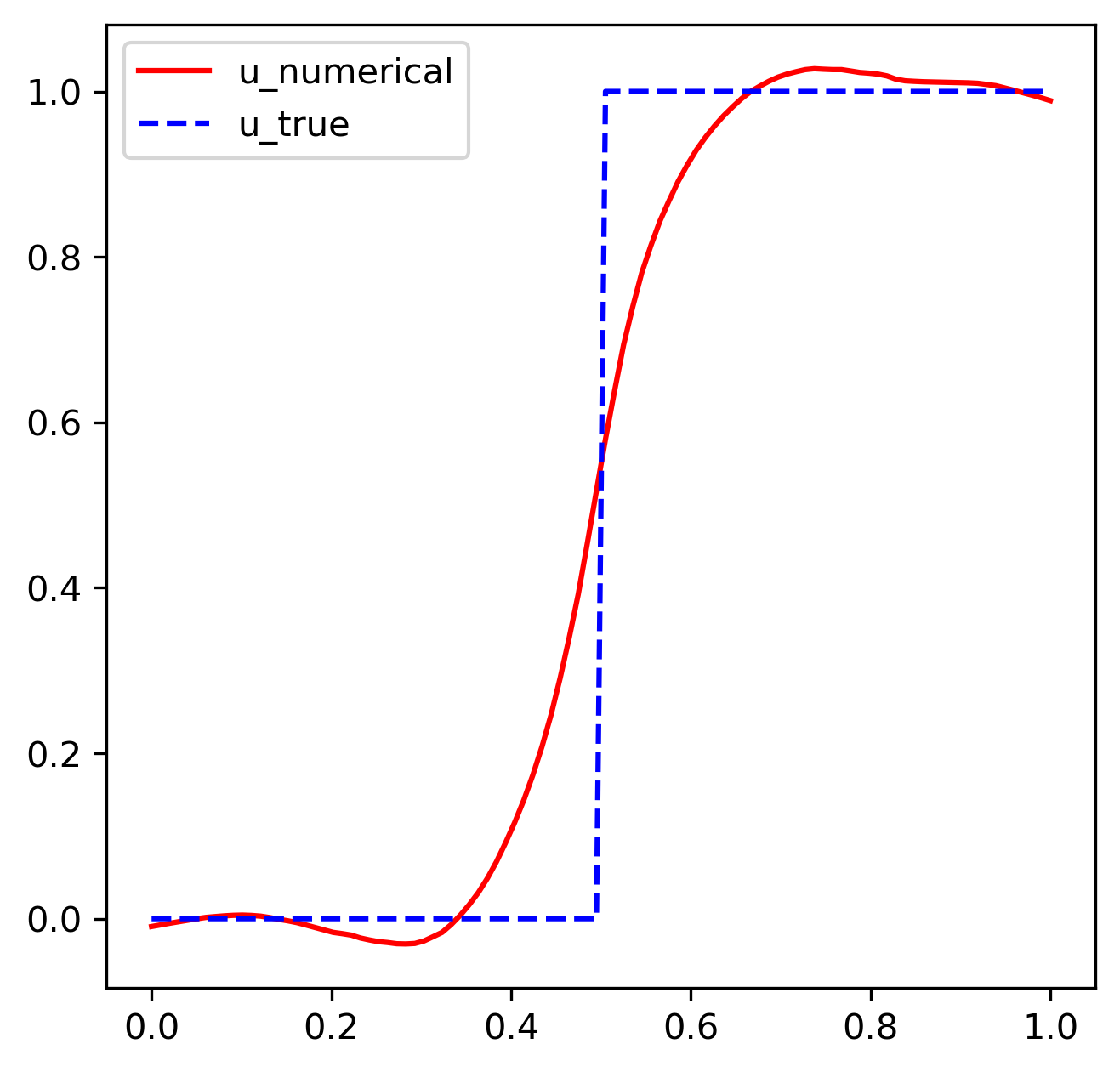}
\end{minipage}%
}%
\hspace{0.2in}
\subfigure[The trace of Figure \ref{comparison3d2l2} on $y=x$\label{vertical3d2l}]{
\begin{minipage}[t]{0.4\linewidth}
\centering
\includegraphics[width=1.8in]{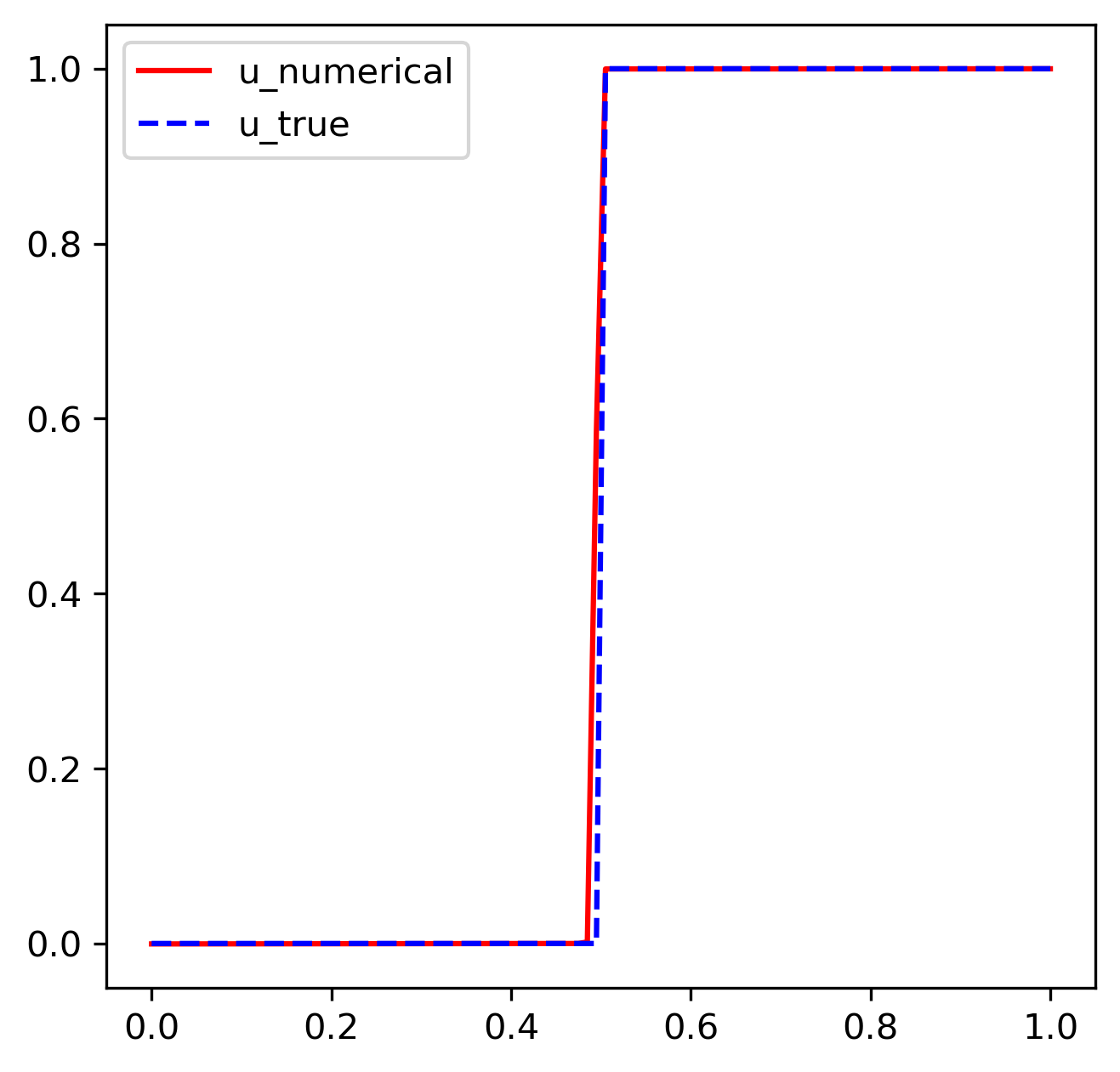}
\end{minipage}%
}%
\\
\subfigure[The breaking hyperplanes of the approximation in Figure \ref{comparison3d2l_one}\label{breaking3d2l_one}]{
\begin{minipage}[t]{0.4\linewidth}
\centering
\includegraphics[width=1.8in]{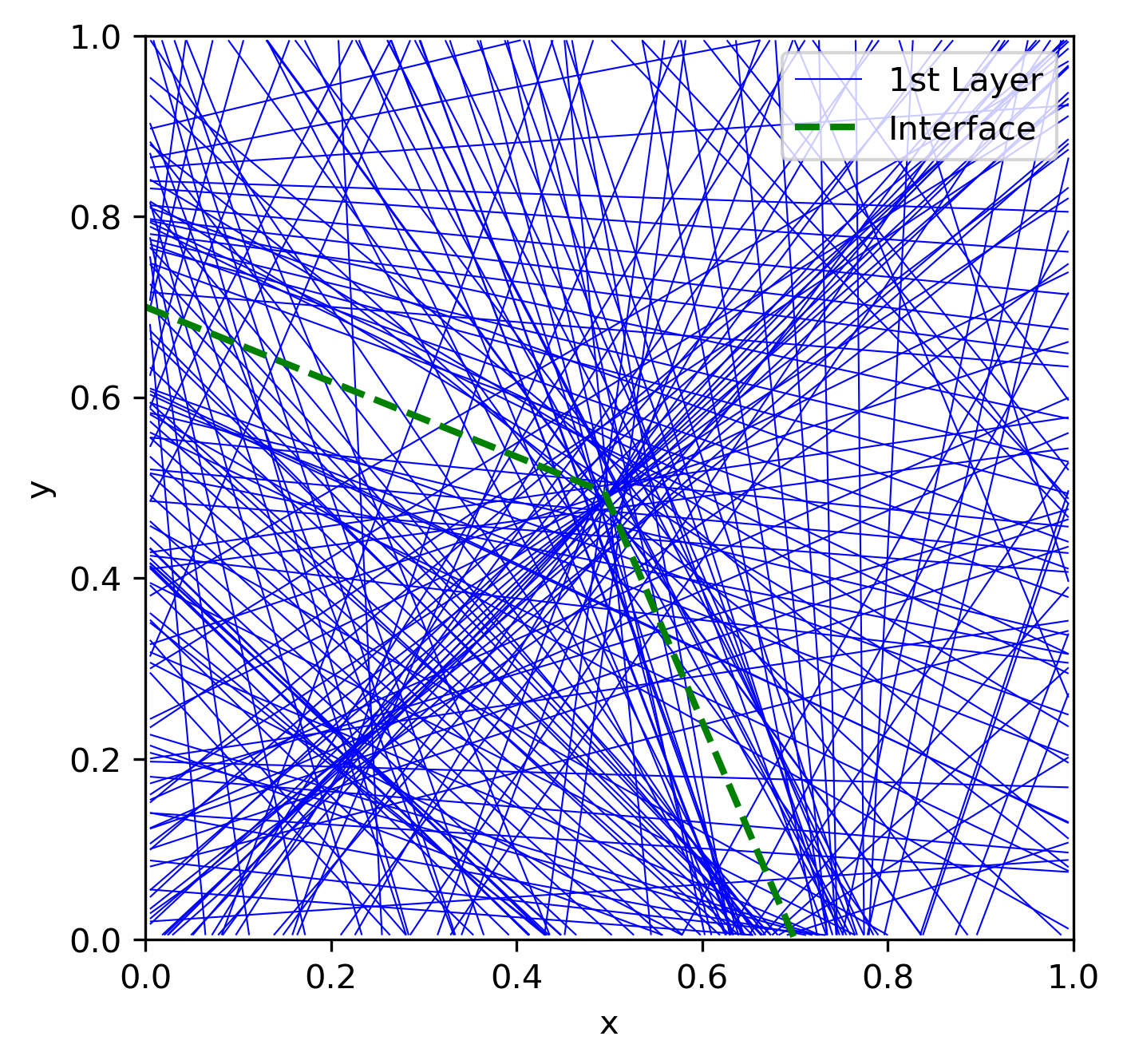}
\end{minipage}%
}%
\hspace{0.2in}
\subfigure[The breaking hyperplanes of the approximation in Figure \ref{comparison3d2l2}\label{breaking3d2l}]{
\begin{minipage}[t]{0.4\linewidth}
\centering
\includegraphics[width=1.8in]{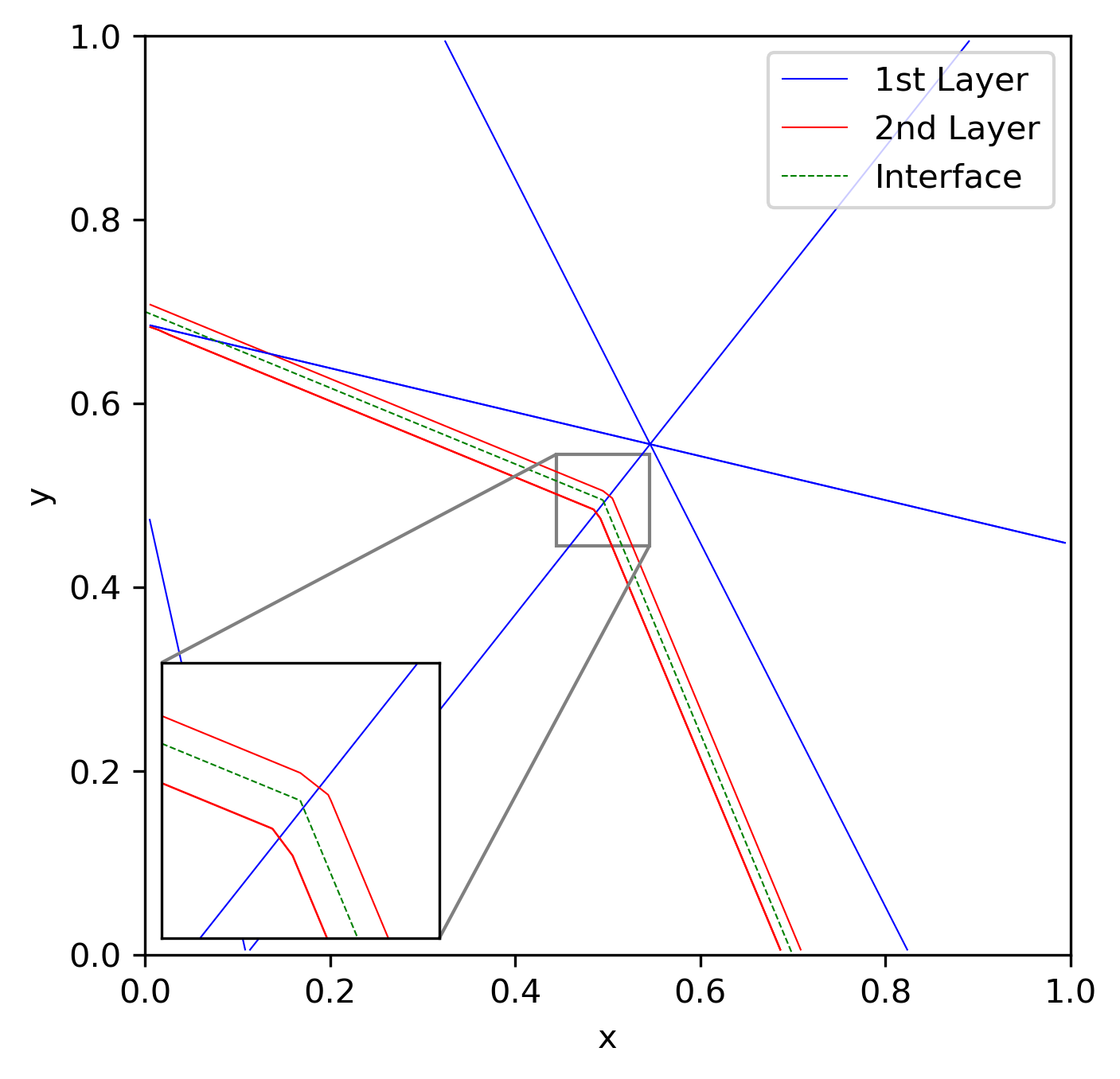}
\end{minipage}%
}%
\caption{Approximation results of the problem in \Cref{3d test1 section}}
\end{figure}

Let $\gamma=f=0$, $\bar{\Omega}=\bar{\Upsilon}_1\cup \bar{\Upsilon}_2$, and 
\[
    \Upsilon_1=\{(x,y,z)\in\Omega:y<x\}\text{ and }
    \Upsilon_2=\{(x,y,z)\in\Omega:y\ge x\}.
\]
The advective velocity field is a piecewise constant field given by \begin{equation}
\bm{\beta}(x,y,z) =\left\{ \begin{array}{rl}
(1-\sqrt{2},1,0)^T,&(x,y,z)\in\Upsilon_1\\[2mm]
(-1,\sqrt{2}-1,0)^T,&(x,y,z)\in\Upsilon_2.
\end{array}\right.
\end{equation}
The inflow boundary and the inflow boundary condition are given by
\begin{eqnarray*}
\Gamma_{-}&=&\{(x,0,z):x,z\in(0,1)\}\cup\{(1,y,z):y,z\in(0,1)\}\\[2mm]
\mbox{and }\,\, g(x,y,z)&=&\left\{ \begin{array}{rl}
 0,& (x,y,z)\in \Gamma^1_-\equiv \{(x,0,z): x\in(0,0.7),\ z\in(0,1)\}, \\[2mm]
 1, &(x,y,z)\in \Gamma^2_-=\Gamma_-\setminus \Gamma_-^1,
\end{array}\right.
\end{eqnarray*} 
respectively. Let
\[
\Omega_1=\{(x,y,z)\in\Omega:y<(1-\sqrt{2})x+0.7,\ y<\tfrac{1}{1-\sqrt{2}}(x-0.7)\}.
\]
The exact solution is a unit step function in three dimensions (see \cref{comparison3d2l1})
\begin{equation}
u(x,y,z)=\left\{ \begin{array}{rl}
 0,& (x,y,z)\in \Omega_1, \\[2mm]
 1, & (x,y,z)\in \Omega_2=\Omega\setminus\Omega_1.
\end{array}\right.
\end{equation}
100000 iterations were implemented with 3--300--1 and 3--5--5--1 ReLU NN functions (depth $\lceil \log_2(d+1)\rceil+1=3$ for $d=3$). The numerical results are presented in \cref{3d test1,3d test1 table}. For three dimensions again, the 2-layer network structure with a large number of parameters generated a transition layer along the discontinuity interface (\cref{interface3d2l,breaking3d2l_one}) but had trouble approximating the solution (\cref{comparison3d2l1}) accurately pointwise (\cref{comparison3d2l_one,vertical3d2l_one}). The 3-layer network structure with 4\% of the number of parameters of the 2-layer network approximated the solution accurately (\cref{3d test1 table}). As explained in Example \ref{2d test1 section}, \cref{comparison3d2l2,vertical3d2l,breaking3d2l} also indicate that the function $p$ in \cref{chi-curve lem} appears to be the approximation and in this example, be contained in $\mathcal{M}(3,10)$.

\begin{table}[htbp]\label{3d test1 table}
\caption{Relative errors of the problem in \Cref{3d test1 section}}
\centering
\begin{tabular}{|l|l|l|l|l|}
\hline
Network structure  &$\frac{\|u-{u}^{_N}_{_\cT}\|_0}{\|u\|_0}$ &$\frac{\vertiii{u-{u}^{_N}_{_\cT}}_{\bm\beta}}{\vertiii{u}_{\bm\beta}}$ & $\frac{\mathcal{L}^{1/2}({u}^{_N}_{_\cT},\bf f)}{\mathcal{L}^{1/2}({u}^{_N}_{_\cT},\bf 0)}$ & Parameters \\ \hline
3--300--1  & 0.185006 & 0.214390 & 0.189820   & 1501\\ \hline
3--5--5--1  & 0.055365 & 0.055370 & 0.045902   & 56\\ \hline
\end{tabular}
\end{table}

\subsubsection{A problem with a cylindrical interface}\label{3d test2 section}
\begin{figure}[htbp]\label{3d test2}
\centering
\subfigure[The interface\label{interface3d}]{
\begin{minipage}[t]{0.4\linewidth}
\centering
\includegraphics[width=1.8in]{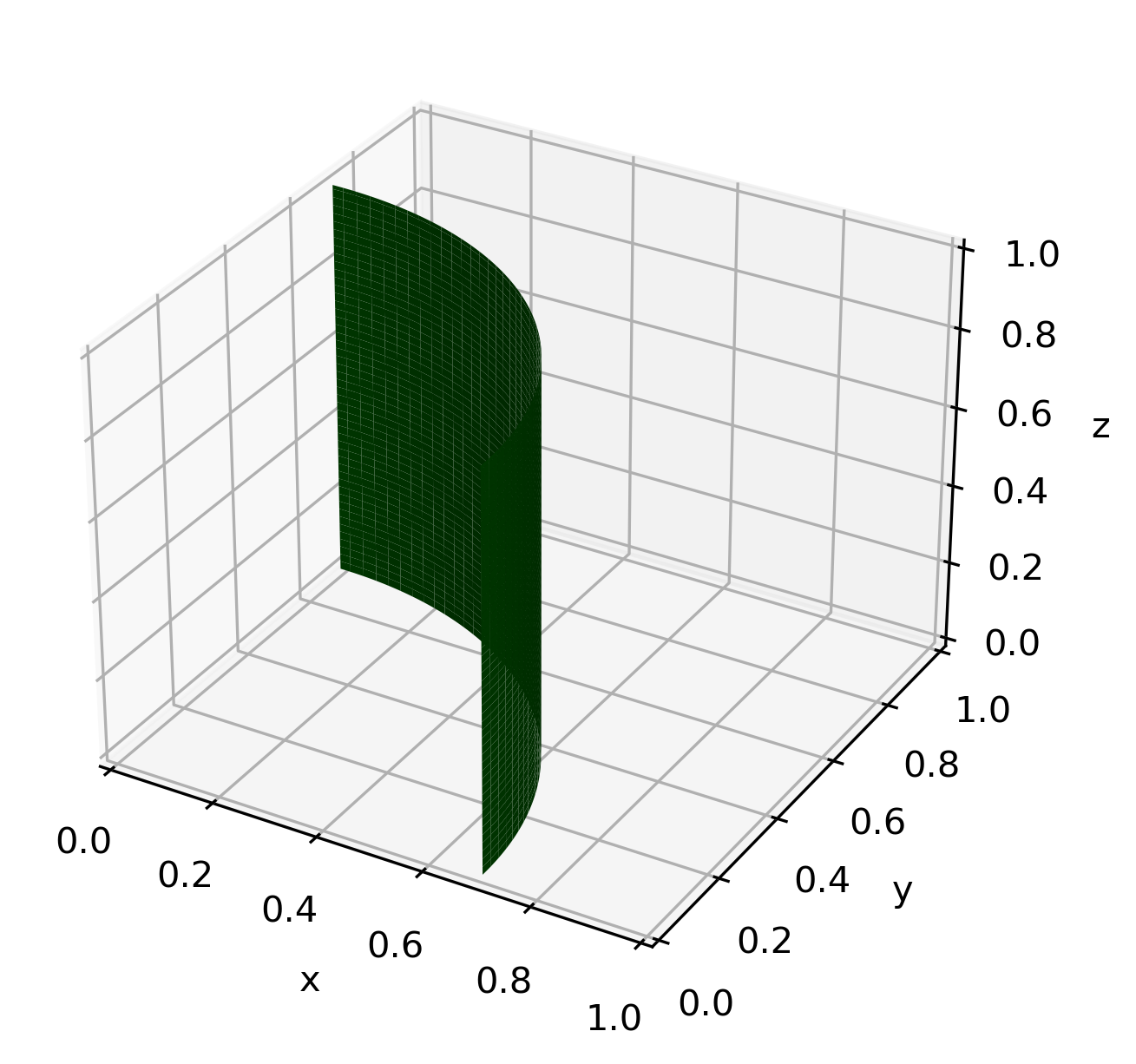}
\end{minipage}%
}%
\hspace{0.2in}
\subfigure[The exact solution\label{comparison3d1}]{
\begin{minipage}[t]{0.4\linewidth}
\centering
\includegraphics[width=1.8in]{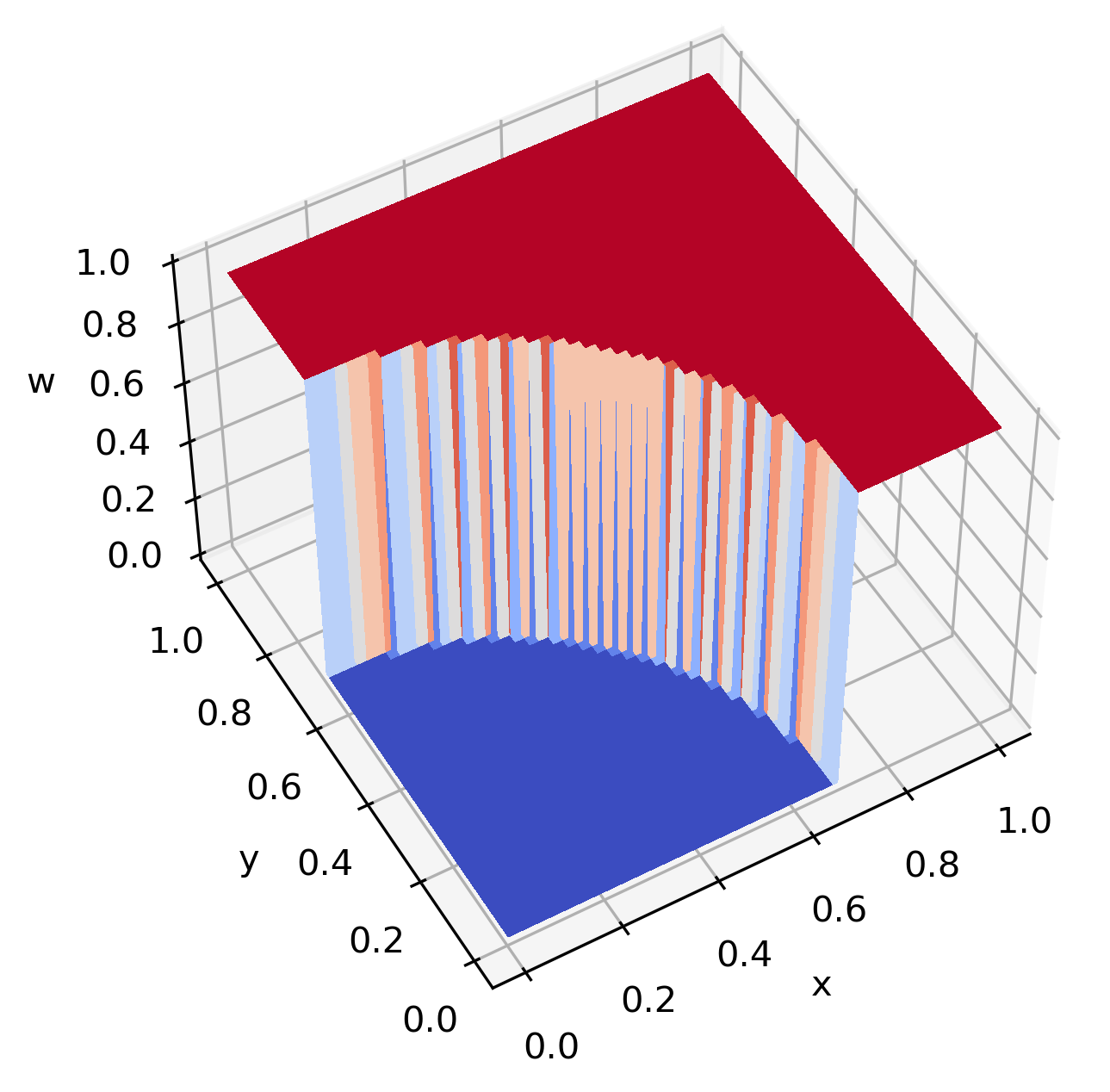}
\end{minipage}%
}%
\\
\subfigure[A 3--1500--1 ReLU NN function approximation\label{comparison3d_one}]{
\begin{minipage}[t]{0.4\linewidth}
\centering
\includegraphics[width=1.8in]{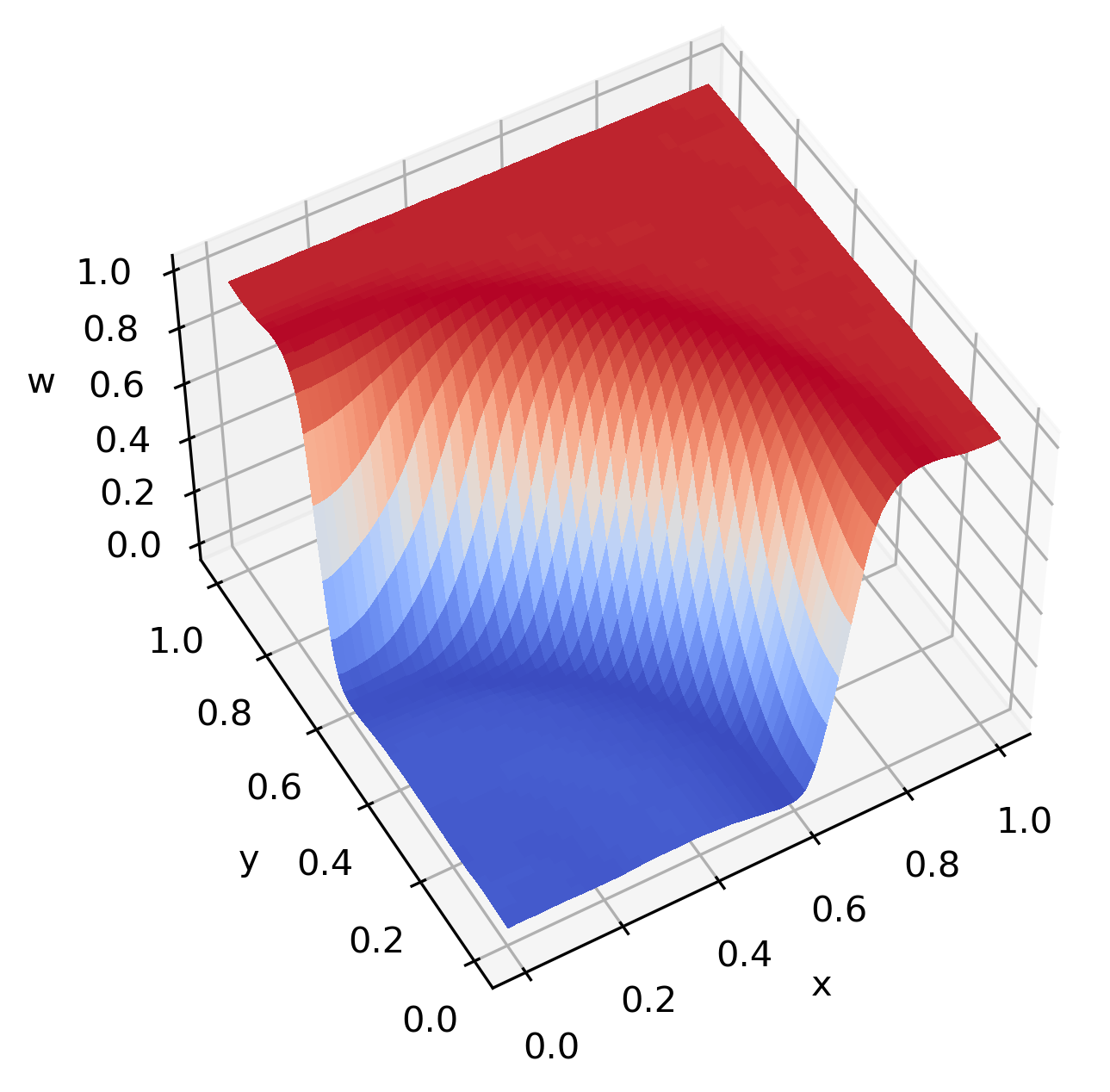}
\end{minipage}%
}%
\hspace{0.2in}
\subfigure[A 3--50--50--1 ReLU NN function approximation\label{comparison3d2}]{
\begin{minipage}[t]{0.4\linewidth}
\centering
\includegraphics[width=1.8in]{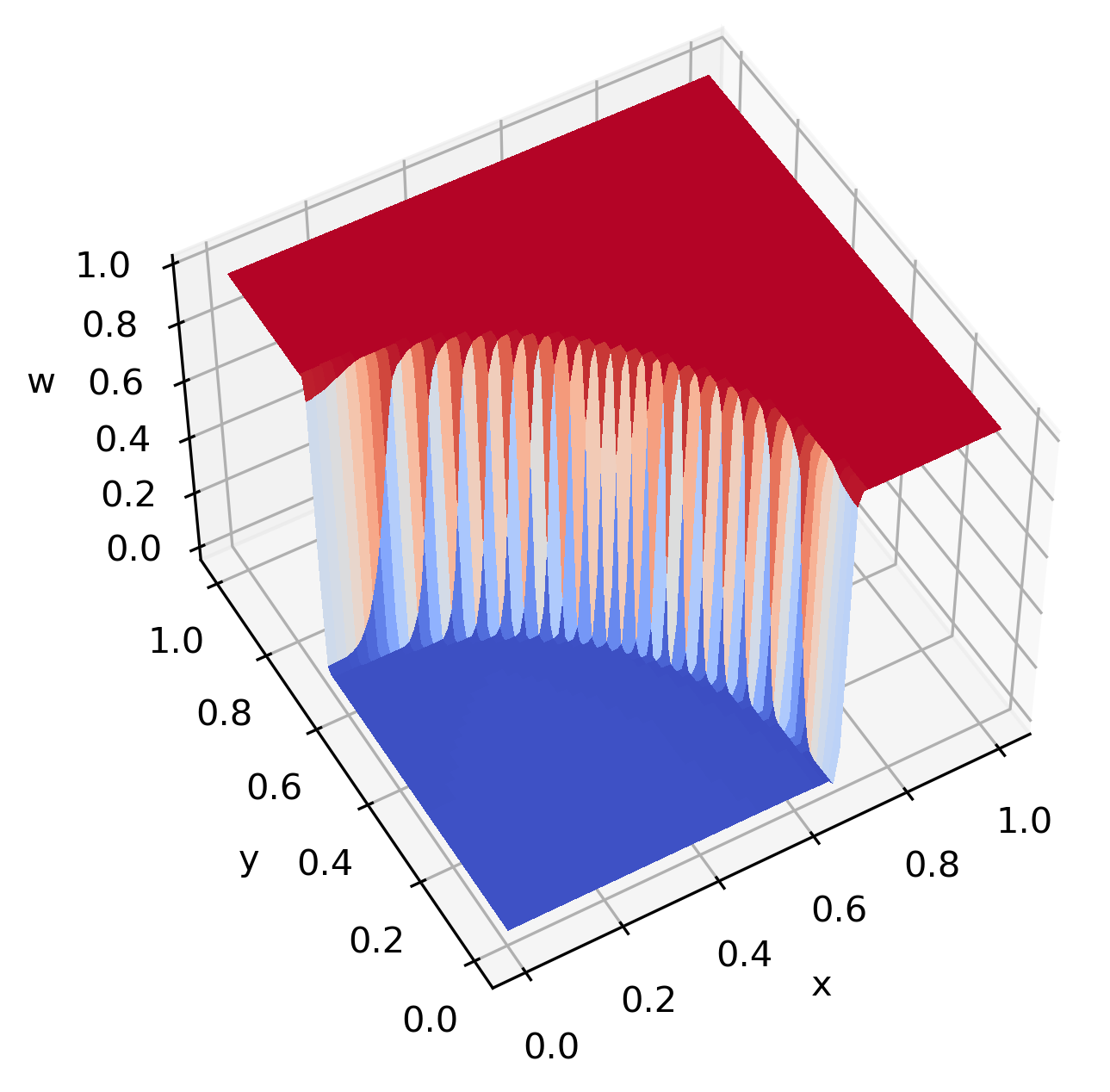}
\end{minipage}%
}%
\\
\subfigure[The trace of Figure \ref{comparison3d_one} on $y=x$\label{vertical3d_one}]{
\begin{minipage}[t]{0.4\linewidth}
\centering
\includegraphics[width=1.8in]{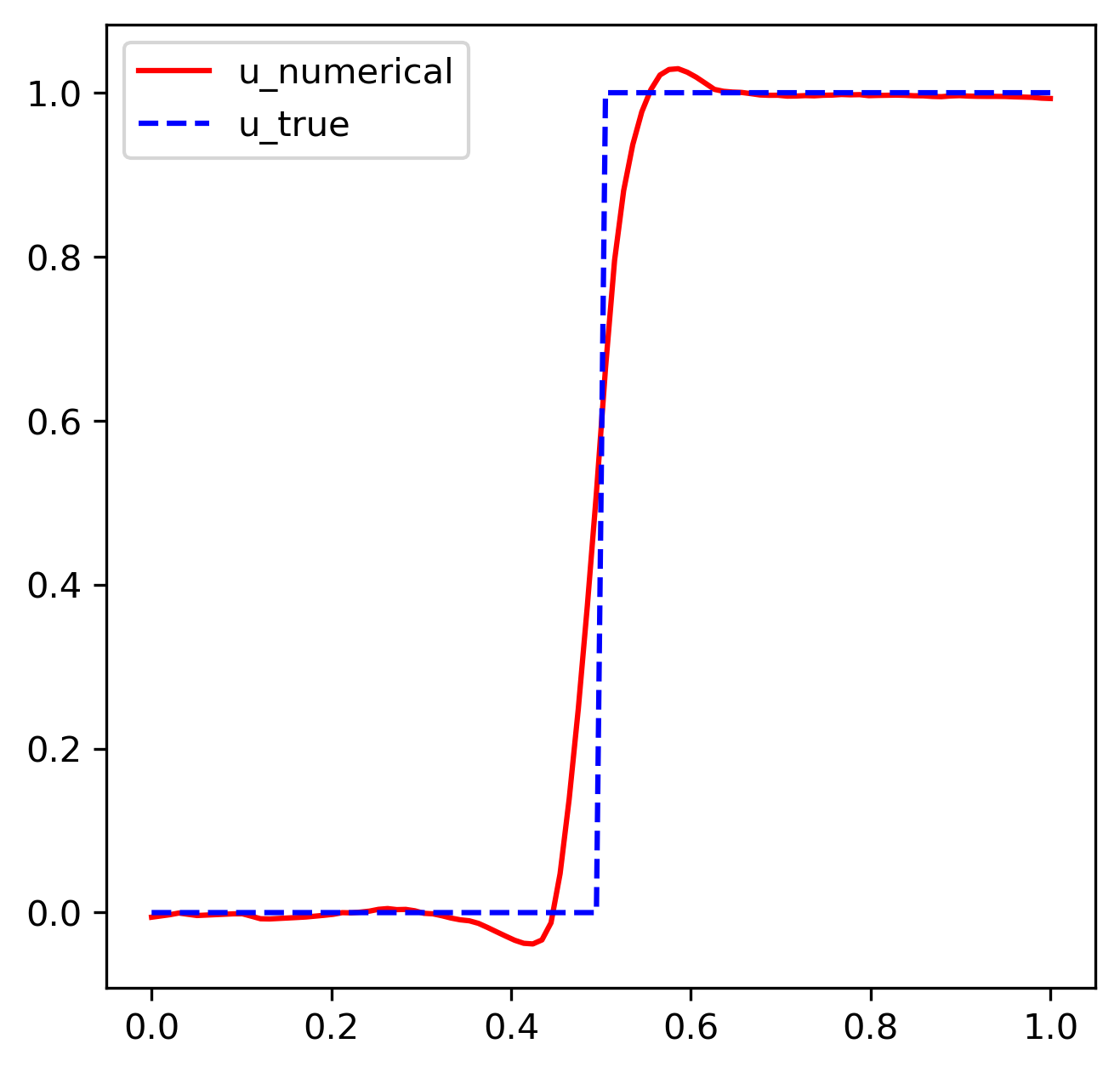}
\end{minipage}%
}%
\hspace{0.2in}
\subfigure[The trace of Figure \ref{comparison3d2} on $y=x$\label{vertical3d}]{
\begin{minipage}[t]{0.4\linewidth}
\centering
\includegraphics[width=1.8in]{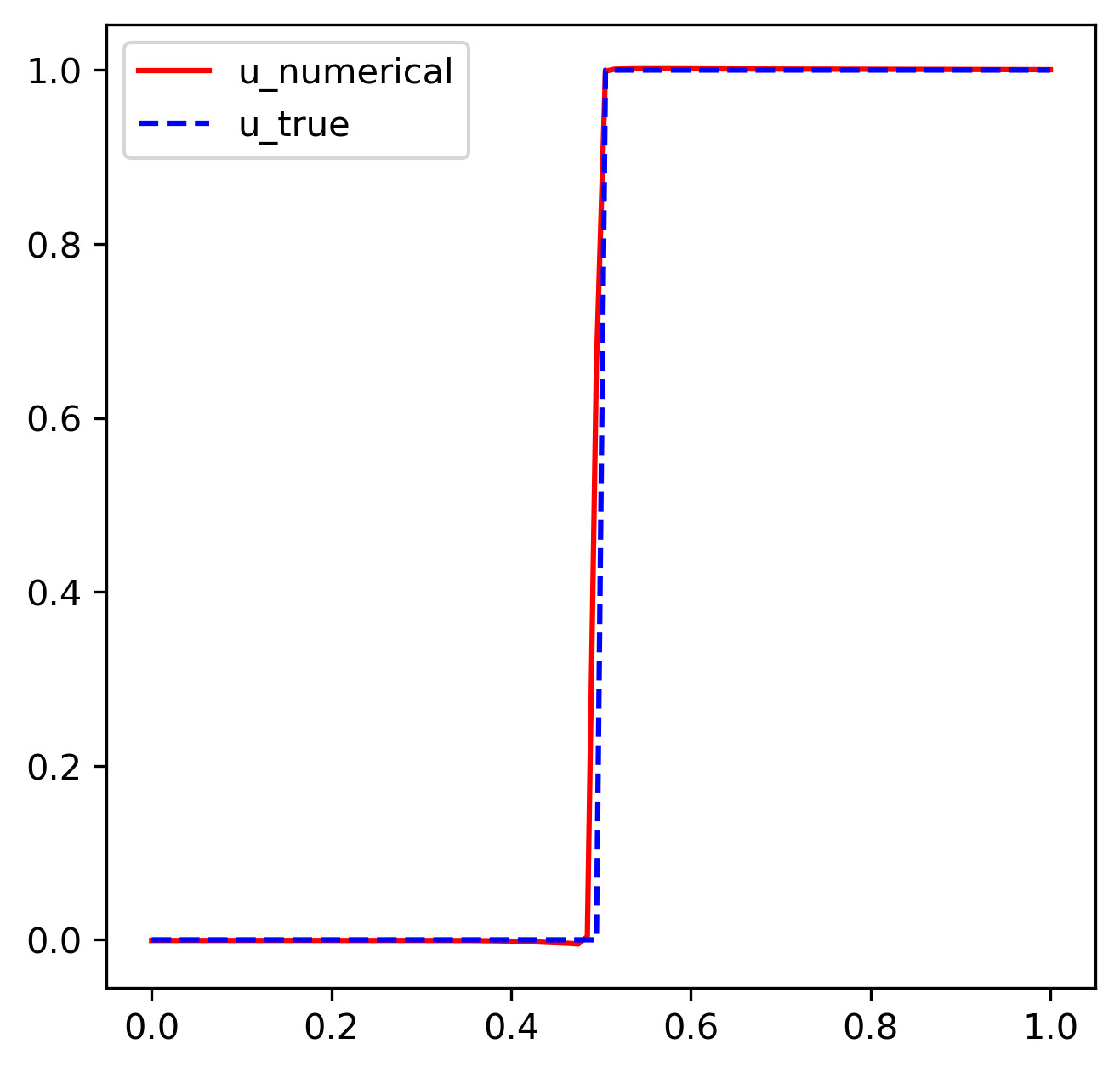}
\end{minipage}%
}%
\\
\subfigure[The breaking hyperplanes of the approximation in Figure \ref{comparison3d_one}\label{breaking3d_one}]{
\begin{minipage}[t]{0.4\linewidth}
\centering
\includegraphics[width=1.8in]{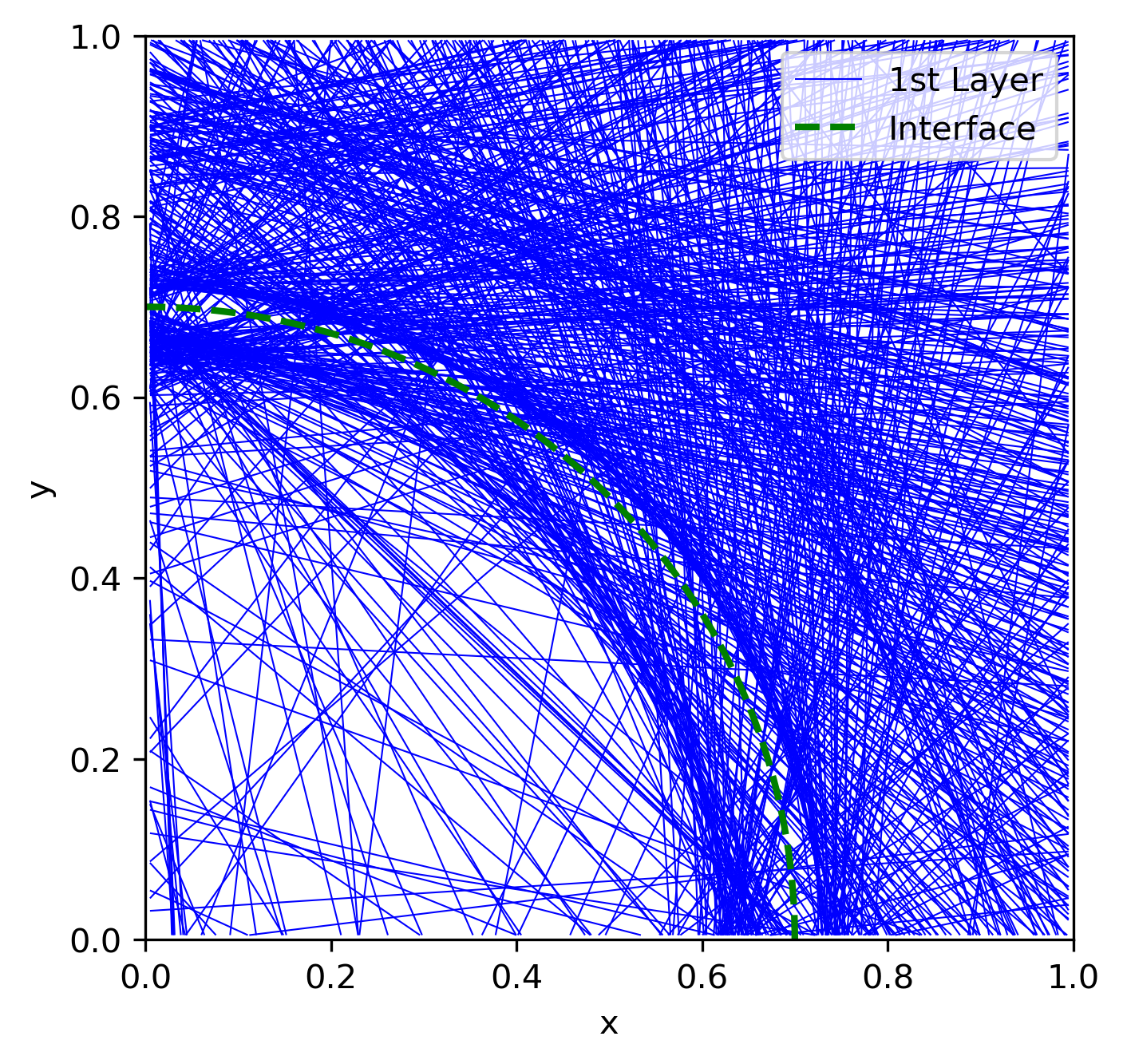}
\end{minipage}%
}%
\hspace{0.2in}
\subfigure[The breaking hyperplanes of the approximation in Figure \ref{comparison3d2}\label{breaking3d}]{
\begin{minipage}[t]{0.4\linewidth}
\centering
\includegraphics[width=1.8in]{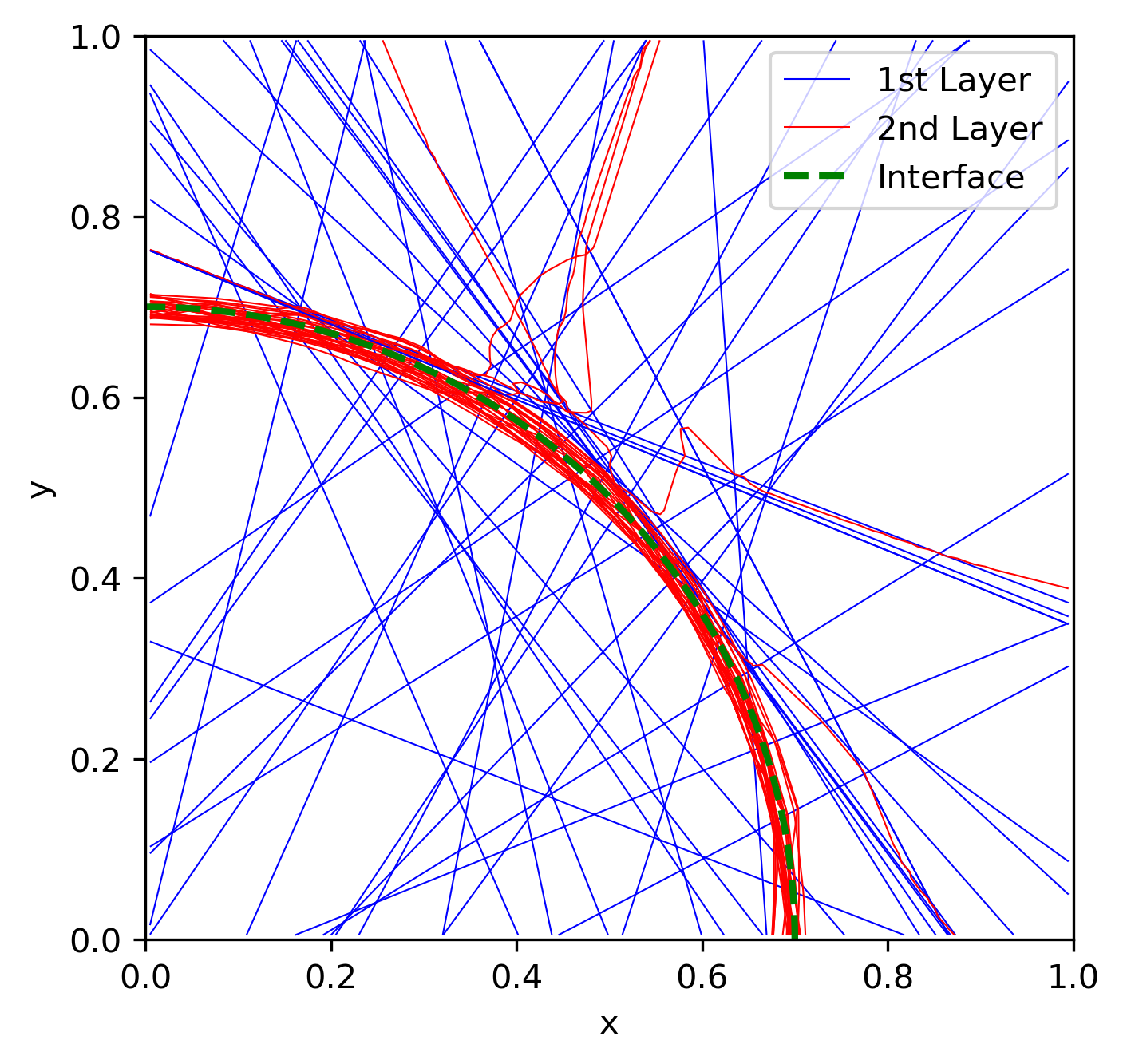}
\end{minipage}%
}%
\caption{Approximation results of the problem in \Cref{3d test2 section}}
\end{figure}

Let $\gamma=1$. The advective velocity field is a variable field given by
\begin{equation}
\bm{\beta}(x,y,z) =
(-y,x,0)^T,\,\, (x,y,z)\in\Omega.
\end{equation}
The inflow boundary and the inflow boundary condition are given by
\begin{eqnarray*}
\Gamma_{-}&=&\{(x,0,z):x,z\in(0,1)\}\cup\{(1,y,z):y,z\in(0,1)\}\\[2mm]
\mbox{and }\,\, g(x,y,z)&=&\left\{ \begin{array}{rl}
 0,& (x,y,z)\in \Gamma^1_-\equiv \{(x,0,z): x\in(0,0.7),\ z\in(0,1)\}, \\[2mm]
 1, &(x,y,z)\in \Gamma^2_-=\Gamma_-\setminus \Gamma_-^1,
\end{array}\right.
\end{eqnarray*} 
respectively. Let 
\[
\Omega_1=\{(x,y,z)\in\Omega:y<\sqrt{0.7^2-x^2}\}.
\]
The following right-hand side function is 
\begin{equation}
f(x,y,z)=\left\{ \begin{array}{rl}
0,& (x,y,z)\in \Omega_1, \\[2mm]
 1, & (x,y,z)\in \Omega_2=\Omega\setminus\Omega_1.
\end{array}\right.
\end{equation}
The exact solution is (see \cref{comparison3d1})
\[u(x,y,z)=f(x,y,z),\,\, (x,y,z)\in\Omega.\]
150000 iterations were implemented with 3--1500--1 and 3--50--50--1 ReLU NN functions. The numerical results are presented in \cref{3d test2,3d test2 table}. Again to minimize the loss function over a larger subset of CPWL functions, we increased the number of hidden neurons. Even though the 2-layer ReLU NN function approximation provides an approximate location of the discontinuity interface (\cref{interface3d,breaking3d_one}), the structure failed to approximate the solution (\cref{comparison3d1}) around the interface accurately (\cref{comparison3d_one,vertical3d_one}). The 3-layer network structure with less than 40\% of the number of parameters of the 2-layer network approximated the solution well, locating the discontinuity interface (\cref{comparison3d2,vertical3d,breaking3d,3d test2 table}).

\begin{table}[htbp]\label{3d test2 table}
\caption{Relative errors of the problem in \Cref{3d test2 section}}
\centering
\begin{tabular}{|l|l|l|l|l|}
\hline
Network structure  &$\frac{\|u-{u}^{_N}_{_\cT}\|_0}{\|u\|_0}$ &$\frac{\vertiii{u-{u}^{_N}_{_\cT}}_{\bm\beta}}{\vertiii{u}_{\bm\beta}}$ & $\frac{\mathcal{L}^{1/2}({u}^{_N}_{_\cT},\bf f)}{\mathcal{L}^{1/2}({u}^{_N}_{_\cT},\bf 0)}$ & Parameters \\ \hline
3--1500--1  & 0.125142 & 0.158393 & 0.117929   & 7501\\ \hline
3--50--50--1  & 0.050217 & 0.073780 & 0.018976   & 2801\\ \hline
\end{tabular}
\end{table}

\subsubsection{A problem with a spherical interface}\label{3d test3 section}

Let $\gamma=1$. The advective velocity field is a variable field given by
\begin{equation}\label{3d test3 beta}
\bm{\beta}(x,y,z) =
(-y-z,x,x)^T,\,\, (x,y,z)\in\Omega.
\end{equation}
The inflow boundary and the inflow boundary condition are given by
\begin{eqnarray*}
\Gamma_{-}&=&\{(x,0,z):x,z\in(0,1)\}\cup\{(1,y,z):y,z\in(0,1)\}\cup\{(x,y,0):x,y\in(0,1)\}\\[2mm]
\mbox{and }\,\, g(x,y,z)&=&\left\{ \begin{array}{rl}
 0,& (x,y,z)\in \Gamma^1_-\equiv \{(x,0,z): 0<z<\sqrt{0.7^2-x^2},\ x\in(0,0.7)\}, \\[2mm]
 0,& (x,y,z)\in \Gamma^2_-\equiv \{(x,y,0): 0<y<\sqrt{0.7^2-x^2},\ x\in(0,0.7)\}, \\[2mm]
 1, &(x,y,z)\in \Gamma^3_-=\Gamma_-\setminus (\Gamma_-^1\cup \Gamma_-^2),
\end{array}\right.
\end{eqnarray*} 
respectively.

\begin{figure}[htbp]\label{3d test3}
\centering
\subfigure[The interface\label{interface3dz}]{
\begin{minipage}[t]{0.4\linewidth}
\centering
\includegraphics[width=1.8in]{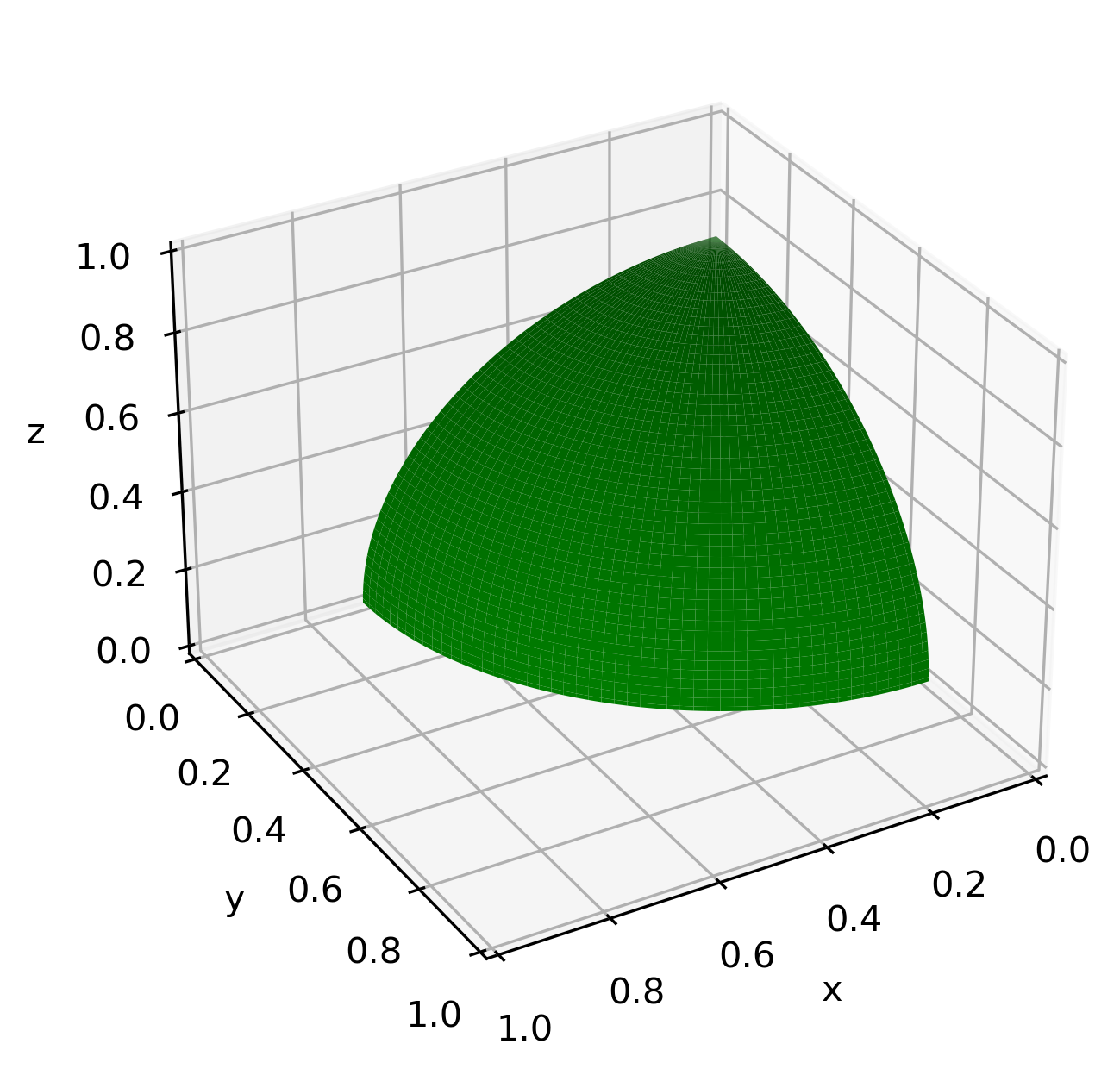}
\end{minipage}%
}%
\hspace{0.2in}
\subfigure[The exact solution\label{comparison3dz1}]{
\begin{minipage}[t]{0.4\linewidth}
\centering
\includegraphics[width=1.8in]{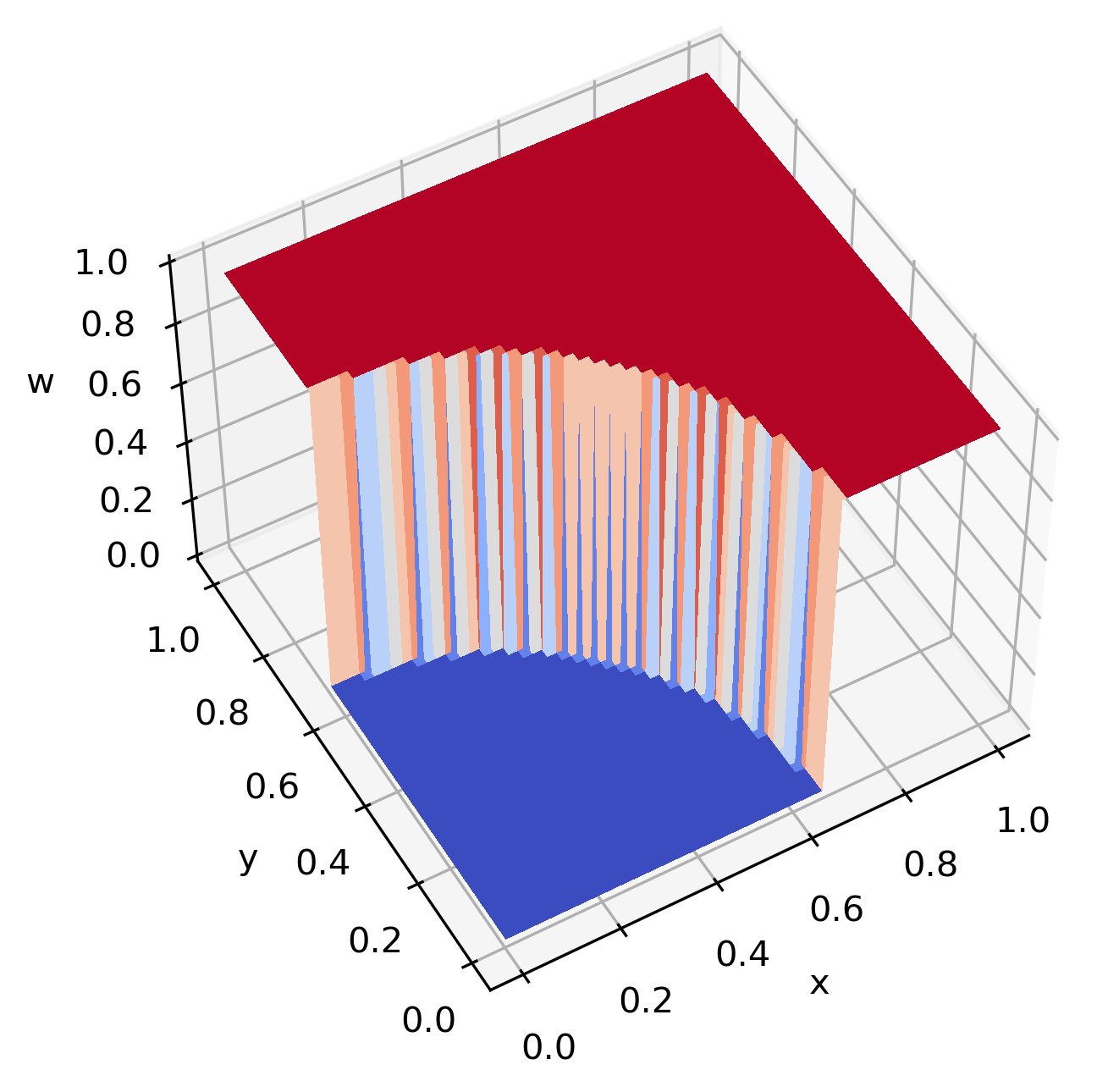}
\end{minipage}%
}%
\\
\subfigure[A 3--1376--1 ReLU NN function approximation\label{comparison3dz_one}]{
\begin{minipage}[t]{0.4\linewidth}
\centering
\includegraphics[width=1.8in]{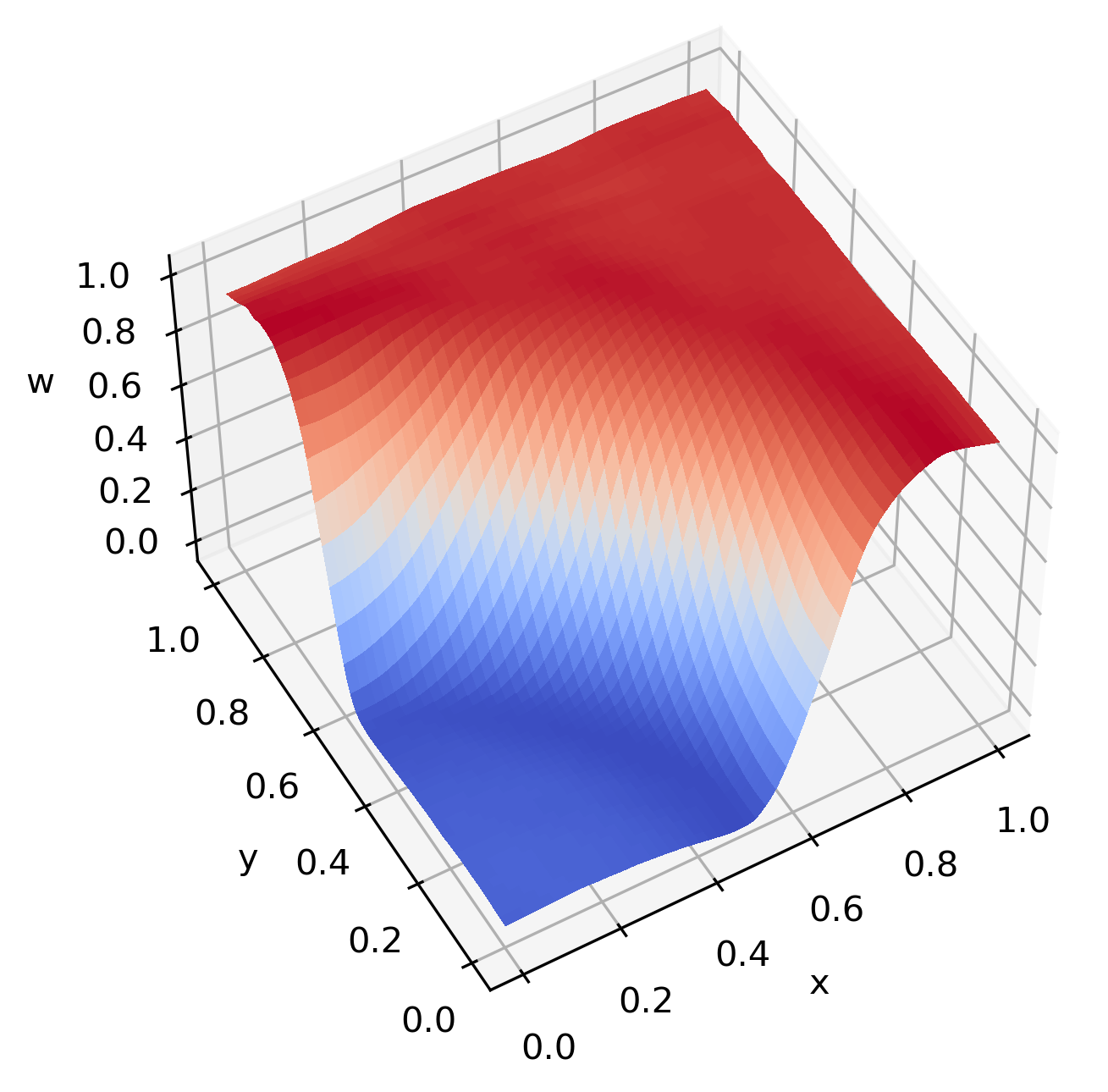}
\end{minipage}%
}%
\hspace{0.2in}
\subfigure[A 3--80--80--1 ReLU NN function approximation\label{comparison3dz2}]{
\begin{minipage}[t]{0.4\linewidth}
\centering
\includegraphics[width=1.8in]{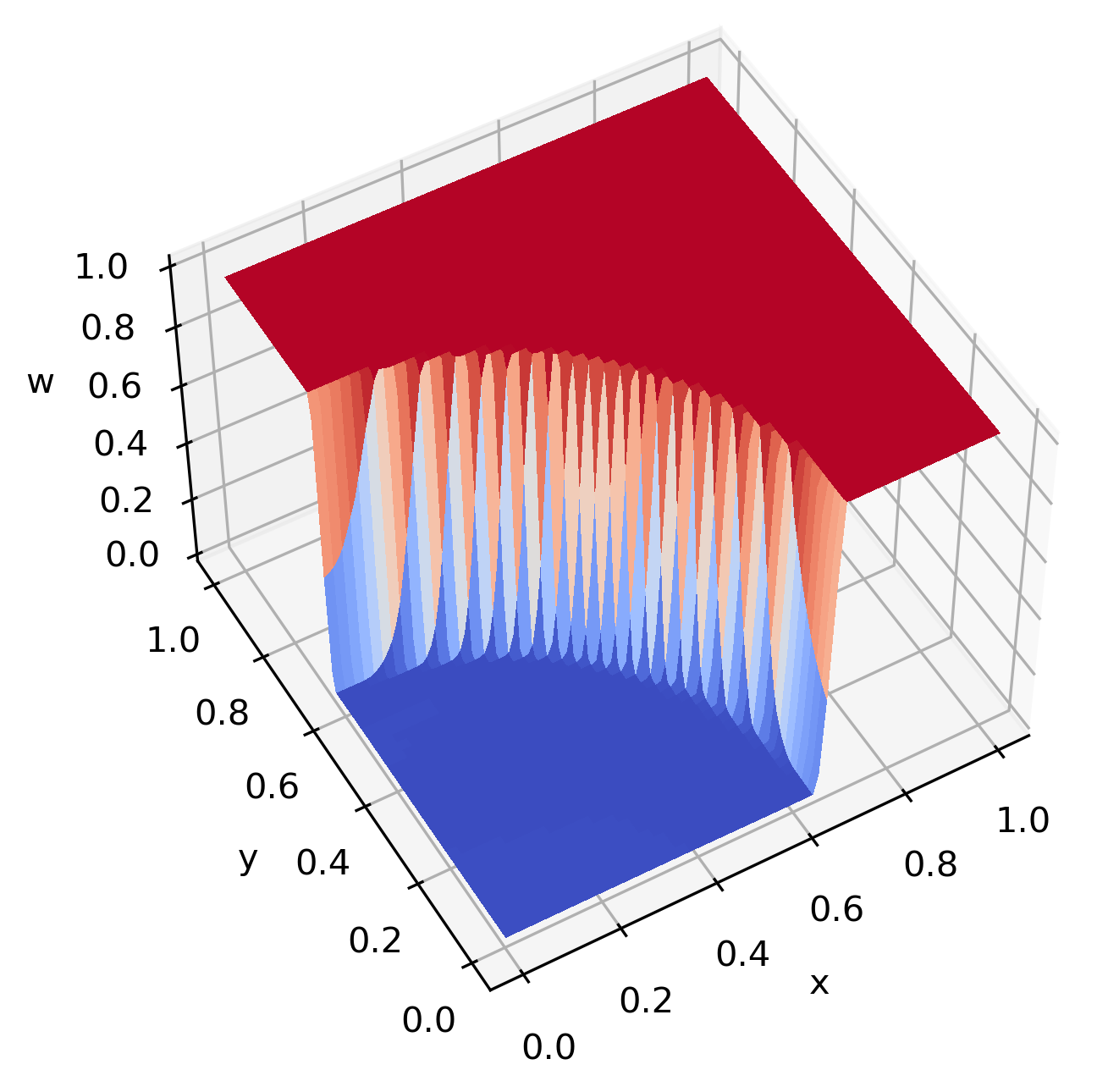}
\end{minipage}%
}%
\\
\subfigure[The trace of Figure \ref{comparison3dz_one} on $y=x$\label{vertical3dz_one}]{
\begin{minipage}[t]{0.4\linewidth}
\centering
\includegraphics[width=1.8in]{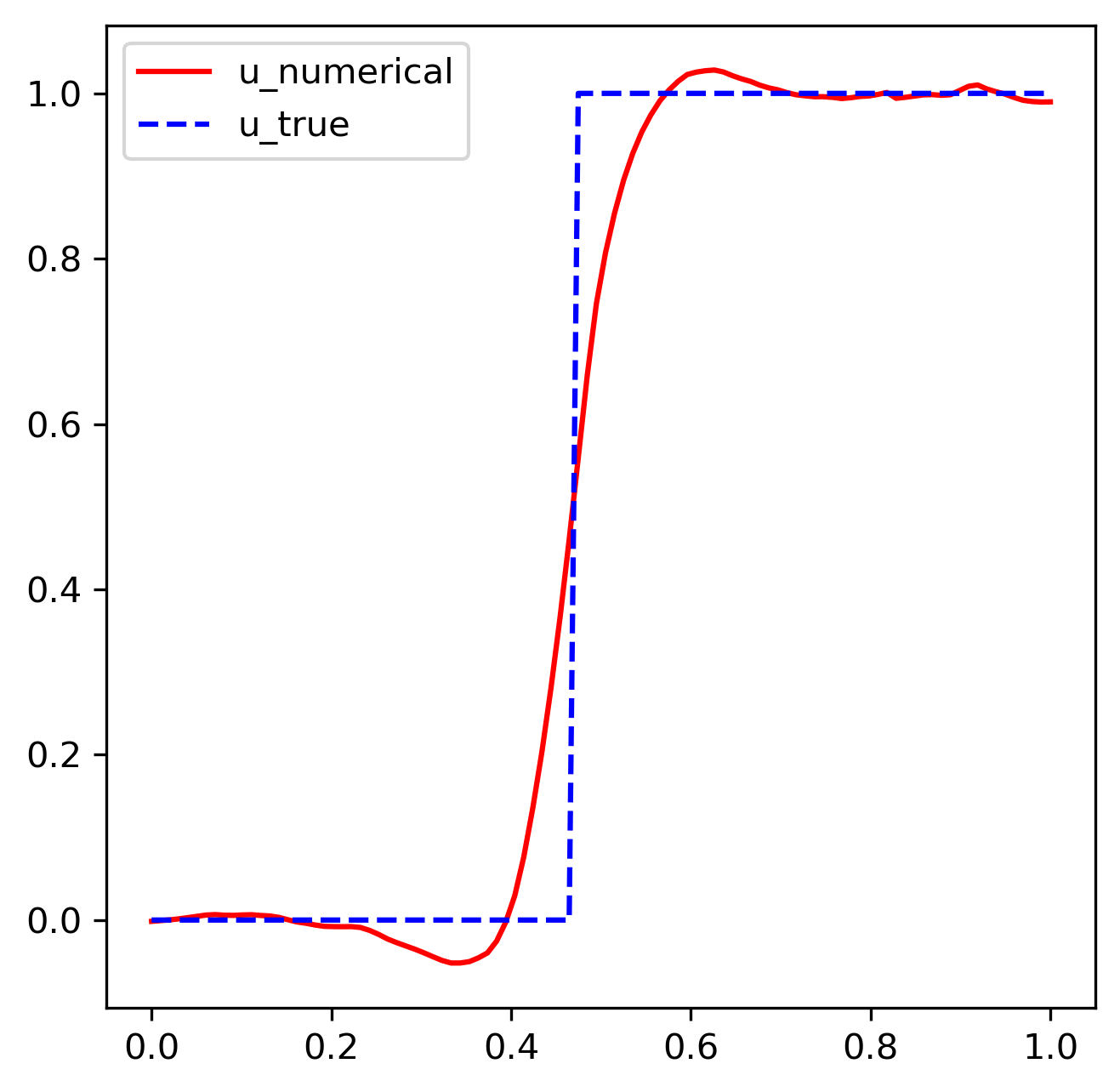}
\end{minipage}%
}%
\hspace{0.2in}
\subfigure[The trace of Figure \ref{comparison3dz2} on $y=x$\label{vertical3dz}]{
\begin{minipage}[t]{0.4\linewidth}
\centering
\includegraphics[width=1.8in]{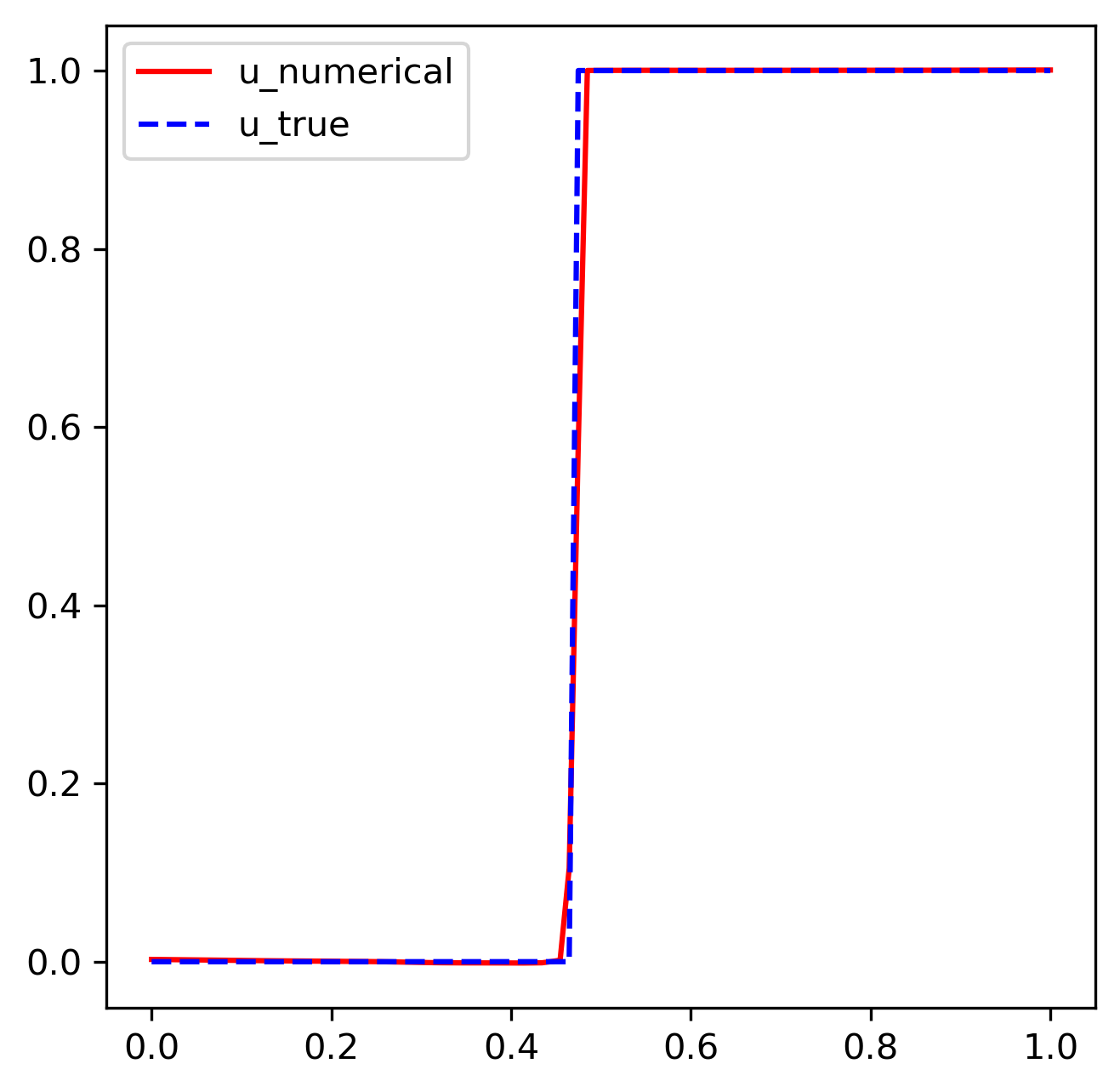}
\end{minipage}%
}%
\\
\subfigure[The breaking hyperplanes of the approximation in Figure \ref{comparison3dz_one}\label{breaking3dz_one}]{
\begin{minipage}[t]{0.4\linewidth}
\centering
\includegraphics[width=1.8in]{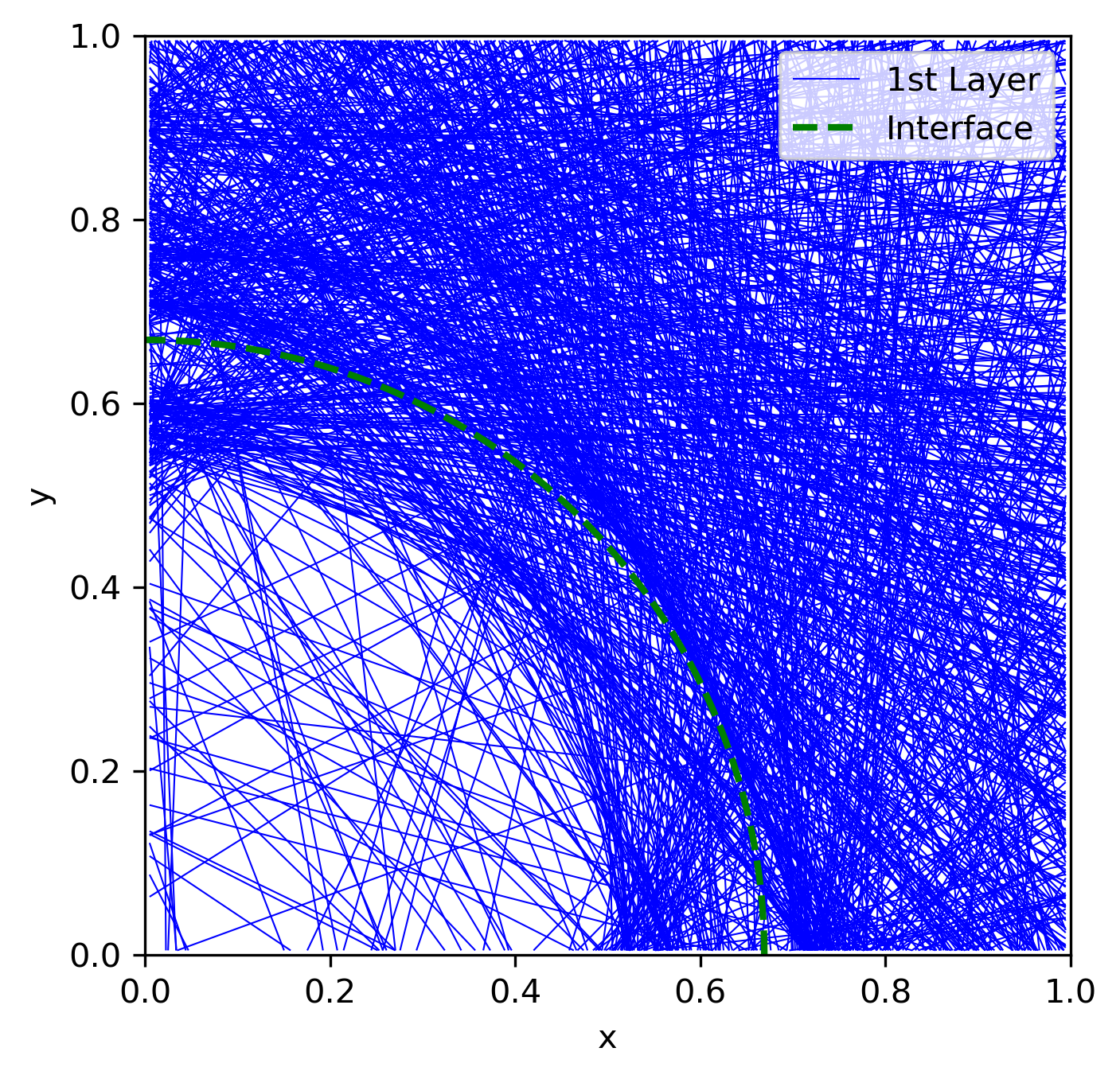}
\end{minipage}%
}%
\hspace{0.2in}
\subfigure[The breaking hyperplanes of the approximation in Figure \ref{comparison3dz2}\label{breaking3dz}]{
\begin{minipage}[t]{0.4\linewidth}
\centering
\includegraphics[width=1.8in]{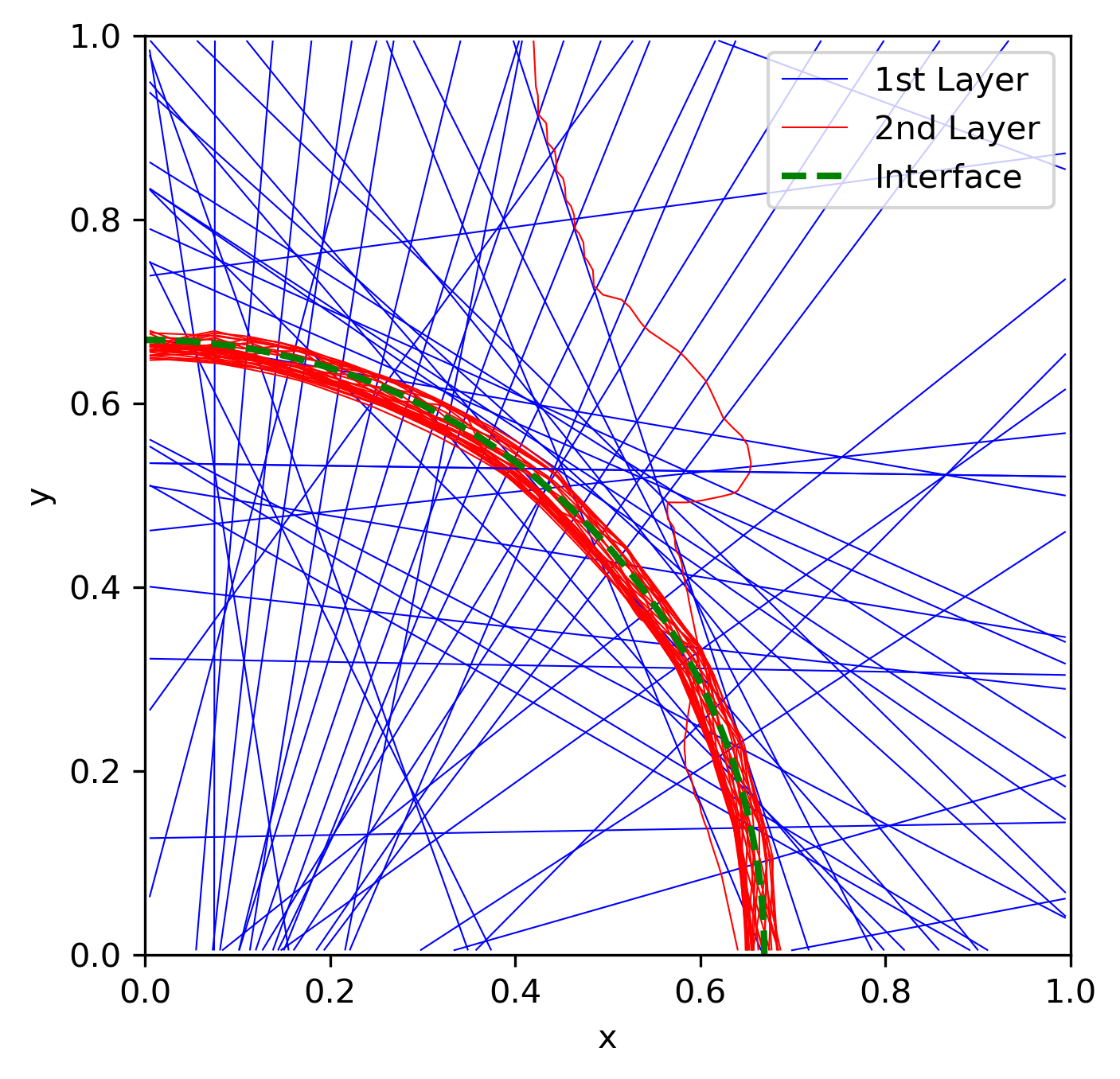}
\end{minipage}%
}%
\caption{Approximation results of the problem in \Cref{3d test3 section}}
\end{figure}

Let 
\[
\Omega_1=\{(x,y,z)\in\Omega:z<\sqrt{0.7^2-x^2-y^2}\}.
\]
The following right-hand side function is 
\begin{equation}
f(x,y,z)=\left\{ \begin{array}{rl}
0,& (x,y,z)\in \Omega_1, \\[2mm]
 1, & (x,y,z)\in \Omega_2=\Omega\setminus\Omega_1.
\end{array}\right.
\end{equation}
The exact solution is (see \cref{comparison3dz1})
\[u(x,y,z)=f(x,y,z),\,\, (x,y,z)\in\Omega.\]

200000 iterations were implemented with 3--1376--1 and 3--80--80--1 ReLU NN functions. The numerical results are presented in \cref{3d test3,3d test3 table}. We increased the number of hidden neurons and $\rho$ was set to $h/12$ in the finite difference quotient in \cref{finite_diff} because of the jump along the more curved interface (\cref{interface3dz}). Unlike the other two three-dimensional test problems, in this test problem, the third dimension actually plays a role as we can see from the advective velocity field \cref{3d test3 beta}. The approximation results are similar to those of the previous examples. Moreover, we note that the two-hidden-layer NN outperforms the one-hidden-layer NN with the same number of parameters, which together with the previous examples, suggests that the one-hidden-layer NN needs more neurons.

\begin{table}[htbp]\label{3d test3 table}
\caption{Relative errors of the problem in \Cref{3d test3 section}}
\centering
\begin{tabular}{|l|l|l|l|l|}
\hline
Network structure  &$\frac{\|u-{u}^{_N}_{_\cT}\|_0}{\|u\|_0}$ &$\frac{\vertiii{u-{u}^{_N}_{_\cT}}_{\bm\beta}}{\vertiii{u}_{\bm\beta}}$ & $\frac{\mathcal{L}^{1/2}({u}^{_N}_{_\cT},\bf f)}{\mathcal{L}^{1/2}({u}^{_N}_{_\cT},\bf 0)}$ & Parameters \\ \hline
3--1376--1  & 0.113045 & 0.144105 & 0.106094   & 6881\\ \hline
3--80--80--1  & 0.042233 & 0.064332 & 0.041935   & 6881\\ \hline
\end{tabular}
\end{table}


\bibliographystyle{siamplain}
\bibliography{references}
\end{document}